\documentclass[11pt]{article}
\usepackage{graphicx}
\usepackage{latexsym}

\usepackage{amsmath}
\usepackage{amsthm}
\usepackage{amssymb}
\usepackage{natbib}

\newtheorem{theorem}{Theorem}[section]
\newtheorem{lemma}[theorem]{Lemma}
\newtheorem{corollary}[theorem]{Corollary}
\newtheorem{proposition}[theorem]{Proposition}

\newtheorem{conjecture}[theorem]{Conjecture}
\newcommand{\goto}{\rightarrow}

\newcommand{\bitem}{\begin{itemize}}
\newcommand{\eitem}{\end{itemize}}

\newcommand{\argmin}{\mbox{arg min}}
\newcommand{\hf}{\textstyle{\frac{1}{2}}}

\newenvironment{weird}{}{}
\newcommand{\bn}{\begin{weird}}
\newcommand{\en}{\end{weird}}

\begin{document}

\title{Adapting to Unknown Sparsity \\
by controlling the \\
False Discovery Rate}
\author{Felix Abramovich, Yoav Benjamini, David L. Donoho and Iain M. Johnstone}
\date{VERSION OF: April 19 2005}
\maketitle

\begin{abstract}
We attempt to recover an $n$-dimensional
vector observed in white noise, where 
$n$ is large and the vector is 
known to be
sparse, but the degree of sparsity is unknown. We consider
three different ways of defining sparsity of a vector:  using the
fraction of nonzero terms; imposing power-law decay bounds on the
ordered entries; and controlling the $\ell_p$ norm for $p$ small.
We obtain a procedure which is asymptotically minimax for
$\ell^r$ loss, simultaneously throughout a range
of such sparsity classes.

The optimal procedure is
a data-adaptive thresholding
scheme, driven by control of the {\it False Discovery Rate} (FDR).
FDR control is a relatively recent innovation in simultaneous testing,
ensuring that at most a certain fraction
of the rejected null hypotheses will correspond to false
rejections.

In our treatment, the FDR control parameter $q_n$
also plays a determining role in asymptotic minimaxity.
If $q  = \lim q_n \in [0,1/2]$ and also
$q_n  > \gamma/\log(n)$ we get sharp asymptotic
minimaxity, simultaneously, over a wide range of
sparse parameter spaces
and loss functions.  On the other hand, 
$ q 
= \lim q_n \in (1/2,1]$, forces
the risk to exceed the minimax risk
by a factor growing with $q$.

To our knowledge, this relation between ideas in simultaneous
inference and asymptotic decision theory is new.

Our work provides a new perspective
on a class of model selection rules
which has been introduced
recently by several
authors. These new
rules impose complexity penalization
of the form $2 \cdot \log( \mbox{ potential model
size } / \mbox{  actual model size } )$. We exhibit a close connection with
FDR-controlling procedures under stringent control
of the false discovery rate.

\end{abstract}


\medskip

\textbf{Short Title}  {\sc False Discovery Rate Thresholding}

\normalsize
\medskip
\textbf{Acknowledgements:} 
ABDJ would like to acknowledge
the 
support of Israel-USA BSF grant 1999441. IMJ
would like to acknowledge support
of NSF grants {\sc dms} 95--05151, 00-77621,
an NIH grant, a Guggenheim Foundation Fellowship,
and the Australian National University. DLD
would like to thank FA and YB for hospitality
at Tel Aviv University during a sabbatical there,
and would like to acknowledge the support of
AFOSR MURI 95--P49620--96--1--0028.

\pagebreak

\section{Introduction}
The problem of {\it model selection} has
attracted the attention of both applied and theoretical
statistics for as long as anyone can remember.
In the setting of the standard linear model,
we have noisy data on a response variable which we wish to predict
linearly using a subset of a large collection of predictor variables.
We believe that good parsimonious models can be constructed
using only a relatively few variables from the available ones.
In the spirit of the modern, computer-driven
era, we would like a simple
automatic procedure which is data adaptive, can find a good parsimonious
model when one exists, and is effective for very
different types of data and model.

There has been an enormous range of contributions to
this problem, so large in fact that it would be
impractical to summarize here.  Some key
contributions, mentioned further below,
include the AIC, BIC, and RIC model selection
proposals \citep{Akaike,Mallows,Schwarz,RIC}.
Key insights from this vast literature are
\begin{itemize}
\item  The tendency of certain rules (notably AIC), when
used in an exhaustive model search mode, to include
too many irrelevant predictors \--- \cite{BreimanFreedman};
\item  The tendency of rules which do not suffer from this problem
(notably RIC) to place evidentiary standards for inclusion
in the model that are far stricter than the time-honored
`individually significant' single coefficient approaches.
\end{itemize}

In this paper we consider a very special
case of the model selection problem
in which a full decision-theoretic
analysis of predictive risk can be carried out.
In this setting, model parsimony can be concretely defined
and utilized, and we exhibit a model selection method
enjoying optimality over a wide
range of parsimony classes.
While the full story is rather technical,
at the heart of the method
is a simple practical method with an
easily understandable benefit: the ability to
prevent the inclusion of too many irrelevant
predictors -- thus improving on AIC -- while setting
lower standards for inclusion -- thus improving on RIC.
The optimality result assures us that in
a certain sense the method is unimprovable.

Our special case is the problem
of estimating  a high-dimensional
mean vector which is sparse,
when the nature and degree of sparsity
are unknown and may vary through a range
of possibilities. We consider three
ways of defining sparsity and
will derive asymptotically
minimax procedures applicable across
all modes of definition.

Our asymptotically minimax procedures
will be based on a relatively recent innovation --
False Discovery Rate (FDR) control in multiple
hypothesis testing. The FDR control parameter plays
a key role in delineating superficially similar cases where
one can achieve asymptotic minimaxity
and where one cannot.

To our knowledge, this connection between developments in these two important
subfields of statistics is new.  Historically,
the multiple hypothesis testing literature has had little to do with
notions like minimax estimation or asymptotic minimaxity in estimation.

The procedures we propose will be very easy to
implement and run quickly on computers.
This is in sharp contrast to certain
optimality results in minimaxity which exhibit optimal procedures
that are computationally unrealistic.  Finally, because of
recent developments in harmonic analysis
-- wavelets, wavelet packets, etc. -- these results are of
immediate practical significance in applied settings.
Indeed, wavelet analysis of noisy signals can result
in exactly the kind of sparse means problem discussed here.

Our goal in this introduction is to make clear to the non
decision-theorist the motivation for these results, the form of a few
select results, and some of the implications.  Later sections
will give full details of the proofs and the methodology
being studied here.

\subsection{Thresholding}

\label{sec:intro}
Consider the standard multivariate normal mean problem:
\begin{equation}
        y_i = \mu_i + \sigma_n z_i,\qquad z_i \stackrel{i.i.d.}{\sim}
        N(0,1), \qquad i = 1, \ldots , n.
        \label{eq:gaussshift}
\end{equation}
Here $\sigma_n$ is known, and the goal is to estimate the unknown vector $\mu$
lying in a fixed set $\Theta_n$.
The index $n$ counts the number of variables and is assumed large. The key extra
assumption, to be quantified later, is that the vector $\mu$ is
\textit{sparse}: only a small
number of components are significantly large, and the indices, or locations of
these large components are not known in advance. In such situations,
thresholding will be appropriate;  specifically,
hard thresholding at threshold $t \sigma_n$, meaning the estimate
$\hat{\mu}$ whose $i^{\mbox{th}}$ component is
\begin{equation}
        \hat{\mu}_i = \eta_H(y_i,t) =
        \begin{cases}
            y_i & |y_i| \geq t \sigma_n  \\
            0 & \text{else}.
        \end{cases}
        \label{eq:hardthresh}
\end{equation}

A compelling motivation for this strategy is
provided by wavelet analysis, since the wavelet representation of many
smooth and piecewise smooth signals is sparse in precisely our sense
\citep{djkp93}.  Consider, for example, the empirical wavelet
coefficients in Figure 1(c).  Model \eqref{eq:gaussshift} is quite
plausible if we consider the coefficients to be grouped \textit{level
by level}. Within a level, the number of large coefficients is small,
though the relative number clearly decreases as one moves from coarse
to fine levels of resolution.

\begin{figure}[htb]
      \begin{center}
        \leavevmode
     \centerline{\includegraphics[width= 6cm, height = 9cm,
     angle = 90]{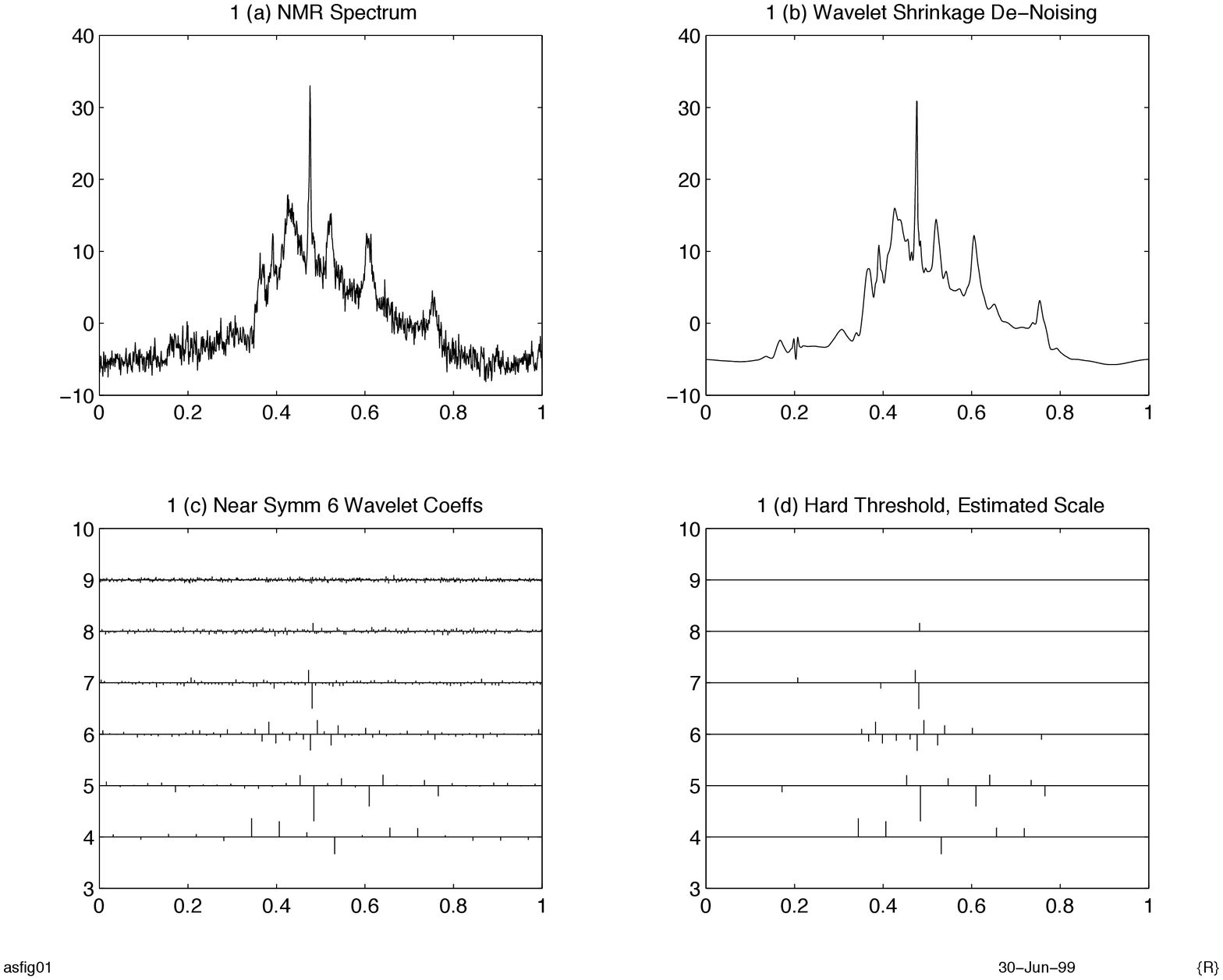}}
\caption{ \small   (a): sample NMR spectrum provided by A. Maudsley and
      C. Raphael, $n=1024$, and discussed in \citet{djkp93}. (c):
      Empirical wavelet coefficients $w_{jk}$ displayed by nominal
      location and scale $j$, computed using a discrete orthogonal wavelet
      transform and the Daubechies near symmetric filter of order $N=6$.
      (d): Wavelet coefficients after hard thresholding using the FDR
      threshold described at \eqref{eq:fdr-def}, with estimated scale
      $\hat{\sigma} = med.abs.dev. (w_{9k})/.6745,$ a resistant estimate
      of scale at level 9 -- for details on $\hat \sigma$, see
      \citet{djkp93}.  (b): Reconstruction using inverse discrete wavelet
      transform. }
        \label{fig:nmr-by-fdr}
      \end{center}
\end{figure}

\subsection{Sparsity}

In certain subfields of signal and image processing,
the wavelet coefficients of a typical
object can be modeled as a sparse vector; 
the interested reader might  consult literature going back to
\citet{field}, extending through \citet{djl92}, \citet{ruderman},
\citet{simoncelli} and \citet{huangmumford}.  A representative result
was given by Simoncelli, who found that in looking at a database
of images, the typical behavior of histograms of wavelet coefficients
at a single resolution level of the wavelet pyramid was highly
structured, with a sharp peak at the origin and somewhat heavy tails.
In short, many coefficients are small in amplitude while
a few are very large.

Wavelet analysis of images is not the only
place where one meets transforms with sparse
cofficients. There are several other signal
processing settings -- for example acoustic signal processing -- where,
when viewed in an appropriate basis, the underlying
object has sparse coefficients \citep{benedetto}.

In this paper we consider several ways to define
sparsity precisely.

The most intuitive notion of sparsity is
simply that there is a {\it relatively small proportion
of nonzero coefficients}.  Define the $\ell_0$ quasi-norm
by $\| x \|_0 = \#\{i : x_i \neq 0 \}$. Fixing a proportion
$\eta$, the collection of sequences with at
most a proportion $\eta$ of nonzero entries is
\begin{equation}
          \ell_0[\eta] = \{ \mu \in \mathbb{R}^n: \| \mu \|_0 \leq \eta n \}
       \label{eq:l0balldef}
\end{equation}
By analogy with night-sky images,
we will call \textit{nearly-black}
a setting where the
fraction of non-zero entries 
$\eta \approx 0$ \citep{djhs92}.

Sparsity can also mean that there is a relatively small
proportion of relatively large entries. Define
the decreasing rearrangement
of the amplitudes of the entries so that
\[
       |\mu|_{(1)} \geq |\mu|_{(2)} \geq ... \geq |\mu|_{(n)};
\]
we control the entries by a termwise power-law bound
on the decreasing rearrangements:
\[
        |\mu|_{(k)} \leq C \cdot k^{-\beta}, k =1,2,\dots .
\]
For reasons which will not be
immediately obvious, we work with $p = 1/\beta$ instead,
and call such a constraint a weak-$\ell_p$ constraint.
The interesting range is $p$ small, yielding substantial sparsity.
One can check whether a vector obeys such a constraint
by plotting the decreasing rearrangement on semilog
axes, and comparing the plot with a straight line of
slope $-1/p$. Certain values
of $p < 2$ provide a reasonable model for wavelet
coefficients of real-world images; \citet{djl92}.

Formally, a \textit{weak $\ell_p$ ball} of radius $\eta$
is defined by requiring that the
ordered magnitudes of components of
$\mu$ decay quickly:
\begin{equation}
          m_p[\eta] = \{ \mu \in \mathbb{R}^n:  |\mu|_{(k)} \leq \eta
          n^{1/p} k^{-1/p} \ \text{for all}  \ k = 1, \ldots , n \}.
       \label{eq:mpballdef}
\end{equation}
Weak $\ell_p$ has a natural `least-sparse'
sequence, namely
\begin{equation}
      \label{eq:extremal}
      \bar \mu_k = \eta n^{1/p} k^{-1/p}, \qquad k = 1, \ldots , n
\end{equation}
(and its permutations).
We also measure sparsity using $\ell_p$ norms with $p$ small:
\begin{equation} \label{eq:lpballdef}
      \| \mu \|_p = ( \sum_{i=1}^n | \mu_i |^p )^{1/p}.
\end{equation}
That small $p$ emphasises sparsity may be seen by noting that the
two vectors
\begin{displaymath}
      (1, 0, \ldots , 0) \qquad \text{and} \qquad ( n^{-1/p}, \ldots , n^{-1/p})
\end{displaymath}
have equivalent $\ell_p$ norms, but when $p$ is small, the components
of the latter, dense, vector are all negligible.
Strong-$\ell_p$ balls of small average
radius $\eta$ are defined so:
\begin{displaymath}
      \ell_p[\eta] = \{ \mu \in \mathbb{R}^n: \tfrac{1}{n} \sum_{i=1}^n
      | \mu_i |^p \leq \eta^p \}.
\end{displaymath}
If we refer to $\ell_p$ without qualification  -- weak
or strong -- we mean strong $\ell_p$.

There are important relationships between these classes.
Note that as $p \rightarrow 0$, the $\ell_p$ norms approach $\ell_0$:
$\| \mu \|_p^p \rightarrow \| \mu \|_0 $.
Weak-$\ell_p$ balls contain the corresponding strong $\ell_p$ balls,
but only just:
\begin{displaymath}
      \ell_p[\eta] \subset m_p[\eta] \not\subset \ell_{p'}[\eta], \qquad p'
      > p.
\end{displaymath}

\subsection{Adapting to Unknown Sparsity}

Estimation of sparse normal means over $\ell_p$ balls
has been carefully studied in \citet{dojo94}, with the result
that much is known about asymptotically minimax strategies
for estimation. In essence, {\it if we know} the degree of
sparsity of the sequence, then it turns out that
thresholding is indeed asymptotically minimax, and there are
simple formulas for optimal thresholds.

\begin{figure}[htb]
      \begin{center}
        \leavevmode
     \centerline{\includegraphics[width= 10cm, height = 8cm,
     angle = 0]{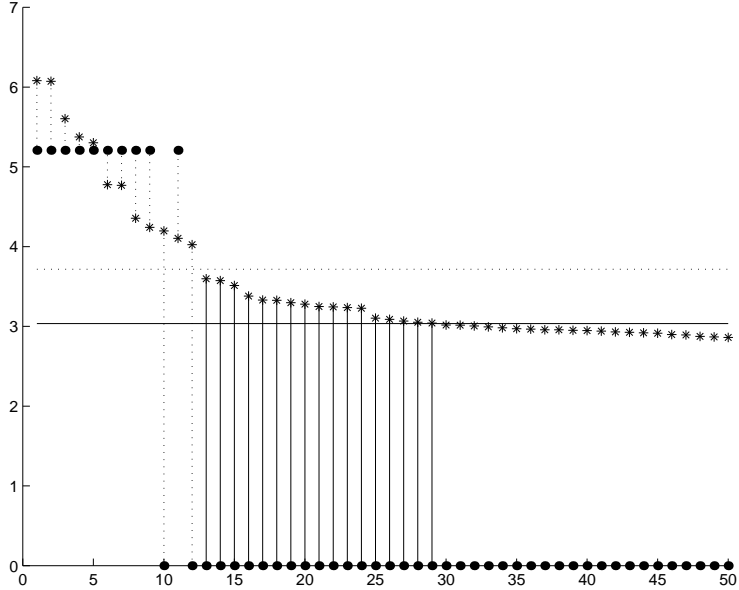}}
\caption{ \small Gaussian shift model \eqref{eq:gaussshift} with
    $n=10,000$ and $\sigma_n = 1$. There are $n_0 = n^{1/4} = 10$
    non-zero components $\mu_i = \mu_0 = 5.21.$ Thus $\beta=1/4.$ Stars
    show ordered data $|y|_{(k)}$ and solid circles the corresponding
    true means.  Dotted horizontal line is ``correct'' threshold
    $t_{1/4} = \sqrt{2(1 - \tfrac{1}{4}) \log n} = 3.72,$ and dotted
    vertical lines show magnitude of the error committed with $t_{1/4}.$
    Solid horizontal line is a 'misspecified' threshold $t_{1/2} =
    \sqrt{2(1-\tfrac{1}{2}) \log n} = 3.03$, which would be the
    appropriate choice for $n_0 = n^{1/2} = 100$ non-zero components.
    Solid vertical lines show the \textit{additional} absolute error
    suffered by using this misspecified threshold.  Quantitatively, the
    absolute error $\| \hat \mu - \mu \|_1$ using the right threshold is
    14.4 versus 70.0 for the wrong threshold.  For $\ell_2$ error $\|
    \hat \mu - \mu \|_2^2$, the right threshold has error 38.8 and the
    wrong one has error 221.1.}
        \label{fig:errorillus}
      \end{center}
\end{figure}

Figure \ref{fig:errorillus} gives an example. One simple model of
varying sparsity levels sets $n_0 = n^\beta$ non-zero components out
of $n$, $ 0 < \beta < 1$. Theory, reviewed in Section 3, suggests
that a threshold of about $t_\beta = \sigma_n\sqrt{2(1-\beta) \log n}$ is
appropriate for such a sparsity level. Suppose that $\beta$ is
unknown, and examine the consequences of using misspecified thresholds
$t_\gamma$, $\gamma \neq \beta$.
The solid lines in Figure \ref{fig:errorillus} show the
\textit{increased} absolute error incurred using $t_{1/2}$ when
$t_{1/4}$ is appropriate -- the total absolute error is five times
worse. For squared error, the misspecified threshold produces a
discrepancy that is larger by nearly a factor of six.

Typically we could not know in advance the degree of
sparsity of the estimand, so
we prefer methods
adapting automatically to the unknown degree of sparsity.

\subsection{FDR-Controlling Procedures}

\citet{beho95} proposed a new principle for design of simultaneous
testing procedures  -- control of the False Discovery Rate (FDR).
In a setting where one is testing many hypotheses, the principle
imposes control on the ratio of the number of erroneously rejected
hypotheses to
the total number rejected. The exact definition and
basic properties of the FDR, as well as 
examples of procedures holding it below a
specified level $q$, are
reviewed in Section \ref{sec:test-est-fdr}.
In the context of \textit{estimation},
a thresholding procedure, which reflects the step-up FDR controlling
procedure in \citet{beho95}, was first proposed in \citet{abbe95}. The
procedure is quite simple:

Form the order statistics of the magnitudes of the observed estimates,
\begin{equation}
            |y|_{(1)} \geq |y|_{(2)} \geq  \ldots |y|_{(k)} \geq  \ldots \geq
|y|_{(n)},
      \label{eq:ordereddata}
\end{equation}
and compare them to the series of right tail Gaussian quantiles $t_k
= \sigma_n
z(q/2\cdot k/ n)$.
Let $\hat{k}_F$ be the largest index $k$ for which $|y|_{(k)} \geq t_k$;
threshold the estimates at (the data dependent) threshold $t_{\hat{k}_F}
= \hat{t}_F$,
\begin{equation}
        \label{eq:fdr-def}
\hat{\mu}_{F,k} =
          \begin{cases}
                y_k, & |y_k| \geq \hat{t}_F \\
                0, & \text{else}.
          \end{cases}
      \end{equation}
The FDR threshold is inherently adaptive to the sparsity level: it
is higher for
sparse signals and lower for dense ones. In the context of model selection,
control of the FDR means that when the
model is discovered to be complex, so that many variables are needed,
we should not be concerned unduly about occasional inclusion of unnecessary
variables; this is bound to happen. Instead, it
is preferable to control the proportion of erroneously included
variables. In a limited simulation study in the context of wavelet
estimation, \citet{abbe96}
demonstrated the good adaptivity properties of the FDR thresholding procedure
as reflected in relative mean square error performance.

\begin{figure}[htb]
      \begin{center}
        \leavevmode
     \centerline{\includegraphics[width= 8cm, height = 8cm,
     angle = 0]{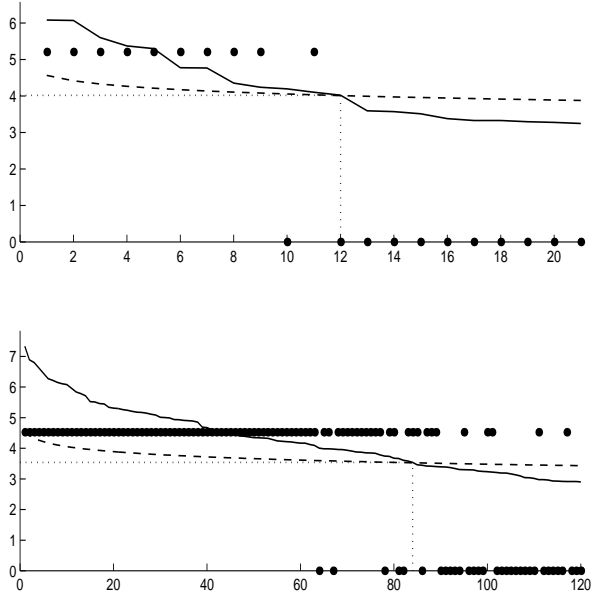}}
\caption{ \small (a) 10 out of 10,000. $\mu_i = \mu_0 \doteq 5.21$ for $i=1,
         \ldots, n_0 = 10$ and $\mu_i=0$ if $i=11, 12, \ldots , n =
         10,000.$  Data $y_i$ from model \eqref{eq:gaussshift},
         $\sigma_n = 1$. Solid
         line: ordered data $|y|_{(k)}$.  Solid circles: true unobserved
         mean value $\mu_i$ corresponding to observed $|y|_{(k)}$.
         Dashed line: FDR quantile
         boundary $t_k = z( q/2 \cdot k/n), \  q= 0.05.$ Last crossing at
         $\hat{k}_F = 12$ producing threshold $\hat{t}_F=4.02$. Thus
         $|y|_{(10)}$ and $|y|_{(12)}$ are false discoveries out of a
         total of $\hat{k}_F=12$ discoveries. The empirical false
         discovery rate $\hat{FDR}=2/12$.
(b) 100 out of 10,000. $\mu_i = \mu_0 \doteq 4.52$ for $i=1,
         \ldots, n_0 = 100;$ otherwise zero. Same FDR quantile boundary,
         $q = 0.05$. Now there are $\hat{k}_F = 84$ discoveries, yielding
         $\hat{t}_F=3.54 $ and $\hat{FDR}=5/84 $.  }
        \label{fig:twosparsities}
      \end{center}
\end{figure}

In order to demonstrate the adaptivity of FDR thresholding,
Figure \ref{fig:twosparsities} illustrates
the results of FDR thresholding at two
different sparsity levels. In the first, sparser, case a higher
threshold is chosen. Furthermore, the fraction of discoveries
(coefficients above threshold) that are false discoveries (coming from
coordinates with true mean $0$) is roughly similar in the two
cases. This is consistent with the fundamental result of
\citet{beho95} that the FDR procedure described above controls the
false discovery rate
below level $q$ whatever be the configurations of means
$\mu \in \mathbb{R}^n$, $n \geq 1$.

\subsection{Certainty-Equivalent Heuristics for FDR-based thresholding}

How can FDR multiple-testing ideas be related to the performance of the
corresponding estimator? Here we sketch a simple heuristic.

Consider an `in-mean' analysis of FDR thresholding.
In the FDR definition, replace the observed data
$|y|_{(k)}$ by the mean values $\bar \mu_k$, assumed to be already 
decreasing. Consider a pseudo-FDR index $k_* (\bar \mu )$,
found, assuming $\sigma_n = 1$, by solving for the crossing point
\begin{displaymath} \label{crossing}
      \bar \mu_{k_*}  = t_{k_*}.
\end{displaymath}

Consider the case where the object of interest obeys
the weak $\ell_p$ sparsity constraint $\mu \in m_p[\eta_n]$.
Weak $\ell_p$ has a natural `extreme' sequence, namely (\ref{eq:extremal}).
    Consider the `in-mean' behavior at this extremal sequence;
the crossing point relation (\ref{crossing}) yields
$$
\eta_n (n/k_*)^{1/p} = t_{k_*} .
$$
Using the relation $t_k \sim \sqrt{2\log n/k}$, valid for $k = o(n)$,
one sees quickly that
\begin{displaymath}
      t_{k_*} \sim \sqrt{2 \log \eta_n^{-p}} ;
\end{displaymath}
the right hand side of this display is asymptotic to 
the correct 
minimax threshold for weak and
strong $\ell_p$ balls of radius $\eta_n$!

Thus, FDR, in a heuristic, certainty-equivalent analysis, is able to
determine the threshold appropriate to a given signal sparsity. Further, this
calculation makes no reference to the loss function, and so we might
hope that the whole range $0 < r \leq 2$ of $\ell_r$ error measures
is covered.

\subsection{Main Results}

Given an $\ell_r$ error measure and $\Theta_n \subset \mathbb{R}^n$,
the worst-case risk of an estimator $\hat \mu$ over $\Theta_n$ is
\begin{equation}
          \bar \rho ( \hat \mu, \Theta_n) = \sup_{\mu \in \Theta_n} 
E_\mu \| \hat
          \mu - \mu \|_r^r
       \label{eq:ellrloss}
\end{equation}
The parameter spaces of interest to us will be those introduced earlier:
\begin{description}
\item{i)} $\Theta_n = \ell_0[\eta_n]$ (``nearly black''),
\item{ii)} $\Theta_n = m_p[\eta_n],\; 0 < p < r$ (weak-$l_p$ balls), and
\item{iii)} $\Theta_n = \ell_p[\eta_n],\; 0 < p < r$ (strong-$l_p$ balls).
\end{description}
In these cases we will need to have $\eta_n \rightarrow 0$ with increasing
$n$, reflective of increasing sparsity.

For a given $\Theta_n$, the minimax risk is the best attainable
worst-case risk:
      \begin{equation}
        \label{eq:minmaxdef}
R_n(\Theta_n) = \inf_{\hat \mu} \bar \rho (\hat \mu, \Theta_n);
      \end{equation}
the infimum covers all estimators (measurable functions of the data).
Any particular estimator, such as FDR, must have $\bar \rho (\hat \mu_F,
\Theta_n ) \geq R_n (\Theta_n)$, but we might ask how inefficient
$\hat \mu_F$ is relative to the ``benchmark'' for $\Theta_n$ provided
by $R_n(\Theta_n)$.

\begin{theorem}
      \label{th:mainopt}
      Let $y \sim N_n(\mu,\sigma_n^2I)$ and the FDR estimator $\hat \mu_F$ be
      defined by (\ref{eq:fdr-def}). In applying the FDR estimator,
      the FDR control parameter ($q_n$, say) may depend on $n$,
 but suppose this has a limit $q \in [0,1)$. In addition, 
    suppose 
$q_n \geq \gamma/{\log(n)}$ for some $\gamma > 0$
    and all $n 
\geq 1$.
     
Use the $\ell_r$ risk measure 
(\ref{eq:ellrloss})
      where $0 \leq p   < r \leq 2$.
      Let $\Theta_n$ be one of the parameter spaces detailed above with
      $\eta_n^p \in [n^{-1}\log^5n, n^{-\delta}],\;\delta>0$.
                     Then as $n \rightarrow \infty$,
     $$
      \sup_{\mu \in \Theta_n} \rho(\hat{\mu}_F,\mu) =
        R_n( \Theta_n) \{ 1 + u_{rp} \frac{(2q -1)_+}{1-q} + o(1) \},
     $$
where $u_{rp}=1$ and $u_{rp} =  1 - (p/r)$ for strong- and weak-$l_p$ balls
respectively.
\end{theorem}

Hence, if the FDR control parameter $q \leq 1/2$,
$\bar \rho ( \hat \mu_F, \Theta_n) \sim R_n (\Theta_n)$
in the sense that the ratio approaches $1$ as $n \rightarrow \infty$.
Otherwise, $\bar \rho ( \hat \mu_F, \Theta_n) \sim c(q) R_n (\Theta_n)$ for
an explicit $c(q) > 1$ growing with $q$.

In short Theorem \ref{th:mainopt} establishes the asymptotic minimaxity of
the FDR estimator in the setting of \eqref{eq:gaussshift} \-- provided we 
control false discoveries so that there are more true 
discoveries
than false ones. Moreover,
this minimaxity is {\em adaptive} across various losses and sparse
parameter spaces.

This exhibits a tighter connection between False Discovery Rate ideas and
adaptive minimaxity than one might have expected.  
The key parameter
in the FDR theory -- the rate itself -- seems to be diagnostic for
performance.

\subsection{Interpretations}

Two remarks help place the above result in context.

\subsubsection{Comparison with other estimators}

The result may be compared to traditional results in
the estimation of the multivariate normal mean.
Summarizing results given in \citet{dojo94}:
\begin{description}
\item{(i)} Linear estimators attain the wrong \textit{rates} of convergence
when $0 < p < r$, over these parameter spaces;

\item{(ii)} The James-Stein estimator, which is essentially a linear
estimator with data-determined shrinkage factor, has the
same defect as linear estimators;

\item{(iii)} Thresholding at a fixed level, say $\sigma_n \cdot
\sqrt{2 \log n}$,
\textit{does} attain the right rates, but with the wrong constants for
$0 < p < r$;

\item{(iv)} Stein's unbiased risk estimator (SURE)
directly optimizes the $\ell_2$ error, and is adaptive for $r=2$ and $1 < p
\leq 2$ \citep{dojo95a}. However, there appears to be a major technical
(empirical process) barrier to extending this result to $p \leq 1$,
and indeed, instability has been observed in such cases
in simulation experiments \citet{dojo95a}.
Further, there is no reason to expect that optimizing an $\ell_2$
criterion should also give optimality for $\ell_r$ error measures, $p < r < 2.$

\end{description}
In short, traditional estimators are not able to achieve the desired level
of adaptation to unknown sparsity.
On the other hand, recent work
by Johnstone and Silverman (2004), triggered by the present paper,
exhibits an empirical Bayes estimator -- EBayesThresh --
which seems, in simulations, competitive with FDR thresholding, although the
theoretical results for sparse cases are currently weaker.

\subsubsection{Validity of Simultaneous Minimaxity}

Minimax estimators are often criticised as being complicated,
counter-intuitive and distracted by irrelevant worst cases.
An often-cited example is
$\hat p = [ x + \sqrt{n}/2]/[n + \sqrt{n}]$ for estimating a success
probability $p \in [0,1]$ from $X \sim Bin(n,p)$. Although this
estimator is minimax for estimating $p$ under squared-error loss,
`everybody' agrees that the common sense estimator $\bar x = x/n$
is `obviously better' -- better at most $p$
and marginally worse only at $p$ near $1/2$.

Perhaps surprisingly, simultaneous (asymptotic) minimaxity
seems to avoid such objections. Instead, to paraphrase an old
dictum, it shifts the focus from an ``exact solution to the wrong
problem'' to ``an approximate solution to the right problem''.
To explain this, note that to develop a standard minimax solution, one
starts with parameter space $\Theta$ and error measure $\| \cdot \|$
and finds a minimax estimator $\hat \mu_{\Theta, \| \cdot \|}$
attaining the minimum in \eqref{eq:minmaxdef}.
This estimator may indeed be unsatisfactory in practice,
for example because it may depend on aspects of $\Theta$ that will not
be known, or may be incorrectly specified.

In contrast, we begin here with an a priori reasonable estimator
$\hat \mu_{F}$ whose definition does not depend on the imposed $\| \cdot
\|$ and
the presumably unknown $\Theta_n$. Adaptive minimaxity
for $q \leq 1/2$ \-- as established for
$\hat \mu_{F}$ in Theorem
\ref{th:mainopt} \-- shows that,
for a large class of relevant parameter spaces $\Theta_n$
and error measures $\| \cdot \|$,
$\bar \rho ( \hat \mu_n, \Theta_n) \sim R_n (\Theta_n)$.
In other words, the prespecified estimator $\hat \mu_n$ is flexible
enough to be approximately an optimal solution in many situations of
very different type (varying sparsity degree $p$, sparsity control $\eta_n$ and
error measure $r$ in the FDR example).

Using large $n$ asymptotics to exhibit approximately minimax solutions for
finite $n$ also renders the theory more flexible. For example in the
binomial setting cited earlier, the standard estimator $\bar x = x/n$,
while not exactly minimax for finite $n$,
\textit{is} asymptotically minimax.  More: if we consider
in the binomial setting the parameter spaces
$\Theta_{[a,b]} = \{ p: a \leq p \leq b \}$, then
$\bar x$ is simultaneously asymptotically minimax for
a very wide range of parameter spaces \-- each $\Theta_{[a,b]}$
for $0 < a < b < 1$ \--
whereas $\hat p$ is asymptotically minimax only for special
cases $a < 1/2 < b $. In short, whereas minimaxity
violates common sense in the binomial case, simultaneous asymptotic
minimaxity agrees with it perfectly.

\subsection{Penalized Estimators} \label{subsec.variational}

At the center of our paper
is the study, not of $\hat{\mu}_F$,
but of a family of {\it complexity-penalized}
estimators.  These yield approximations to FDR-controlling procedures,
but seem far more amenable to direct mathematical analysis.
Our study also allows us to exhibit connections of FDR control to
several other recently proposed model selection methods.

A penalized estimator is a minimizer of $\tilde \mu
\mapsto K(\tilde \mu,y)$, where
\begin{equation}
      \label{eq:penalized}
      K(\mu ,y ) = \| y - \mu \|_2^2 + Pen(\mu).
\end{equation}
If the \textit{penalty term} $Pen(\mu)$ takes an $\ell_p$ form,
$Pen(\mu) = \lambda \| \mu \|_p^p$, familiar estimators result:
$p=2$ gives linear shrinkage $\hat \mu_i = (1+\lambda)^{-1}y_i$;
while $p=1$ yields soft thresholding $\hat \mu_i = ( \text{sgn}\, y_i) ( |y_i|
- \lambda /2 )_+$; for $p=0$, $Pen(\mu) = \lambda \| \mu \|_0$,
gives hard thresholding $\hat \mu_i = y_i I \{ |y_i| \geq \lambda\}$.

\textit{Penalized FDR} results from modifying the penalty to
\begin{displaymath}
      Pen(\mu) = \sum_{l=1}^{\| \mu \|_0} t_l^2\
\end{displaymath}
Denote the resulting minimizer of \eqref{eq:penalized} by $\hat
\mu_2$. For small $\| \mu \|_0$, $Pen(\mu) \sim \ t_{ \| \mu \|_0}^2
\cdot \| \mu \|_0$.
It therefore has the flavor of an $\ell_0$-penalty, but with the
regularization parameter $\lambda$ replaced by the squared Gaussian
quantile appropriate to the complexity $\| \mu \|_0$ of $\mu$.
Further, $\hat \mu_2$ is indeed a variable hard threshold rule. If $
\hat k_2$ is a minimizer of
\begin{displaymath}
      S_k = \sum_{l=k+1}^n y_{(l)}^2 + \sum_{l=1}^k t_l^2,
\end{displaymath}
then $\hat \mu_{2,i} = y_i I \{ |y_i| \geq t_{\hat k_2} \}.$

The connection with original FDR arises as follows: $ \hat k_2$ is the
location of the \textit{global} minimum of $S_k$, while the FDR index
$\hat k_F$ is the rightmost \textit{local} minimum.
Similarly, we define $\hat k_G$ as the leftmost local minimum of
$S_k$: evidently $ \hat k_G \leq \hat k_2 \leq \hat k_F$.  For future
reference, we will call $\hat{k}_G$ the {\it Step-Down FDR} index. In practice,
these indices are often identical.
For theoretical purposes, we show (Proposition \ref{prop:sandwichbd}
and Theorem \ref{th:pen-to-fdr})
that $\hat k_F - \hat k_G$ is uniformly small enough on our sparse
parameter spaces $\Theta_n$  that asymptotic minimaxity
conclusions for $\hat \mu_2$ can be carried over to $\hat \mu_F.$

To extend this story from $\ell_2$ to $\ell_r$ losses, we make 
a straightforward translation:
\begin{equation}
          \hat \mu_r = \text{argmin}_\mu \ \| y - \mu \|_r^r + \sum_{l=1}^{\|
          \mu \|_0} t_l^r.
       \label{eq:fdr-r-def}
\end{equation}
Again it follows that $\hat k_r \in [ \hat k_G , \hat k_F].$
Our strategy is, first, to prove $\ell_r$-loss optimality results
using $\hat \mu_r$, and later, to draw parallel conclusions for the
original FDR rule $\hat \mu_F.$

Why is the penalized form helpful? In tandem with the definition of
$\hat \mu_r$ as the minimizer of an empirical complexity $ \tilde \mu
\mapsto K( \tilde \mu , y)$, we can define the minimizer $\mu_0$
of the \textit{theoretical} complexity $\tilde \mu \mapsto K
( \tilde \mu, \mu)$ obtained by replacing $y$ by its expected value $\mu$.
By the very definition of $\hat \mu_r$, we have
$K(\hat \mu_r, y) \leq K( \mu_0, y)$, and by simple manipulations one
arrives (in the $\ell_2$ case  here) at the basic bound,
valid \textit{for all} $\mu \in \mathbb{R}^n$:
\begin{equation}
      \label{eq:keybd}
      E \| \hat \mu_2 - \mu \|^2 \leq K( \mu_0, \mu) + 2 E \langle \hat \mu_2
      - \mu, z \rangle -E_\mu Pen(\hat \mu_2).
\end{equation}

Analysis of the individual terms on the right side is very revealing.
Consider the theoretical complexity term
$K(\mu_0,\mu)$. For $\Theta_n$ of type
(i)-(iii) in the previous section,
it turns out that the worst-case theoretical
complexity is asymptotic to the minimax risk!
Thus:
\begin{equation}
      \label{eq:complexity-to-mmx}
      \sup_{\mu \in \Theta_n} K(\mu_0,\mu) \sim R_n ( \Theta_n), \quad n \goto \infty.
\end{equation}
The argument for this relation is rather
easy, and will be given below in Section \ref{ssec:max-theoretical-complexity}.
The remaining term $2 E \langle \hat \mu_2
      - \mu, z \rangle -E_\mu Pen(\hat \mu_2)$ in \eqref{eq:keybd} has 
the flavor
of an error term of
lower order.
Detailed analysis is actually
rather hard work, however. Section \ref{ssec:error-term} overviews
a lengthy argument, carried out in the immediately
 following sections, showing that this error term is indeed
negligible over $\Theta_n$ if $q \leq 1/2$, and of the order of
$R_n ( \Theta_n)$ otherwise.

Plausibility for simultaneous asymptotic minimaxity of FDR
is thus aid out for us very directly within the penalized FDR point of view.
A full justification requires study of the
theoretical complexity and the error term respectively.
This fact permeates the architecture of the arguments to follow.

\medskip

\subsection{Penalization by $2 k \log(n/k)$}

Penalization connects
our work with a vast literature on model selection.
Dating back to Akaike (1973), it has been
popular to consider model selection rules
of the form
\[
      \hat{k} = \argmin_k  RSS(k)
       + 2 \sigma^2 \cdot k \cdot \lambda  ,
\]
where $\lambda$ is the penalization parameter
and $RSS(k)$ stands for ``the best residual
sum of squares $\| y -  m \|_2^2$ for a model
$m$ with $k$ parameters''.  The AIC model
selection rule takes $\lambda = 1$.
Schwarz' BIC model selection rule takes
$\lambda = \log(n)/2$, where $n$ is the sample size. Foster and George's
RIC model selection rule takes $\lambda = \log(p)$,
where $p$ is the
number of variables available
for potential inclusion in the model.

Several independent groups of researchers
have recently proposed model selection rules
with variable penalty factors.  For convenience,
we can refer to these as $2 \log(n/k)$ factors,
yielding rules of the form
\begin{equation} \label{log-n-k}
      \hat{k } = \argmin_k  RSS(k)
       + 2 \sigma_n^2 \cdot k \cdot \log(n/k)  .
\end{equation}

\begin{itemize}
\item
\citet{fost97} arrived at a penalty $\sigma^2\sum_1^{k} 2 \log( n/j)$ from
information-theoretic considerations.
Along sequences of $k$ and $n$ with $n \goto \infty$ and $ k/n \rightarrow 0,$
$2 k \log(n/k) \sim  \sum_{j = 1}^k 2 \log(n/j).$

\item
For prediction problems, \citet{tikn97} propose model selection using
a covariance inflation criterion which adjusts the training error by
the average covariance of predictions and responses on permuted
versions of the dataset. In the case of orthogonal regression, their
proposal takes the form of complexity-penalized residual sum of
squares, with the complexity penalty approximately of the above
form, but larger
by a factor of $2$:  $2 \sigma_n^2 \sum_{j = 1}^k 2 \log(n/j).$
There are intriguing parallels between the covariance expression for
the optimism \citep{efro86} in \citet[formula (6)]{tikn97} and the
complexity bound \eqref{eq:keybd}.

\item
\citet{foge00} adopt an empirical Bayes approach, drawing
the components $\mu_i$ independently from a mixture prior $(1-w) \delta_0 +
w N(0,C)$ and then estimating the hyperparameters $(w,C)$ from the
data $y.$
They argue that the resulting estimator penalizes the addition of a $k$th
variable by a quantity close to $2 \log (\frac{n+1}{k} - 1).$

\item
\citet{BirgeMassart}
studied complexity-penalized model selection for a class of
penalty functions, including penalties
of the form $2 \sigma^2_n k \log (n/k).$ They develop
non-asymptotic risk bounds for such procedures over $\ell_p$ balls.

\end{itemize}

Evidently, there is substantial interest
in  the use of variable-complexity penalties.
There is also an {\it extensive similarity}
of $2k\log(n/k)$ penalties to FDR penalization.
Penalized FDR $\hat{\mu}_2$ from (\ref{eq:fdr-r-def}) can be written
in penalized form with a variable-penalty factor $\lambda_{k,n}$:
\[
      \hat{k}_2 = \argmin_k  RSS(k)
       + 2 \sigma_n^2  k \lambda_{k,n}  ,
\]
where
\[
     \lambda_{k,n} = \frac{1}{2k} \sum_{l=1}^{k} z^2\left(\frac{lq}{2n}\right)\;\sim \;
     z^2\left(\frac{kq}{2n}\right)/2
\;\sim\; \log(n/k)   - \frac{1}{2} \log \log(n/k)  + c(q,k,n)
\]
for large $n$, $k=o(n)$, and bounded remainder $c$ (compare (12.7) below).
FDR penalization is thus slightly weaker than $2k \log(n/k)$ penalization.
We could also say that $2k \log(n/k)$ penalties 
have a formal algebraic similarity to FDR penalties, but require 
a variable $q = q(k,n)$ that is both small and decreasing with $n$.
This perspective on $2k\log(n/k)$ penalties,
suggests the following conjecture:
\begin{conjecture}
In the setting of this paper, where `model selection' means
adaptive selection of nonzero means, and the underlying estimand $\mu$
belongs to one of the parameter spaces as detailed in
Theorem \ref{th:mainopt}, the procedure (\ref{log-n-k})
is asymptotically minimax, simultaneously over the full
range of parameter spaces and losses covered in that theorem.
\end{conjecture}

In short, although the $2k \cdot \log(n/k)$ rules were not proposed
from a formal decision-theoretic, they might well
exhibit simultaneous asymptotic minimaxity. We suspect that the methods
developed in this paper may be extended to yield a proof of this conjecture.

\subsection{Take-Away Messages}

The theoretical results in this paper
suggest the following two messages:
\begin{description}
\item{TAM 1.} FDR-based thresholding gives an optimal
way of adapting to unknown sparsity: choose  $q \leq 1/2$.
In words, aiming for fewer false discoveries than true ones yields sharp asymptotic
minimaxity.
\item{TAM 2.} Recently proposed $2k\log(n/k)$ penalization schemes,
when used in a sparse setting, may be viewed as similar to FDR-based
thresholding.
\end{description}

\subsection{Simulations}

We tested FDR thresholding
and related procedures in simulation experiments.
The outcomes support TAMs 1 and 2.

Table \ref{tbl:performance} displays results from simulations
at the so-called least-favorable
case $\mu_k = \min \{  n^{-1/2} k^{-1/p}, \sqrt{ (2-p) \log n } \}$
for the weak-$\ell_p$ parameter ball (compare remark following (\ref{eq:kstareqn})).
Here $p = 1.5$, $r=2$, $n=1024$
and $n = 65536$, $\sigma = 1$.
The table records the ratio of
squared-error risk of FDR
to squared-error risk
of the asymptotically optimal threshold $t^* = t^*(p,n) =  \sqrt{(2 - p) \log n}$ 
for that parameter ball (compare Section 3.3 below).
All results derive from $100$ repeated experiments. The standard errors of the MSEs 
were between 0.001-0.003 for $n=1024$ and between 0.0005-0.0007 for $n=65536$.

These results should be compared with the behavior
of $2 \log(n/k)$-style penalties. For the
estimator of \citet{fost97},  minimizing
$RSS + \sigma^2\sum_{j=1}^k 2\log(n/j)$,
we have that for $n=1024$,  $MSE/MSE(t^*) = 1.2308$
while for $n=65536$, $MSE/MSE(t^*) = 1.2281$.
This is consistent with behavior that would result from FDR control
with  $q=.3$
for $n=1024$ and $q = .25$ for $n=65536$.

\begin{center}
\begin{table}
\begin{tabular}{|l|c|c|c|c|}
\hline
& q & step-up FDR & penalized FDR & step-down FDR \\
\hline
n=1024 & 0.01 & 1.3440 & 1.3440 & 1.3440 \\
         & 0.05 & 1.3283 & 1.3293 & 1.3334 \\
         & 0.25 & 1.2473 & 1.2482 & 1.2512 \\
         & 0.40 & 1.2171 & 1.2171 & 1.2173 \\
         & 0.50 & 1.2339 & 1.2335 & 1.2321 \\
         & 0.75 & 1.4159 & 1.4132 & 1.4100 \\
         & 0.99 & 1.9810 & 1.9744 & 1.9687 \\
\hline
n=65536& 0.01 & 1.3370 & 1.3372 & 1.3374 \\
         & 0.05 & 1.3178 & 1.3180 & 1.3183 \\
         & 0.25 & 1.2276 & 1.2277 & 1.2277 \\
         & 0.40 & 1.1889 & 1.1889 & 1.1890 \\
         & 0.50 & 1.1937 & 1.1936 & 1.1936 \\
         & 0.75 & 1.5122 & 1.5118 & 1.5114 \\
         & 0.99 & 4.0211 & 4.0189 & 4.0174 \\
\hline
\end{tabular}
\caption{Ratios of MSE(FDR)/MSE($t^*(p,n)$), $p=3/2$.}
\label{tbl:performance}
\end{table}
\end{center}

\begin{center}
\begin{figure}
\centering
\includegraphics[scale=.40]{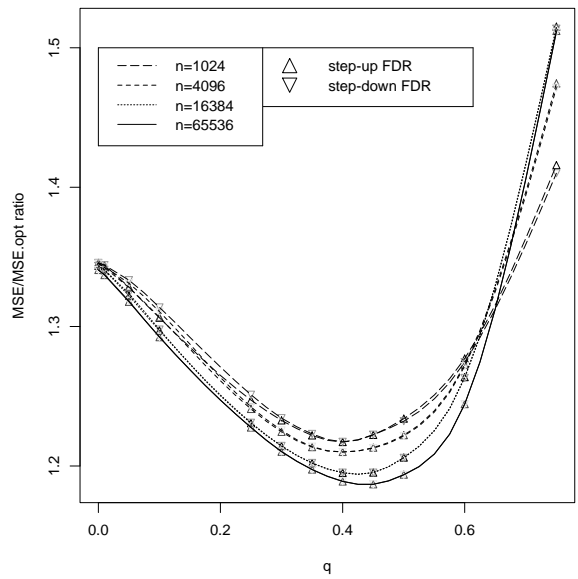}  
\caption{Ratios of MSE(FDR)/MSE($t^*(p,n)$), $p=3/2$.}
\label{fig:performance}
\end{figure}

\end{center}

In Figure \ref{fig:performance}, we display simulation results
under a range of sample sizes.  Apparently  the minimum
MSE occurs somewhere below $q = 1/2$.
We propose the following
interpretations:
\begin{description}
\item{INT 1.} FDR procedures with $q \leq 1/2$
have a risk which is a reasonable multiple of
the `ideal risk' based on the threshold which would have been
optimal for the given sparsity of the object. The ratios in
Table \ref{tbl:performance} do not differ much for various $q \leq 1/2$ that
demonstrates robustness of the FDR procedures towards the choice of $q$.

\item{INT 2.} An FDR procedure with $q$ near $1/2$ appears
to outperform $q$-small procedures at this configuration, and achieving
risks which roughly comparable to the ideal risk.

\item{INT 3.} Avoid FDR procedures with large $q$, in favor of
$q \leq 1/2$.
\end{description}

\subsection {Contents}

The paper to follow is far more technical than the introduction;
in our view necessarily so, since much of the work concerns refined
properties of fluctuations
in the extreme upper tails of the normal order statistics.
However, Sections 2-4 should be accessible on a first reading.
They review pertinent information about FDR-controlling procedures,
about minimax estimation over $\ell_p$ balls, and parse our main
 result into an Upper Bound result and a Lower Bound result.
 Section 4 then gives an overview of the paper to follow,
 which carries out rigorous proofs of the Lower Bound (Sections 5-8, 13)
 and the Upper Bound (Sections 9-11).

\section{The False Discovery Rate}
\label{sec:test-est-fdr}
\setcounter{equation}{0}

The field of Multiple Comparisons has developed
many techniques to control
the increased rate of type I error when testing a family of $n$
hypotheses  $H_{0i}$ versus $H_{1i},$ $i=1,2,\ldots ,n$.
The traditional approach is to
control the {\it familywise error rate} at some level $\alpha$, that
is, to use a
testing procedure that controls at level $\alpha$
the probability of erroneously
rejecting even one true $H_{0i}$.  The venerable Bonferroni
procedure tests ensures this by testing each hypothesis at the
$\alpha/n$ level.

The Bonferroni procedure is criticised as being too conservative, since it
often lacks power to detect the alternative hypotheses. Much research has
been devoted to devise more powerful procedures: tightening the probability
inequalities, and incorporating the dependency structure when it is known.
For surveys, see \citet{hota87} and \citet{shaf95}.
In one fundamental sense the success has been limited. Generally the
power deteriorates substantially when the problem is large. As a result,
many practitioners avoid altogether using any multiplicity adjustment to
control for the increased type I errors caused by simultaneous inference.

\citet{beho95} argued that the control of the familywise error-rate
is a very conservative goal which is not always necessary.
They proposed to control the expected ratio of the number
of erroneously rejected hypotheses to
the number rejected  - the False Discovery Rate (FDR).
Formally,
for any fixed configuration of true and false hypotheses,
let $ V$ be the number of true null hypotheses erroneously rejected,
among the $R$ rejected hypotheses. Let $Q$ be $V/R$ if $R > 0$, and $
0$ if $R=0$; set
$FDR = E \{ Q\}$,
where the expectation is taken according to the same configuration.
The FDR is equivalent to the familywise error-rate when all tested
hypotheses are true, so an FDR-controlling procedure  at level $q$ controls
the probability of making even one erroneous discovery in such a situation.
Thus for many problems the value of $q$ is naturally chosen at the
conventional levels for tests of significance. The FDR of a multiple-testing
procedure is never larger than the familywise error-rate.
Hence controlling FDR admits more powerful
procedures.

Here is a simple {\it step-up} FDR-controlling procedure.
Let the individual $P-$values for the hypotheses $H_{0i}$, be arranged
in ascending order: $P_{[1]} \leq \ldots \leq P_{[n]}.$ Compare the
ordered $P-$values to a linear boundary $i/n q$, and note the last crossing
time:
\begin{equation}\label{FDRdef}
           \hat{k}_F = \max \{ k: P_{[k]} \leq (k/n) q \}.
\end{equation}
The FDR multiple-testing procedure is to reject all hypotheses $H_{(0i)}$
corresponding to the indices $i = 1, \ldots , \hat{k}_F.$
If $\hat{k}_{B}$ denotes the number of $P-$values below the
Bonferroni cutoff $q/n$ it is apparent that $\hat{k}_F \geq
\hat{k}_{B}$ and hence that the FDR test conducted at the same level is
necessarily less conservative.

\citet{beho95} considered the above testing procedure
in the situation of independent hypothesis tests
on many individual means. They considered the two-sided $P$-values
from testing that each individual mean was zero. They found that the false
discovery rate of the above multiple testing procedure is bounded by $q$
\textit{whatever be the number of true null hypotheses $ n_0 $ or the
configuration of the means under the alternatives}:
\begin{equation} \label{FDRdefQ}
            FDR = E_{\mu} \{ Q\} = q n_0/n \qquad\leq q\qquad \qquad \text{for
all } \mu \in
           \mathbb{R}^n.
\end{equation}

The multiple-testing procedure (\ref{FDRdef})
was proposed informally by Elkund, by
\citet{seeg68} and much later independently by \citet{sime86}.
Each time it was neglected
because it was shown not to control the familywise error-rate
[\citet{seeg68},\citet{homm88}]. In the absence of the FDR
concept, it was not understood why  this procedure could be
a good idea. After introduction of the FDR concept, 
it was recognized that
$\hat{k}_F$ had the  FDR 
property, but also that other procedures
offered FDR control \-- most importantly for us, 
the step-down estimator $\hat{k}_G$; \citet{Sarkar2002}. 
This rule, 
introduced in Section 1, will also be referred to 
frequently below,
and our theorems are also applicable to 
thresholding
estimators based upon it.

As noted in the introduction, \citet{abbe95} adapted FDR testing
to the setting of estimation, in particular of
wavelet coefficients of an unknown regression function.
In this setting, given $n$ data on a unknown function
observed in Gaussian white noise,
we are testing $n$ independent null hypotheses on function's wavelet
coefficients $\mu_i = 0$.
Using the above formulation with two-sided P-values, we obtain
(\ref{eq:fdr-def}).

Previously in the same setting of wavelet estimation,
\citet{dojo94} had proposed to estimate wavelet coefficients
by setting to zero all coefficients below a certain ``universal threshold''
$\sqrt{2 \log(n)} \sigma_n$. A key
observation in \cite{dld95} and \cite{djkp93},
about this threshold is that, with high probability,
every truly zero wavelet coefficient is estimated by zero.

Using ideas from simultaneous
inference we can look at universal thresholding differently.
The likelihood ratio test of the null
hypothesis $H_{0i}: \mu_i =0$ rejects if and only if $|y_i| > t \sigma,$
and the Bonferroni method at familywise level $\alpha$
sets the cutoff for rejection $t$ at $t_{BON} = \sigma
z(\alpha/2n)$. Now very roughly, $z(1/n) \sim \sqrt{2 \log(n)}$;
much more precise results are derived below and lie at the center
of our arguments. Hence, Bonferroni at any reasonable
level $\alpha$ leads us to set a
threshold not far from the universal threshold.
Put another way, universal thresholding may be viewed
as precisely a Bonferroni
procedure, for $\alpha = \alpha_n^U$. We can derive
$\alpha_n^U \asymp 1/\sqrt{ \log(n) }$ as $n \goto \infty$.

As was emphasised by \cite{abbe96}, the FDR estimator can choose
lower thresholds than $\sigma_n \cdot \sqrt{2 \log n}$ when
$\hat{k}_F$ is relatively
large. It thus offers the possibility of adapting to the unknown mean
vector by adapting to the data, choosing less conservative thresholds when
significant signal is present. It is this possibility we explore here.

{\it Pointers to the FDR literature more generally.}
The above 
discussion of FDR thresholding
emphasizes just that `slice' of the 
FDR literature needed
for this paper, so it is highly selective.
The 
literature of FDR methodology is growing rapidly,
and is too diverse 
to adequately summarize here.
Recent papers have illuminated the FDR 
from different points
of view: asymptotic, Bayesian, Empirical Bayes, and as the limit of
empirical processes (Efron et al. (2001), Storey (2002), Genovese and Wasserman
(2002).

Another line of work, starting with Benjamini and Hochberg (2000),
addresses the factor $n_0/n$ in (\ref{FDRdefQ});
many methods have been offered to estimate this, followed by
the step-up procedure with the adjusted (larger) $q$.
The results of Benjamini et al. (2001) and Storey et al. (2004)
assure FDR control under independence.
When $n_0/n$ is close to 1, as is our case in this paper,
such methods are close to the original step-up procedure.

An immediate next step
beyond this paper would be to study
dependent situations. The 
FDR-controlling property 
of the step-up procedure under positive 
dependence has
been established in Benjamini and Yekutieli (2001),
and similar results were derived for the step-down version
in Sarkar (2002).  Since much of the formal structure below
is based  on marginal properties of the observations,
this raises the 
possibility that our estimation results would
extend to a broader 
class of situations involving
dependence in the noise terms $z_i$.

\section{Minimax estimation on $\ell_0, \ell_p, m_p$}
\label{sec:mmx-on-lp}
\setcounter{equation}{0}

\medskip

As a prelude to the formulation of the adaptive minimaxity results,
we review information
\citep{djhs92,dojo94,john94}
on minimax estimation over $\ell_0, \ell_p$ and weak
$\ell_p$ balls in the sparse case: $0 < p < 2$ and with normalized
radius $\eta_{n} \rightarrow 0$ as $n \rightarrow \infty.$ Throughout
this section, we suppose a shift Gaussian model
\eqref{eq:gaussshift} with unit noise level $\sigma_n = 1.$
We will denote the risk of an estimator $\hat{\mu}$ under $\ell_r$
loss by
\begin{displaymath}
       \rho(\hat{\mu},\mu) = E_{\mu} \| \hat{\mu} - \mu \|_r^r.
\end{displaymath}

Particularly important classes of estimators are obtained by
thresholding of individual co-ordinates: hard thresholding was defined
at \eqref{eq:hardthresh}, while soft
thresholding of a single
co-ordinate $y_1$ is given by $\eta_S(y_1,t) = \text{sgn}(y_1) (|y_1| - t)_+.$
We use a special notation for the risk function of thresholding on a
single scalar observation $y_1 \sim N(\mu_1,1)$:
\begin{displaymath}
           \rho_{S}(t,\mu_1) = E_{\mu_1} |\eta_S(y_1,t) - \mu_1|^r
\end{displaymath}
with an analogous definition of $\rho_H(t,\mu_1)$ for hard thresholding.

\subsection{$\ell_0$ balls}
\label{ssec:l0balls}

Asymptotically least-favorable configurations for $\ell_0$
balls $\ell_0[\eta_n]$ can be built by drawing the $\mu_i$ i.i.d.  from
sparse two-point prior distributions
\begin{displaymath}
           \pi = (1 - \eta_n) \delta_0 + \eta_n \delta_{\mu_n}, \qquad
           \mu_n \sim (2 \log \eta_n^{-1})^{1/2}.
\end{displaymath}
The precise definition of $\mu_n$ is given in the Remark below.
The expected number of non-zero components $\mu_i$ is $k_n = n
\eta_n$.
The prior is constructed so that the corresponding Bayes estimator
essentially estimates zero even
for those \( \mu_i \) drawn from the atom at \( \mu_n \), and so the
Bayes estimator has an $l_r$ risk of at least $k_n \mu_n^r.$
A corresponding asymptotically minimax estimator is given by soft or
hard thresholding at threshold
$\tau_\eta = \tau(\eta_n) := (2 \log
\eta_n^{-1})^{1/2} \sim \mu_n$ as $ n \rightarrow \infty$.
This estimator achieves the precise asymptotics of the minimax
risk, namely:
          \begin{equation}
           R_n(\ell_0[\eta_n]) \sim k_n \mu_n^r = n \eta_n \mu_n^r
          \sim n \eta_n (2 \log \eta_n^{-1} )^{r/2}.
                   \label{l0mmx}
           \end{equation}

\medskip

\textit{Remark.}
Given a sequence $a_n^2 = o(\log \eta^{-1}_n)$ that
increases slowly to $\infty$, $\mu_n$ is defined as the solution of the
equation $\phi(a_n + \mu_n) = \eta_n \phi(a_n),$ where $\phi$ denotes
the standard Gaussian density.
Equivalently,
$\mu_n^2 + 2 a_n \mu_n = 2 \log \eta_n^{-1} = \tau_\eta^2$
giving the more precise relation
\begin{equation}
           \tau_\eta = \mu_n + a_n + o(a_n).
           \label{eq:texpans}
\end{equation}
Thus $\tau_\eta - \mu_n \rightarrow \infty$, which  for both
soft or hard thresholding at $\tau_\eta$ indicates $\rho(\tau_\eta;\mu_n)
\sim \mu_n^r.$ [Note also that, to simplify notation, we
are using $\mu_n$ to denote a sequence of constants rather than the
$n^{\text{th}}$ component of the vector $\mu$].

\subsection{$\ell_p$ balls}
\label{ssec:lpballs}

Again, asymptotically least favorable configurations for
$\ell_{p}[\eta_n]$ are obtained by i.i.d. draws  from
$ \pi = (1 - \beta_n) \delta_0 + \beta_n \delta_{\mu_n}$,
where now the mass of the non-zero atom and its location 
are,
informally given by the pair of properties
           \begin{equation}
                   \beta_n = \eta_n^p \mu_n^{-p}, \qquad \qquad \mu_n \sim (2
\log
                   \beta_n^{-1})^{1/2}.
                   \label{betadef}
           \end{equation}
More precisely  $\mu_n = \mu_n(\eta_n,a_n;p)$ is now
the solution of $\phi(a_n + \mu_n) =
\beta_n \phi(a_n)$, which implies that
\begin{equation}
           \mu_n \sim \tau_\eta = ( 2 \log \eta_n^{-p})^{1/2} \qquad \qquad
           n  \rightarrow \infty,
           \label{munasymp}
\end{equation}
and then that \eqref{eq:texpans} continues to hold for $\ell_p$ balls.
The expected number of non-zero components $\mu_i$ is now
$n \beta_n = n \eta_n^p \mu_n^{-p}$. For later use, we define
\begin{equation}
   \label{eq:kndef}
   k_n = n \eta_n^p \tau_\eta^{-p},
\end{equation}
since $\mu_n \sim \tau_\eta$, we have $k_n \sim n \beta_n$, and so
$k_n$ is effectively the non-zero number.
With similar heuristics for the Bayes estimator, the exact asymptotics
of minimax risk becomes
\begin{equation}
   \label{eq:lpmmxrisk}
     R_n(\ell_p[\eta_n]) \sim k_n \mu_n^r =
        n \eta_n^p \tau_\eta^{r-p}
        = n \eta_n^p ( 2 \log \eta_n^{-p})^{(r-p)/2}.
\end{equation}
Asymptotic minimaxity
is had by thresholding
at $t_{\eta_{n}} = (2 \log \eta_n^{-p})^{1/2} \sim (2 \log
n/k_n)^{1/2}$.

\subsection{Weak $\ell_p$ balls}
\label{ssec:wk-lpballs}

The weak $\ell_p$ ball $m_p[\eta_n]$ contains the corresponding strong
$\ell_p$ ball $\ell_p[\eta_n]$ with the same radius, and the
asymptotic minimax risk is larger by a constant factor:
           \begin{equation}
            R_n(m_p[\eta_n]) \sim (r/(r-p)) R_n(\ell_p[\eta_n]), \qquad \qquad n
                                   \rightarrow \infty.
                           \label{eq:rmpn}
                   \end{equation}
Let $F_p (x) = 1 - x^{-p}, x \geq 1$ denote the distribution function of the
Pareto($p$) distribution and let $X$ be a random variable having this
law. Then an asymptotically least favorable distribution for
$m_p[\eta_n]$  is given by drawing $n$ i.i.d. samples from the
univariate law
           \begin{equation}
                   \pi_1 = \mathcal{L} ( \min(\eta_n X, \mu_n) ),
                   \label{weak-least-fav}
           \end{equation}
where $\mu_n$ is defined exactly as in the strong case.  The mass of
the prior probability atom at $\mu_n$ equals $\int_{\mu_n}^{\infty}
F_p(dx/\eta_n) = \eta_n^p \mu_n^{-p} = \beta_n$, again as in the
strong case.  Thus, the weak prior can be thought of as being obtained
from the strong prior by smearing the atom at 0 out over the interval
$[\eta_{n},\mu_n]$ according to a Pareto density with scale $\eta_n.$
One can see the origin of the extra factor in the minimax risk
from the following outline (for details when $r=2$, see \citet{john94}).
The minimax theorem says that $R(m_p[\eta_n])$ equals the Bayes risk
of the least favorable prior. This prior is roughly
the product of $n$ copies of $\pi_1$, and the corresponding  Bayes
estimator is approximately (for large $n$) soft thresholding at
$\tau_\eta,$ so
\begin{displaymath}
           R_n(m_p[\eta_n]) \sim n \int \rho_S(\tau_\eta,\mu) \, \pi_1(d \mu).
\end{displaymath}
Now consider an approximation to the risk function of soft
thresholding, again at threshold $t_\eta.$ Indeed, using the estimate
$\rho_S(t,\mu) \doteq \rho_S(t,0) + |\mu|^r$, appropriate in the
range $0 \leq \mu
\leq \mu_n,$ ignoring the term \( \rho_S(t,0) \) and reasoning as before
\eqref{eq:lpmmxrisk} we find
\begin{align}
     R_n(m_p[\eta_n]) & \sim n \int_{\eta_n}^{\mu_n} \rho_S(t_\eta,\mu)
     F_p(d\mu/\eta_n) + n \beta_n  \rho_S(t_\eta,\mu_n)  \label{risk1} \\
     & \doteq n \int_{\eta_n}^{\mu_n} \mu^r \cdot p \eta_n^p \mu^{-p-1}
     d\mu +  k_n \mu_n^r \label{risk1.5} \\
     & \doteq [ \frac{p}{r-p} + 1] n \eta_n^p \mu_n^{r-p} \sim
     \frac{r}{r-p} R(\ell_p[\eta_n]).  \label{risk2}
\end{align}
Comparison with \eqref{eq:lpmmxrisk} shows that the second term in
\eqref{risk1} - \eqref{risk1.5} corresponds exactly to
$R(\ell_p[\eta_n])$; the first term is contributed by
 the Pareto density in the weak-$\ell_p$ case.

\section{Adaptive Minimaxity of FDR Thresholding}
\label{sec:adaptive-minimaxity}
\setcounter{equation}{0}

We now survey the path to our main result, providing
in this section  an overview of the remainder of the paper
and the arguments to come. 

What we ultimately prove is broader than the result given in 
the introduction, and the argument will develop several 
ideas seemingly of broader interest.

\subsection{General Assumptions}
\label{sec:GenAssump}

Continually below we invoke a collection
of assumptions (Q),(H),(E), and (A) defined here. 
\bitem
\item {\it False discovery control.}  We allow 
false-discovery rates
to depend on $n$, but approach a limit as $n 
\goto \infty$.  Moreover,
if the limit is zero, rates should not go to 
zero very fast. Formally 
define the assumption:

\textbf{(Q)} \ \ %
Suppose that
$q_n \rightarrow q \in [0,1)$.
If $q = 0,$ assume that $q_n \geq b_1/ \log n.$

The constant $b_1 > 0$ is arbitrary; its value
could be important at a specific $n$.

\item {\it Sparsity of the estimand.} 
We consider only parameter
sets which are sparse, and we place 
quantitative upper bounds
keeping them away from the `dense' case. 
Formally  define:

\textbf{(H)} \ \  
$\eta_n$ (for $\ell_0[\eta_n]$) and
$\eta_n^p$ (for $m_p[\eta_n]$) lie in the interval $[n^{-1} \log^5 n,
b_2 n^{-b_3}]$.

Here the constants $b_2 > 0$ and $b_3 > 0$ are arbitrary;
their 
chosen values could again be important in finite samples 

\item {\it 
Diversity of Estimators.}  Our results apply not just to
the 
usual FDR-based estimator $\hat{\mu}_F$ of (1.7) but also
the 
penalty-based estimator $\hat{\mu}_r$ of (1.11). 
More generally, 
recall  the terms $[\hat{t}_F, \hat{t}_{G}]$ defined in
Section \ref{subsec.variational}. Under formal
assumption {\bf (E)}, 
we  consider any estimator $\hat{\mu}$
obeying 
\begin{equation}
           \hat{\mu} = \mbox{ hard thresholding at \ }
            \hat{t} \in [ \hat{t}_F,
           \hat{t}_{G}] \ \ \text{w.p. \ } 1.
           \label{eq:muhatlim}
\end{equation}

\item {\it Notation.} We introduce a sequence $\alpha_n$ 
which often appears in estimates in Sections 7 and 8 and in dependent 
material. We also define constants $q'$ and $q^{\prime\prime}$. 
Formally:\\
\textbf{(A)} \ \ Set $\alpha_n = 1/(b_4 
\tau_\eta)$, with $b_4 =
(1-q)/4$. Also set $q' = (q+1)/2$ and $q^{\prime\prime}= (1-q)/2 = 1 -q'.$

\eitem

Finally, as a global matter, we suppose that  our 
observations $y \sim N_n(\mu,I)$; thus $\sigma_n^2 \equiv 1$.
For estimation of $\mu$, we consider $\ell_r$ risk
\eqref{eq:ellrloss}, $0 < r \leq 2$, and minimax risk $R_n(\Theta_n)$
of \eqref{eq:minmaxdef}. Here the parameter spaces are $\Theta_n =
\ell_0[\eta_n]$ or   $\ell_p[\eta_n]$ or
$m_p[\eta_n]$ defined by 
\eqref{eq:l0balldef}, \eqref{eq:lpballdef} and
\eqref{eq:mpballdef} respectively, with $0 < p < r$.

\subsection{Upper Bound Result}

Our argument for Theorem 
\ref{th:mainopt}
splits into two
parts, beginning with an upper bound 
on minimax risk:

\begin{theorem}
\label{th:mainres}

Assume (H),(E),(Q).
Then, as $n \rightarrow \infty$,
\begin{equation}
\sup_{\mu \in \Theta_n} \rho(\hat{\mu},\mu) \leq
      R_n( \Theta_n) \{ 1 + u_{rp} \frac{(2q_n -1)_+}{1-q_n} + o(1) \},
            \label{eq:mainbound_mp}
\end{equation}
where $u_{rp} =  1 - (p/r)$ if $\Theta_n = m_p[\eta_n]$
and $u_{rp} = 
1$ if $\Theta_n = \ell_p[\eta_n]$ or 
$\ell_0[\eta_n]$.
\end{theorem}

 The bare bones of our 
strategy for proving the
 upper bound result were described in the 
Introduction.
 The global idea is to study the penalized FDR 
estimator 
 $\hat{\mu}_2$ of (\ref{subsec.variational}) and then 
compare
 to the behavior of $\hat{\mu}_F$.  To make this
 work, 
numerous technical facts
 will be needed concerning the behavior of 
 hard thresholding, the mean and fluctuations
 of the threshold 
exceedance process, and so on.   As it turns
 out, those same 
technical facts form the core of
 our lower bound on the risk 
behavior of $\hat{\mu}_F$.
  As a result,  it is convenient for us to 
study the
  lower bound and associated technical machinery first,
in 
Sections 5-8 (with some details deferred to
Sections 12 and 13), and then later, in Sections 9-11, to
prove the upper 
bound, using results and viewpoints
established 
earlier.

\subsection{Lower Bound Result}

Theorem
\ref{th:mainopt} 
is completed by a lower bound
on the behavior of the FDR 
estimator.

\begin{theorem}
\label{th:lowerbound}
Suppose (H),(Q). With notation as in Theorem \ref{th:mainres},
\begin{equation}
\sup_{\mu \in \Theta_n} \rho(\hat{\mu}_F,\mu) \geq
      R_n( \Theta_n) \{ 1 + u_{rp} \frac{(2q_n -1)_+}{1-q_n} + o(1) 
\},
  \quad n \rightarrow \infty,
            \label{eq:lowbound_mp}
\end{equation}
where $u_{rp} =  1 - (p/r)$ for $\Theta_n = m_p[\eta_n]$,
and $u_{rp}=1$  for  $\Theta_n = \ell_p[\eta_n]$.
\end{theorem}

This bound  establishes the importance of $q$; showing 
that if $q > 1/2$,
then certainly FDR cannot be asymptotically 
minimax.
We turn immediately to its proof.

\section{Proof of the 
Lower Bound}

The proof involves three 
technical but significant ideas.
First,  it bounds the number of 
discoveries
made by FDR, as a function of the underlying means 
$\mu$.
Second, it studies the risk of ordinary hard thresholding
with 
non-adaptive threshold in a specially chosen, 
(quasi-) 
least-favorable one-parameter
subfamily of $\Theta_n$. Finally, it 
combines these
elements to show that, on this least-favorable 
subfamily,  
$\hat{\mu}_F$ behaves like hard thresholding
with a 
particular threshold.  The lower bound result
then follows. 
Unavoidably, the results in this section will
invoke lemmas and 
corollaries only proven in later sections.

Beyond simply proving the 
lower bound,
this section introduces some basic
viewpoints and 
notions. These include
\bitem
\item  A threshold exceedance function 
$M$,
which counts the number of threshold exceedances
as a function 
of the underlying means vector.
\item A special `coordinate system' 
for thresholds,
mapping thresholds $t$ onto the scale of 
relative
expected exceedances.
\item A special one-parameter (quasi-) 
least-favorable family
for FDR, at which the lower bound is 
established.
\eitem
These notions will be used heavily in later 
sections.

\subsection{Mean Exceedances and Mean Discoveries}
\label{ssec:mean-exc-disc}

Define the exceedance number \( N(t_{k}) = \# \{ i : |y_{i}| \geq
|t_{k}| \} \).
Since $|y_{(k)}| \geq |t_k|$ if and only if $N(t_k) \geq k$,
Thus, we are interested in the values of $k$ for which $N(t_k)
\approx k$.
(See Section \ref{sec:ldproof} for details).

Throughout the paper we will refer to the
  mean threshold exceedance 
function,  counting the mean number of
exceedances over threshold $t_k$ as $k$ varies:
         \begin{equation}
                 M(k;\mu) = E_\mu N(t_k) = \sum_{l=1}^n P_\mu ( |y_l| \geq t_k )
                         = \sum_{l=1}^n \Phi( [ \mu_l-t_k,\mu_l +t_k]^c ).
                 \label{eq:mnkmu}
         \end{equation}
[Here $\Phi(A)$ denotes the probability of event $A$ under the standard
Gaussian probability distribution, and $k$ is extended from positive
integer to positive \textit{real} values.]
If $\mu = 0,$ then $M(k;\mu) = 2n \tilde{\Phi} (t_k) = qk,$
reflecting the fact that in the null case, the fraction of
exceedances is always governed by the FDR parameter $q$. If $\mu \neq 0$,
we expect that $\hat{k}_F $ will be close to the \textit{mean discovery number}
\begin{equation}
     k(\mu) = \inf \{ k \in \mathbb{R}^+ : M(k;\mu) = k \}.
     \label{eq:meandiscovery}
\end{equation}
The existence and uniqueness of $k(\mu)$ when $\mu \neq 0$ follows
from facts to be established in Section \ref{subsec.meanexceed}: 
that
(taking $k$ as real-valued), the function $k \rightarrow
M(k;\mu)/k$ decreases strictly and continuously from a limit of $+ 
\infty$ as $k
\rightarrow 0$ to a limit $< 1$, as $k \rightarrow n$.

A key point is that the mean discovery number is bounded
over the 
parameter spaces $\Theta_n$.
The mean discovery number is monotone in $\mu$: if $|\mu_1|_{(k)} \geq
|\mu_2|_{(k)}$ for all $k$, then $k(\mu_1) > k(\mu_2).$
Thus, on $\ell_0[\eta_n],$ the largest mean discovery number 
$\bar{M}$, say,
is obtained by
taking $k_n = [n \eta_n]$ components to be very large. Writing this out:
\begin{align*}
         \overline{M}_n(k)  = \sup_{\ell_0[\eta_n]} M(k;\mu)
         & = k_n + 2(n-k_n) \tilde{\Phi}(t_k) \\
         & = k_n + (1 - k_n/n) k q_n  \sim k_n + k q_n,
\end{align*}
using the definition of $t_k = z( k q_n / 2 n)$ and $\eta_n \sim k_n/n
\rightarrow 0.$
The first term corresponds to ``true'' discoveries, and the second to
``false'' ones.
Solving $\overline{M}(k) = k$ yields a solution
\begin{equation}
     \tilde{k} = k_n/(1 - (1 - n^{-1}k_n) q_n) \sim k_n/(1 - q_n).
     \label{eq:ktilde}
\end{equation}
In particular, for all $\mu \in \ell_0[\eta_n],$ we have $k(\mu) \leq 
k_n/(1 - q_n (1+ o(1)).$

\medskip

\small
\textbf{Weak $\ell_p$.} \ On $\Theta_n = m_p[\eta_n],$ $E_\mu N(t_k)$
is maximized by taking the components of $\mu$ as large as possible -
i.e. at the coordinatewise upper bound $\bar{\mu}_l = \eta_n (n/l)^{1/p}.$
Thus now
\begin{displaymath}
           \overline{M}_{n}(k) = \sup_{m_n[\eta_n]} E_\mu N(t_k) = 
E_{\bar{\mu}} N(t_k).
\end{displaymath}
To approximate $M(k;\bar{\mu}),$ note first that the summands in
\eqref{eq:mnkmu} are decreasing from nearly $1$ for $\mu_l$ large to
$2 \tilde{\Phi}(t_k)$ when $\mu_l$ is near $0$.
With $k$ held fixed, break the sum into two parts using
$k_n = n \eta_n^p \tau_\eta^{-p}$.
[This choice is explained in more detail after \eqref{eq:kstareqn} below.]
For $l \leq k_n$, the summands are mostly well
approximated by $1,$ and for $l \geq k_n$
predominantly by $2 \tilde{\Phi}(t_k) = q
k /n.$ Since $k_n/n \approx 0,$ we have
\begin{align*}
           M(\nu ;\bar{\mu}) & \approx k_n + (n - k_n) q_n \nu /n \\
                            & \approx n \eta_n^p \tau_\eta^{-p} + q_n \nu .
\end{align*}
Again the first term tracks ``true'' discoveries and the second
``false'' ones.
Solving $M(\nu; \bar \mu) = \nu$ based on this approximation suggests
that, just as in \eqref{eq:ktilde},
\begin{equation}
   \label{eq:kmubar}
   k( \bar \mu) \leq k_n/(1-q_n) (1+o(1)).
\end{equation}
The full proof is given in Section 
\ref{sec:kmubound}.

\textbf{Strong $\ell_p$.}  Since  $\ell_p[\eta_n] \subset 
m_p[\eta_n]$, (\ref{eq:kmubar}) applies here as well.
\normalsize

\subsection{Typical behavior of $\hat{k}_F$ and $\hat{k}_G$}
\label{ssec:thresh-lower-bound}

We turn to the stochastic quantities $\hat k_F$ and $\hat
k_G$. These are defined in terms of the exceedance numbers $N(t_k)$,
which themselves depend on independent (and non-identically
distributed) Bernoulli variables. This suggests the use of
bounds on $\hat k_F$ and $\hat
k_G$ derived from large deviations properties of $N(t_k)$.
Since we are concerned mainly with relatively high thresholds $t_k$,
results appropriate to Poisson regimes are required. Details are in
Section \ref{sec:ldproof}.

To describe the resulting bounds on $[\hat{k}_G, \hat{k}_F],$ we first
introduce some terminology. We say that an event $A_n(\mu)$ is 
$\Theta_n$-likely
  if there exist constants $c_0, c_1$ not depending on
$n$ and $\Theta$ such that
\begin{displaymath}
     \sup_{\mu \in \Theta_n} P_\mu \{ A_n^c (\mu) \} \leq
     c_0 \exp \{ - c_1 \log^2 n \}.
\end{displaymath}

With $\alpha_n$ as in assumption (A), define
\begin{equation} 
\label{eq:kminusdef}
   k_-(\mu) =
   \begin{cases}
     k(\mu) - \alpha_n k_n  &  k(\mu) \geq 2 \alpha_n k_n \\
     0                      &  k(\mu) < 2 \alpha_n k_n ,
   \end{cases}
\end{equation}
and
\begin{equation} \label{eq:kplusdef}
   k_+(\mu) = k(\mu) \vee \alpha_n k_n + \alpha_n k_n.
\end{equation}

\begin{proposition} \label{prop:sandwichbd}
Assume (Q), (H) and (A). 
   For each of the parameter spaces $\Theta_n$, it is
   $\Theta_n-$likely that
\begin{align*}
   k_-(\mu) \leq \hat k_G & \leq \hat k_F \leq k_+(\mu) .
\end{align*}
\end{proposition}

Thus all the penalized 
estimates $\hat{k}_r$ (and any $\hat{k} \in [ 
\hat{k}_G, \hat{k}_F]$)
are with uniformly high probability bracketed between
$k_-(\mu)$ and 
$k_+(\mu)$.
In particular, note that
\begin{equation}
   \label{eq:kpminuskm}
   k_+(\mu) - k_-(\mu) \leq 3 \alpha_n k_n,
\end{equation}
and so the fluctuations in $\hat{k}_r$ are typically small 
compared
to the maximal value over $\Theta_n$.

Here and below it is 
convenient to have a notational variant
for $t_k$, used especially 
when the subscript would be very complicated;
so define
\[
   t[k] = 
z( 2n/kq);
\]
keep in mind that $t$ depends implicitly on $q = q_n$ 
and $n$.
We occasionally use $t_k$ when the subscript is very simple.

Giving this notation its first
workout, the thresholds $\hat{t}_r$ 
are bracketed between
\begin{equation} \label{eq:tplusminusdef}
     t_+(\mu) = t[k_-(\mu)] \quad \mbox{and} \quad t_-(\mu) = t[k_+(\mu)].
\end{equation}
Note that $t_+ > t_-,$ but from \eqref{eq:kpminuskm} and
\eqref{eq:tnubd}, it follows that $t_+/t_- \leq 1 + 3
\alpha_n/t_-^2$.

\subsection{Risk of Hard Thresholding}
\label{sec:hardthreshrisk}

We now study the error of \textit{fixed}
thresholds as a prelude to the study of the data-dependent FDR
thresholds.  We define one-parameter families of configurations and of
thresholds which exhibit key transitional behavior. As might be
expected, these are concentrated around the critical threshold
$\tau_\eta = \sqrt{ 2 \log \eta_n^{-1}}$ corresponding to sparsity
level $\eta_n$.

Consider first a family of (quasi-) least-favorable means 
$\mu_\alpha$.
The coordinates take one of two values, most being 
zero, and
the others amounting to a fraction $\eta$ with value roughly
$\tau_\eta + \alpha.$ Specifically, for $\alpha \in \mathbb{R}$, set
\begin{equation}
   \label{eq:mualphal}
   \mu_{\alpha,l} =
   \begin{cases}
     t[k_n] + \alpha & l \leq k_n , \\
     0    &   k_n < l \leq n.
   \end{cases}
\end{equation}

In a sense $\mu_{0,k_n}$
is right at the FDR 
boundary, while with $\alpha > 0$,
$\mu_{\alpha,k_n}$ is above the 
boundary and
with $\alpha < 0$ is it below the boundary.

Next, consider a `coordinate system' for measuring
the height of 
thresholds in the vicinity of the FDR boundary.
Think of thresholds 
$\{ t : t > 0\}$ as generated by $\{ t[a k_n], a > 0 \}$,
with $a$ fixed while $n$ and $k_n$ increase.
For $a=1$,  we are on 
the
FDR boundary at $k_n$, so that $a < 1$ is above the
boundary and 
$a > 1$ is below the boundary.  The `coordinate'
$a$ will be heavily 
used in what follows.

In fact, these thresholds vary only slowly with $a$: for $a$ fixed, as
$n \rightarrow \infty$,
\begin{equation}
   \label{eq:takkbd}
   | t[a k_n] - t[k_n] | \leq c(a) \tau_\eta^{-1}.
\end{equation}
Nevertheless, the effect of $a$ is visible in the leading term of the
risk:

\begin{proposition}
   \label{prop:hardthreshl0}
Let $\alpha \in \mathbb{R}$ and $a > 0$ be fixed. Let the
configuration $\mu_\alpha \in \ell_0[\eta]$ be defined by
\eqref{eq:mualphal}. For $\ell_r$ loss, the risk of hard thresholding
at $t[ak_n]$ is given, as $n \rightarrow \infty$, by
\begin{equation}
   \label{eq:hardthreshl0}
   \rho( \hat \mu_{H,t[ak_n]}, \mu_\alpha) =
    [ \tilde \Phi (\alpha) + a q_n ]  \cdot k_n \tau_\eta^r  \cdot (1 + o(1)).
\end{equation}
\end{proposition}

Here $k_n \tau_\eta^r$ is asymptotic to the minimax risk for 
$\ell_0[\eta_n]$ \--
compare \eqref{l0mmx} \-- and so defines the benchmark for comparison.

The two leading terms in \eqref{eq:hardthreshl0} reflect false
negatives and false positives respectively. The proof is given in
Section \ref{sec:lowerbds}. Here we aim only to explain how these
terms arise.

The false-negative term $\tilde \Phi (\alpha) k_n \tau_\eta^r$
decreases as $\alpha$ increases. This is natural, as the signals with
mean $m_\alpha = t[k_n] + \alpha$ become easier to detect as $\alpha$
increases -- whatever be the threshold $t[a k_n]$. More precisely,
the $\ell_r$-error due to non-detection, $|y_l| \leq t[a k_n]$, on
each of the $k_n$ terms with mean $m_\alpha$ contributes risk
\begin{equation}
   \label{eq:term20}
   k_n m_\alpha^r P_{m_\alpha}( |y_l| < t[a k_n])
    \sim k_n \tau_\eta^r \Phi( t[a k_n] - m_\alpha),
\end{equation}
since $m_\alpha \sim \tau_\eta$ as $n \rightarrow \infty.$
Finally, \eqref{eq:takkbd} shows that
\begin{equation}
   \label{eq:malpha}
   m_\alpha - t[a k_n] = \alpha + t[k_n] - t[a k_n] = \alpha +
   O(\tau_\eta^{-1}) ,
\end{equation}
so that \eqref{eq:term20} is approximately $\tilde
\Phi(\alpha) k_n \tau_\eta^r $.

The false-positive term shows a relatively subtle
dependence on threshold $a \rightarrow t[a k_n]$. There are $n - k_n$
means that are exactly $0$, and so the risk due to false discoveries
is
\begin{align}
   (n - k_n) E \{ |z|^r : |z| > t[a k_n] \}
    &  \sim 2 n  t[a k_n]^r \tilde \Phi( t[a k_n] )  \label{eq:zerocptrisk}  \\
    & =  a k_n q_n t[a k_n]^r  \label{eq:usetdef} \\
    & \sim a q_n k_n \tau_\eta^r.  \label{eq:useasymp}
\end{align}
[\eqref{eq:zerocptrisk} follows from \eqref{eq:hardriskat0} below,
while \eqref{eq:usetdef} uses the definition of the
FDR boundary, $t[k_n] = \tilde \Phi^{-1}( k_n q_n/ 2n)$, and finally
\eqref{eq:useasymp} follows from \eqref{eq:takntilde} below.]

\bigskip

\textbf{Weak $\ell_p$.}
In this case, we replace the two-point configuration by Winsorized
analogs in the spirit of Section \ref{ssec:wk-lpballs}:
\begin{equation}
   \label{eq:mualpham}
   \mu_{\alpha,l} =\bar \mu_l \wedge m_\alpha, \qquad \qquad
   m_\alpha = t[k_n] + \alpha.
\end{equation}
Now an extra term appears in the risk of hard thresholding when using
thresholds $t[ a k_n]$:
\begin{proposition}
   \label{prop:risk-hard-thresh}
Adopt the setting of Proposition \ref{prop:hardthreshl0}, replacing
only \eqref{eq:mualphal} by \eqref{eq:mualpham}. Then
\begin{displaymath}
   \rho( \hat \mu_{H,t[ a k_n]}, \mu_\alpha) =
     [\tilde \Phi(\alpha) + \frac{p}{r-p} + a q_n ] \cdot k_n 
\tau_\eta^r  \cdot (1 +
     o(1)).
\end{displaymath}
\end{proposition}

The same phenomena as for $\ell_0[\eta]$ apply here, \textit{except}
that the $p/(r-p)$ term arises due to the cumulative effect of missed
detections of means $\bar \mu_l$ that are smaller than $m_\alpha$ but
certainly not $0$. This term decreases as $p$ becomes smaller,
essentially due to the increasingly fast decay of $\bar \mu_l = c_n
l^{-1/p}$. The term disappears in the $p \rightarrow 0$ limit, and we
recover the $\ell_0-$ risk \eqref{eq:hardthreshl0}. This result also
is proved in Section \ref{sec:lowerbds}.

\subsection{FDR on the Least Favorable Family}
\label{sec:fdr-least-favorable}

To track the response of the FDR estimator to members of the family
$\{ \mu_\alpha, \alpha \in \mathbb{R} \}$, we look first at the mean
discovery numbers.  In  Section \ref{sec:lowerbds} we prove:

\begin{proposition}
   \label{prop:kmualpha}
Assume (Q), (H) and (A). Fix $\alpha \in \mathbb{R}$ and define
$\mu_\alpha$ by \eqref{eq:mualphal} and \eqref{eq:mualpham} for
$\ell_0[\eta_n]$ and $m_p[\eta_n]$ respectively. Then as $n
\rightarrow \infty$,
\begin{equation}
   \label{eq:kmualpha}
   k(\mu_\alpha) \sim (1 - q_n)^{-1} \Phi(\alpha) k_n.
\end{equation}
\end{proposition}

Heuristically, for $\ell_0[\eta_n]$, we approximate
\begin{align}
   M(k; \mu_\alpha) & =
    k_n [ \tilde \Phi (t_k - m_\alpha) + \Phi( -t_k - m_\alpha)]
     +2(n-k_n) \tilde \Phi(t_k) \label{eq:mkmualpha}  \\
   & \sim k_n \Phi( m_\alpha - t_k) + q_n k
     \sim k_n \Phi(\alpha) + q_n k,  \notag
\end{align}
from \eqref{eq:malpha}. Solving $M(k;\mu_\alpha) = k$ based on this
approximation leads to \eqref{eq:kmualpha}. The same approach works
for $m_p[\eta_n]$, but with more attention needed to bounding the
component terms in $M(k;\mu_\alpha)$ as detailed in Section
\ref{sec:mean-detect-funct}.

Proposition \ref{prop:kmualpha} suggests that at configuration
$\mu_\alpha$, FDR will choose a threshold close to $t[k(\mu_\alpha)]$,
which is of the form $t[ak_n]$ with $a \sim (1-q)^{-1} \Phi(\alpha)$.
Thus, as $\alpha$ increases, and with it the non-zero components of
$\mu_\alpha$,  the FDR threshold \textit{decreases},
albeit modestly.

The risk incurred by FDR at $\mu = \mu_\alpha$
corresponds to that of hard thresholding at
$t[k(\mu_\alpha)]$. In Section 13 below we prove:
\begin{proposition}
   \label{prop:riskmualpha}
  Assume (Q), (H) and (A). Fix $\alpha \in \mathbb{R}$ and consider
$\mu_\alpha$ defined by \eqref{eq:mualphal} for
$\ell_0[\eta_n]$. Then as $n
\rightarrow \infty$,
\begin{equation}
   \label{eq:riskmualpha0}
   \rho(\hat \mu_F, \mu_\alpha) =
    \Bigl[ \tilde \Phi(\alpha) +  \Phi(\alpha) \frac{q_n}{1 -
    q_n}\Bigr] \cdot k_n \tau_\eta^r  \cdot (1 + o(1)).
\end{equation}
On the other hand, define $\mu_\alpha$ using 
\eqref{eq:mualpham} 
for $m_p[\eta_n]$; 
then
\begin{equation}
   \label{eq:riskmualpham}
   \rho(\hat \mu_F, \mu_\alpha) =
    \Bigl[ \tilde \Phi(\alpha) + \frac{p}{r-p} + \Phi(\alpha) \frac{q_n}{1 -
    q_n}\Bigr] \cdot k_n \tau_\eta^r  \cdot (1 + o(1)).
\end{equation}
\end{proposition}

Formula \eqref{eq:riskmualpha0} shows visibly
the role of the FDR 
control parameter $q$. Note that

\begin{equation} \label{SimpRemark}
 \sup_\alpha  \tilde 
\Phi(\alpha) + \Phi(\alpha) q/(1-q) = \left\{ 
    \begin{array}{ll}
 
1  & q \leq 1/2 \\
          \frac{2q-1}{1-q} & q > 1/2 
 
\end{array} \right.
\end{equation}
Consider the implications 
of this in 
(\ref{eq:riskmualpha0}) in the $\ell_0[\eta_n]$ case.
The minimax risk 
$\sim k_n \tau_\eta^r$, and so the minimax risk
is exceeded asymptotically whenever  $q>1/2$.

We interpret this further.
When $q < 1/2$, the worst configurations in $\{ \mu_\alpha\}$
correspond to $\alpha$ large negative, and yield essentially the
minimax risk.
Indeed, only $\Phi(\alpha)$ of the true non-zero means are
discovered. Each missed mean contributes risk $\sim \mu_\alpha^r\sim
\tau_\eta^r$ and so the risk due to missed means is given roughly by
$\tilde \Phi(\alpha) k_n \tau_\eta^r.$
The risk contribution due to false discoveries, being controlled by
$\Phi(\alpha)$, is negligible in these configurations.

When $q > 1/2$, the worst configurations in $\{ \mu_\alpha \}$
correspond to $\alpha$ large and positive. Essentially all of the
$k_n$ non-zero components are correctly discovered, along with a
fraction $q$ of the $k_n(\mu_\alpha) \sim (1-q)^{-1} \Phi(\alpha) k_n$
which are false discoveries.
In the $\ell_0$ case, the false discoveries
dominate the risk, yielding an error of order
$\Phi(\alpha) [q/(1-q)] k_n \tau_\eta^r.$

When $q = 1/2$, the risk $\rho( \hat \mu_F, \mu_\alpha) \sim k_n
\tau_\eta^r$ \textit{regardless of} $\alpha$, so that all
configurations $\mu_\alpha$ are equally bad, even though the
fraction $\Phi(\alpha)$ of risk due to false discoveries changes from
$0$ to $1$ as $m_\alpha = t[k_n] + \alpha$ increases from values below
$t[k_n]$ through values above.

\subsection{Proof of the Lower 
Bound}
\newcommand{\eps}{\epsilon}

Our interpretation of Proposition 
(\ref{prop:riskmualpha})
in effect gave away the idea for the proof 
of (\ref{eq:lowbound_mp}).  We now fill in 
a few details.

In the 
$\ell_0[\eta_n]$
case, fix $\eps > 0$. Choosing $\alpha = 
\alpha(\eps; q)$ sufficiently large
positive or negative according 
as $q > 1/2$ or $q < 1/2$, 
we get
\[
    \tilde{\Phi}(\alpha(\eps)) 
+ \frac{q}{1-q} \Phi(\alpha(\eps)) > (1-\eps/2) \cdot 
\sup_\alpha   \tilde{\Phi}(\alpha) + \frac{q}{1-q} 
\Phi(\alpha).
\]
\eqref{eq:riskmualpha0}  gives that for large 
$n$,
\[
 \rho(\hat \mu_F, \mu_{\alpha(\eps)})  \geq
       [ 1 + \frac{(2q -1)_+}{1-q}] \cdot k_n \tau_\eta^r  \cdot (1 - 
\eps).
\]
But the family $\mu_\alpha \in \ell_0[\eta_n]$, 
and $R_n( 
\ell_0[\eta_n])  \sim k_n \tau_\eta^r $
hence
\[
    \bar{\rho}(\hat 
\mu_F, \ell_0[\eta_n]) \geq 
     \rho(\hat \mu_F, 
\mu_{\alpha(\eps)}) \geq  
       [ 1 + \frac{(2q -1)_+}{1-q}] \cdot 
R_n( \ell_0[\eta_n])  \cdot (1 - \eps).
\]
As this is true for all 
$\eps > 0$,  the $\ell_0[\eta_n]$ case of
(\ref{eq:lowbound_mp}) 
follows.

For the $\ell_p[\eta_n]$ case, fix $\eps > 0$, and
choose 
again $\alpha = \alpha(\eps;q)$ as in the $\ell_0$ case.
Note that 
$m_\alpha$
is implicitly a function $m_\alpha[k_n]$ 
of the number of 
nonzeros.
Define the pair
$\beta'_n$ and $k'_n$ informally 
as the 
joint solution of
\[
\beta'_n = \eta_n^p (m_\alpha[k'_n])^{-p}, 
\qquad k'_n = n \beta'_n.
\]
(A formal definition can be made using 
the approach in Section 3.2).
Now the mean vector $\mu'_\alpha$ with 
$k'_n$ nonzeros each taking value
$m_\alpha[k'_n]$ 
gives, again by \eqref{eq:riskmualpha0}, that for large 
$n$,
\[
 \rho(\hat \mu_F, \mu_{\alpha})  \geq
       [ 1 + \frac{(2q -1)_+}{1-q}] \cdot k'_n (m_\alpha[k'_n])^r 
\cdot (1 - \eps).
\]
Now from Section 3.2,
$R_n( \ell_p(\eta_n) ) 
\sim n \eta_n^p \tau_{\eta}^{r-p}$, while
\[
  k'_n 
(m_\alpha[k'_n])^r = n \beta'_n  (m_\alpha[k'_n])^{-p} = n  \eta_n^p 
(m_\alpha[k'_n])^{r-p} \sim n \eta_n^p \tau_{\eta}^{r-p}, \quad n 
\goto \infty.
\]
At the same time,  
$\mu_{\alpha} \in 
\ell_p[\eta_n]$.
Hence
\[
    \bar{\rho}(\hat \mu_F, \ell_p[\eta_n]) 
\geq 
     \rho(\hat \mu_F, \mu_{\alpha(\eps)}) \geq  
       [ 1 + 
\frac{(2q -1)_+}{1-q}] \cdot R_n( \ell_p[\eta_n])  \cdot (1 - \eps) 
.
\]
Again this holds for all $\eps > 0$,  and the $\ell_p[\eta_n]$ 
case of
(\ref{eq:lowbound_mp}) follows.

The argument for 
(\ref{eq:lowbound_mp}) in the $m_p[\eta_n]$ case
is entirely 
parallel, only using \eqref{eq:riskmualpham} 
and the modified 
definition of $\mu_\alpha$ for the $m_p$ case.

\subsection{Infrastructure for the Lower Bound}
\label{sec.completeLower}

We look ahead now to 
the arguments supporting
the Propositions of this 
section.

Propositions \ref{prop:hardthreshl0}-\ref{prop:kmualpha} 
will be proved in 
Section \ref{sec:lowerbds} at the very end of the 
paper.
The proofs exploit
viewpoints and estimates set forth 
in
Sections 
\ref{sec:mean-detect-funct}-\ref{sec:threshlems}.
Section 
\ref{sec:mean-detect-funct} sets out properties of
  the mean detection function $\nu \rightarrow M(\nu;\mu)$ of
  \eqref{eq:mnkmu} and its derivatives, with a view to deriving
  information and bounds on the discovery number $k(\mu)$
of Section 5.1.

Section \ref{sec:ldproof} applies these bounds in combination with the
large deviation bounds to prove Proposition \ref{prop:sandwichbd} 
and show that $\hat k_F - \hat k_G \leq 3 \alpha_n k_n$.
Section \ref{sec:threshlems} collects various bounds on the risk of
fixed thresholds, and the risk difference between two data dependent
thresholds. 

All these sections frequently invoke a very useful 
appendix,
Section \ref{sec:gauss-quantiles}, which
assembles needed 
properties of the Gaussian
distribution, of the
quantile function $z(\eta) = \tilde{\Phi}^{-1}(\eta)$ and of
implications of the FDR boundary $t_k$.

\section{The Mean Detection Function}
\label{sec:mean-detect-funct}
\setcounter{equation}{0}

\subsection{Comparing weak $\ell_p$ with $\ell_0$: the effective
   non-zero fraction}
\label{sec:an-analogy-between}


A key feature of $\ell_0[\eta_n]$ is that only $k_n = n \eta_n$
coordinates may be non-zero.
Consequently, the number of `discoveries' at threshold $t[\nu]$ from $n
- k_n$ zero co-ordinates is at most linear in $\nu$ with slope $q_n$:
\begin{equation}
   \label{eq:simplel0}
   (n-k_n) p_\nu(0) = n( 1- \eta_n) \nu q_n / n \leq q_n \nu ,
\end{equation}
since $p_\nu (0) = 2 \tilde \Phi (t[\nu]) =  q_n \nu / n$.

In the case of weak $\ell_p$, the discussion around 
\eqref{supstep}-\eqref{eq:kstareqn} showed that for certain 
purposes, the index $k_n =
n \eta_n^p \tau_\eta^{-p}$ counts the maximum number of
`significantly' non-zero co-ordinates.

In this section, we will see that an alternate, slightly larger index,
$k_n^\prime = n \eta_n^p \tau_\eta^p$ yields for weak $\ell_p$ the
property analogous to \eqref{eq:simplel0}: the number of discoveries
at $t_\nu$ from the $n - k_n^\prime$ smallest means $\mu_l$ is not
essentially larger than $q_n \nu$.
At least for the extremal configuration $( \bar \mu_l )$, the range of
indices $[ k_n, k_n^\prime]$ constitutes a `transition' zone between
`clearly non-zero' means and `effectively zero' ones: this is
discussed further in Section \ref{sec:lipschitz-result}.

To state the result, we need 
a certain error-control function; let 

\begin{equation}
   \label{eq:deltapdef}
   \delta_p(\epsilon) = p \epsilon \int_\epsilon^1 w^{p-2} dw
         \leq
         \begin{cases}
           p |p-1|^{-1} \epsilon^{p \wedge 1}  &  p \neq 1 \\
           ( \log \epsilon^{-1} ) \epsilon   &  p = 1.
         \end{cases}
\end{equation}

\begin{lemma}
   \label{lem:tail-exceed}
Assume hypotheses (Q) and (H).
Let $\tau_\eta^2 = 2 \log \eta_n^{-p}$ and $\epsilon_n = \eta_n
  \tau_\eta$ and $\delta_p(\epsilon)$ be defined as above.
There exists $c = c( b_1,b_3) > 0$ such that for $\nu$ with
$t_\nu^2 \geq 2$, we have,
uniformly in $\mu \in m_p[\eta_n]$ that
   \begin{equation}
     \label{eq:tailex}
     [1 - \epsilon_n^p] q_n \nu \leq
     \sum_{l > k_n^\prime} p_\nu (\mu_l) \leq
              [ 1 + c \delta_p (\epsilon_n) ] q_n \nu.
   \end{equation}
\end{lemma}
The proof is deferred to Section \ref{sec:lipschitz-result} -- compare the
proof of \eqref{eq:mnbd} there.

\medskip
\textit{Remark.} Suppose $q_n \rightarrow q \in [0,1).$
Then for $n$ sufficiently large (i.e. $n$ larger than some $n_0$
depending on $p, q$, and the particular sequence $\eta_n$), it follows
that
\begin{equation}
   \label{eq:qtildedef}
   [1 + c \delta_p(\epsilon_n)] q_n \leq \tilde q_n :=
   \begin{cases}
     (3/2) q_n  & \mbox{if} \ \ q_n \leq 1/2 \\
     (1+q)/2    & \mbox{if} \ \ 1/2 < q_n < 1,
   \end{cases}
\end{equation}
in particular, the right side is strictly less than $1$.

We have just defined $\tilde q_n$ for the case of $m_p$.
If  in the 
nearly-black case
($\ell_0[\eta_n]$) 
we agree to define $\tilde q_n = q_n$ then we may 
write both conclusions
(\ref{eq:simplel0}) and (\ref{eq:tailex}) in one unified form:
\begin{equation}
   \label{eq:tailsum}
   \sum_{l>k_n^\prime} p_\nu (\mu_l) \leq \tilde q_n \nu.
\end{equation}

\textit{The ``true positive'' rate.}
Adopt for a moment the language of diagnostic testing and call
those means with $\mu_l \neq 0$ ``positives'', and those with $\mu_l
=0$ ``negatives''.
In the nearly-black case there are typically $k_n$ positives out of
$n$. In the weak-$\ell_p$ case, there are formally as many as $n$
positives, but as argued above, there are \textit{effectively} $k_n^\prime =
n \eta_n^p t_\eta^p$ positives, and it is this interpretation that we
take here.
If it is assumed (without loss of generality) that $\mu_1 \geq \mu_2
\geq \ldots \geq \mu_n \geq0$, then the {\it true positive rate} using
threshold $t[\nu]$ is defined by
\begin{equation}
   \label{eq:truepos}
   \bar \pi_\nu(\mu) = (1/k_n^\prime) \sum_1^{k_n^\prime} p_\nu (\mu_l).
\end{equation}
In our sparse settings, if the mean discovery number for $\mu$ exceeds
$\nu$, then there is a lower bound on the true positive rate at $\mu$
using threshold $t[\nu]$.

\begin{corollary}
\label{cor:truepos}
Assume (Q) and (H) and define $\tilde q_n$ by \eqref{eq:qtildedef}.
If $n$ is sufficiently large, then uniformly over $\mu$ in $m_p[\eta_n]$ or
   $\ell_o[\eta_n]$ for which $k(\mu) > \nu$,
the true positive rate using threshold $t[\nu]$
   satisfies
   \begin{displaymath}
     \bar \pi_\nu(\mu) \geq (1 - \tilde q_n)(\nu/k_n^\prime).
   \end{displaymath}
\end{corollary}
\begin{proof}
   From the definition of the mean exceedance number, we have
   \begin{displaymath}
     k_n^\prime \bar \pi_\nu(\mu) = M(\nu,\mu) - \sum_{l>k_n^\prime} 
p_\nu(\mu_l).
   \end{displaymath}
Since $\nu < k(\mu)$ we have $M(\nu;\mu) \geq \nu$, and to bound the
sum we use (\ref{eq:tailsum}). Hence
$  k_n^\prime \bar \pi_\nu(\mu) \geq \nu - \tilde q_n \nu, $
as required.
\end{proof}

\subsection{Convexity properties of exceedances}
\label{sec:conv-prop-exce}


The goal of this subsection, Corollary \ref{cor:m1mulower}, shows
that a lower bound on the mean discovery number $k(\mu)$ forces a
lower bound on the mean threshold function $\nu \rightarrow
M(\nu;\mu)$ at least on sparse parameter sets. The idea is to
establish convexity of a certain power function associated with
testing individual components $\mu_l$ and then to use the convexity to
construct two-point configurations providing the needed lower bounds.

Let $N(t_k) = \# \{ l: |y_l|\geq t_k \}$, and as before
$M(k;\mu) = E_\mu N(t_k) = \sum_{\ell=1}^n p_k(\mu_l)$.
Here, the exceedance probability for threshold $t_k$ is given by
\begin{displaymath}
   p_k(\mu) = P\{ |Z+\mu| > t_k \} = \tilde \Phi(t_k - \mu) + \Phi(-t_k
   - \mu),
\end{displaymath}
and we note that as $\mu$ increases from $0$ to $\infty$, $p_k$
increases from $p_k(0)= 2 \tilde \Phi(t_k) = qk/n$ to $p_k(\infty) =
1.$ It has derivative
\begin{displaymath}
   p_k'(\mu) = \phi(t_k - \mu) - \phi(t_k+\mu) > 0, \qquad \qquad
\mu \in (0,\infty).
\end{displaymath}
Since $\mu \rightarrow p_k(\mu)$ is strictly monotone, the inverse
function $\mu_k(\pi) = \mu[\pi;k]
= p_k^{-1}(\pi)$ exists for $qk/n \leq \pi \leq
1$.
In the language of testing, consider the two-sided test of $H_0:
\mu=0$ that rejects at $t_k$. Then $\mu_k(\pi)$ is that alternative $\mu$
at which the test has power $\pi$.
In addition
\begin{displaymath}
   \frac{d}{d \pi} p_k^{-1}(\pi) = \frac{1}{p_k'(p_k^{-1}(\pi))} =
   \frac{1}{p_k'(\mu)}.
\end{displaymath}

\textit{The bi-threshold function.} Given indices $\nu < k$, so that
$t_\nu > t_k$, consider minimizing $M(\nu;\mu)$ over $\mu$ subject to
the constraint that $M(k;\mu)$ stay fixed.
Introduce variables $\pi_l = p_k(\mu_l)$: we wish to minimize
\begin{displaymath}
   M(\nu;\mu) = \sum_l p_\nu(\mu_l) = \sum_l p_\nu( p_k^{-1}(\pi_l))
\qquad \mbox{subject to} \  \sum_l \pi_l = m .
\end{displaymath}
Define a bi-threshold function
\begin{displaymath}
   g(\pi) = g_{\nu, k}(\pi) = p_\nu(p_k^{-1}(\pi)), \qquad \quad
   qk/n \leq \pi \leq 1.
\end{displaymath}
Thus, $g_{\nu, k}(\pi)$ gives the power of the test based on
the 
$t_\nu$-threshold at the
alternative where the $t_k$-threshold has power $\pi$.
As $\nu < k$, 
$g_{\nu,k}(\pi) \leq \pi$.

\begin{lemma}
\label{lem:g-convex}
   If $\nu < k$ then $\pi \rightarrow g_{\nu,k}(\pi)$ is convex and
   increasing from $q \nu /n$ to $1$ for $\pi \in [qk/n,1]$.
\end{lemma}

\begin{proof}
   Setting $\mu = p_k^{-1}(\pi)$, we have
   \begin{displaymath}
     g'(\pi) = \frac{p_\nu'(\mu)}{p_k'(\mu)}
     = \frac{\phi(t_\nu - \mu) - \phi(t_\nu + \mu)}{\phi(t_k - \mu) -
     \phi(t_k + \mu)}
     = e^{(t_k^2 - t_\nu^2)/2} \frac{\sinh t_\nu \mu}{ \sinh t_k \mu} >
     0.
   \end{displaymath}
To complete the proof, we show that if $t > u$, then the function
$G(\mu) = f(t,\mu)/f(u,\mu)$ is increasing for $f(t,\mu) = \sinh
t\mu$. First note that
\begin{displaymath}
   G'(\mu) = G(\mu) \int_u^t D_s D_\mu (\log f)(s, \mu) ds,
\end{displaymath}
and that, on setting $y = 2 s \mu$,
\begin{displaymath}
   D_\mu \log \sinh(s \mu) = s \frac{\cosh s\mu}{\sinh s \mu} =
  \frac{1}{2\mu}[y + \frac{2y}{e^y -1}].
\end{displaymath}
Finally, $D_y \{y  + 2 y/(e^y -1) \}$ has numerator proportional to
$(e^y - y)^2 - 1 - y^2 \geq  y^4/4 > 0.$
\end{proof}

\bigskip
\bigskip

For weak $\ell_p$, define $\tilde q_n$ as in (\ref{eq:qtildedef});
while for the nearly-black case, set $\tilde q_n = q_n.$ Let $a_n$ be
positive constants to be specified. Since \eqref{eq:qtildedef}
guarantees that $\tilde q_n < 1$, define
\begin{displaymath}
   \pi_n = (1 - \tilde q_n) a_n.
\end{displaymath}

\begin{proposition}
   \label{prop:exceed-bound}
Assume (Q) and (H).
As before, let $k_n^\prime = n \eta_n$ (for $\ell_0$) or $n \eta_n^p
\tau_\eta^p$ (for weak $\ell_p$).
Define $\pi_n = a_n (1 - \tilde q_n)$ as above.
Then, uniformly in $\mu$ for which $k(\mu) \geq a_n k_n^\prime$, we have

\noindent (a) for $\nu \leq a_n k_n^\prime$,
\begin{equation}
   \label{eq:mnulower}
   M(\nu; \mu) \geq k_n^\prime \bar \Phi( t_\nu - \mu[\pi_n; a_n k_n^\prime] ).
\end{equation}
(b) In particular, for $\nu = 1$, and $a_n \geq b_5( \log n)^{-r},$
\begin{equation}
   \label{eq:m1logbd}
   M(1;\mu) \geq c ( \log n )^{\gamma - r - 1}.
\end{equation}
\end{proposition}
\textit{Remarks.} 1. The lower bound of (\ref{eq:m1logbd}) is valid
for all sparsities $\eta_n^p$ in the range
$[ n^{-1} \log^\gamma n, n^{-b_3}]$;
clearly it is far from sharp for $\eta_n^p$ away from the lower limit.
If needed, better bounds for specific cases would follow from
(\ref{eq:knvnsn}) and (\ref{eq:basicbd}) in the proof below.

2. We shall need only the lower bound for $\nu = 1$, but the methods
    used below would equally lead to bounds for larger $\nu$, working
    from the intermediate estimate (\ref{eq:mnulower}).

\begin{corollary}
\label{cor:m1mulower}
Let $k_n = n \eta_n^p \tau_\eta^{-p}$, and
$\alpha_n$ as in assumption (A);  then uniformly in $\mu$ for
which $k(\mu) \geq \alpha_n k_n$,
\begin{displaymath}
   M(1;\mu) \geq c ( \log n)^{\gamma - p - 3/2}.
\end{displaymath}
\end{corollary}

\begin{proof}
For convenience, we abbreviate $a_n k_n^\prime$ by $\kappa$. 
The bithreshold function $g = g_{\nu,\kappa}$ is convex (Lemma
   \ref{lem:g-convex}), and so
   \begin{displaymath}
     M(\nu;\mu) = \sum_{l=1}^n g_{\nu,\kappa}(\pi_l)
                \geq \sum_{l=1}^{k_n^\prime} g(\pi_l)
                \geq k_n^\prime g (\bar \pi_\kappa(\mu) ),
   \end{displaymath}
where $\bar \pi_\kappa(\mu) = (1/k_n^\prime) \sum_1^{k_n^\prime} \pi_l$ is
the true positive rate defined at \eqref{eq:truepos}.
Since $k(\mu) \geq \kappa$,
Corollary \ref{cor:truepos} bounds $\bar \pi_\kappa(\mu) \geq (1 - \tilde
q_n) (\kappa/k_n^\prime) = \pi_n$ and so from the monotonicity of $g$,
\begin{displaymath}
   g(\bar \pi_\kappa(\mu) ) \geq g(\pi_n) = p_\nu (\mu_\kappa(\pi_n)) \geq
   \bar \Phi (t_\nu - \mu_\kappa(\pi_n)).
\end{displaymath}
This establishes part (a). For (b), we seek an upper bound for
\begin{equation}
   \label{eq:bd1}
   t_\nu - \mu_\kappa(\pi_n) =
   t_\nu - t_\kappa + t_\kappa -  \mu_\kappa,
\end{equation}
where we abbreviate $\mu_\kappa(\pi_n)$ by $\mu_\kappa$ for convenience.
First note that
\begin{displaymath}
   \pi_n = p_\kappa( \mu_\kappa) =
   \Phi( -t_\kappa -  \mu_\kappa) + \bar \Phi (t_\kappa -  \mu_\kappa),
\end{displaymath}
which shows that $\bar \Phi (t_\kappa -  \mu_\kappa) \leq \pi_n \leq 2
\bar \Phi( t_\kappa -  \mu_\kappa)$
from which we get
\begin{displaymath}
   t_\kappa -  \mu_\kappa \leq z( \pi_n/2) \leq
   \sqrt{2 \log( 2 \pi_n^{-1})} ,
\end{displaymath}
using (\ref{eq:zetabd}). Since $\pi_n =  (1-\tilde q_n) a_n \geq
c_3/ (\log n )^{-r}$, we conclude that
\begin{equation}
   \label{eq:bd2}
   t_\kappa -  \mu_\kappa(\pi_n) \leq \sqrt{ 2r \log \log n} + c_4.
\end{equation}

 From (\ref{eq:tnubd}), we have
\begin{equation}
   \label{eq:bd3}
     0 \leq t_\nu - t_\kappa \leq \sqrt{2 \log (n/\nu)} - \sqrt{2 \log
   (n/\kappa)} + c(b_1,b_3).
\end{equation}
Since $\kappa = a_n k_n^\prime$ with $a_n \leq 1$, the right side 
only increases
if we replace $\kappa$ by $k_n^\prime$.

At this point, we specialize to the case $\nu = 1$ and set
$v_n = \sqrt{ 2 \log n} - \sqrt{ 2 \log n/k_n^\prime}$.
Combining (\ref{eq:bd1}), (\ref{eq:bd3}) and (\ref{eq:bd2}), we find that
\begin{displaymath}
   t_\nu - \mu_\kappa(\pi_n) \leq v_n + c(b_1,b_3) + \sqrt{ 2r \log 
\log n} + c_4.
\end{displaymath}
For $n \geq n(b),$ the last three terms are bounded by $s_n = \sqrt{
   (2r + 1)  \log \log n}$. So, from (\ref{eq:mnulower}),
\begin{equation}
     \label{eq:knvnsn}
M(1;\mu) \geq k_n^\prime \bar \Phi (v_n+s_n).
\end{equation}

We may rewrite $k_n^\prime$ in terms of $v_n$, obtaining
\begin{displaymath}
   \log k_n^\prime = v_n \sqrt{ 2 \log n} - v_n^2/2.
\end{displaymath}
The bound $\bar \Phi (w) \geq \phi(w)/(2w) $ holds for $w > \sqrt 2$;
applying this we conclude
\begin{equation}
   \label{eq:basicbd}
   k_n^\prime \bar \Phi (v_n+s_n) \geq \frac{e^{-s_n^2/2}}{2(v_n+s_n)} 
\exp \{ v_n ( \sqrt{
   2 \log n} - v_n -s_n ) \}.
\end{equation}
Since $e^{-s_n^2/2} = (\log n)^{-r - \hf}$ and $v_n \leq \sqrt{ 2 \log n}$, the
first factor is bounded below by $c_0 ( \log n)^{-r-1}.$ To bound the
main exponential term, set $g(v) = v ( \sqrt{ 2 \log n} - v -s_n ).$
We note that $n \eta_n^p \in [ \log^\gamma n, n^{1-b_3}]$ and so
$\tau_\eta^2 = 2 \log \eta_n^{-p} \leq 2 \log n$ and so
$\tau_\eta^p \in [1, 2 \log n]$ and so $k_n^\prime \in [ \log^\gamma n, (2
\log n) n^{1-b_3}]$.
For $\ell_0[\eta_n]$, $k_n^\prime \in [ \log^\gamma n, n^{1-b_3}]$; we
shall see shortly that the difference between the two cases doesn't
matter here.

We now estimate the values of $v$ and $g(v)$ corresponding to these
bounds on $k_n^\prime$.  At the lower end, $k_n^\prime = \log^\gamma 
n$, then (using
$\sqrt a - \sqrt{a - \epsilon} \geq \epsilon/(2 \sqrt a)) $,
\begin{displaymath}
   v_n
   \geq \frac{\gamma \log \log n}{\sqrt{ 2 \log n}} =: v_{1n}
\end{displaymath}
and one checks that $g(v_{1n}) \geq \gamma \log \log n - 1$ for $n \geq
n(b)$.
At the upper end, if $k_n^\prime = (2 \log n) n^{1 - b_3}$, then
\begin{displaymath}
   v_n \leq ( 1 - \sqrt b_3 ) \sqrt{ 2 \log n} + c =: v_{2n}
\end{displaymath}
and one checks that $g(v_{2n}) = ( \sqrt b_3 - b_3) ( 2 \log n
) ( 1 + o(1)) \geq g(v_{1n})$
for $n$ large.

Since $g(v)$ is a concave quadratic polynomial with maximum between
$v_{1n}$ and $v_{2n}$, it follows that for $k_n^\prime$ in the range
indicated above,
\begin{displaymath}
    e^{g(v_n)}  \geq e^{g(v_{1n})} \geq e^{-1} ( \log n)^\gamma.
\end{displaymath}
Combined with the bound on the first factor in (\ref{eq:basicbd}), we
get
\begin{displaymath}
   k_n^\prime \tilde \Phi (v_n+s_n) \geq c_0 e^{-1} ( \log n)^{\gamma -
   r -1},
\end{displaymath}
as required for part (b).

For part (c), we define $a_n$ by the equation $a k_n^\prime = \alpha_n
k_n$, so that $a_n = \alpha_n \tau_\eta^{-2p} \geq \alpha_n (2 \log
n)^{-p}$. Since $\alpha_n \geq b_4 (\log n)^{-1/2}$, we apply part (b)
with $r = p+ 1/2$.
\end{proof}

\subsection{Properties of the Mean Detection Function}\label{subsec.meanexceed}

This subsection collects some properties of
\begin{displaymath}
   M(\nu;\mu) = \sum_l \tilde \Phi(t_\nu - \mu_l) + \Phi( -t_\nu - \mu_l)
\end{displaymath}
as a function of $\nu$, considered as a real variable in $\mathbb{R}_+$.
Writing $\dot{M}, \ddot{M}$ for partial derivatives w.r.t. $k$,
calculus shows
that
\begin{align}
         \partial t_{\nu} / \partial \nu & = - q_n/ (2n \phi(t_{\nu})) 
\notag \\
         \dot{M}(\nu;\mu) & = (- \partial t_\nu/ \partial \nu) \sum_l
   \phi(t_\nu - \mu_l) + \phi(t_\nu + \mu_l) \label{eq:mdot1}\\
          & = (q_n/n) \sum_{l} e^{-\mu_{l}^{2}/2}
         \cosh (t_{\nu} \mu_{l}) > 0 \label{eq:mdot} \\
     \ddot{M}(\nu;\mu)  & =
         - q_n^{2} / (2 n^{2} \phi(t_{\nu})) \sum_{l} \mu_{l}
         e^{-\mu_{l}^{2}/2}
         \sinh(t_{\nu} \mu_{l} ) \leq 0,   \label{pkdd}
\end{align}
with strict inequality unless $\mu \equiv 0.$
Finally, since $M(0;\mu) = 0,$ there exists $\tilde{\nu} \in [0,\nu]$
such that the threshold exceedance function
$\nu^{-1} M(\nu;\mu) = \nu^{-1} (M(\nu;\mu) - M(0;\mu)) =
\dot{M}(\tilde{\nu},\mu),$ and hence, for each $\mu$, the exceedance
function is decreasing in $\nu$:
\begin{equation}
         \frac{\partial}{\partial \nu} \Bigl( \frac{M(\nu;\mu)}{\nu} 
\Bigr) = \frac{1}{\nu}
         [\dot{M}(\nu;\mu) - \dot{M}(\tilde{\nu},\mu)] \leq 0.
         \label{eq:pkoverk}
\end{equation}

Let us focus now on $\ell_0[\eta_n]$. In this case
\begin{displaymath}
   M(\nu;\mu) = \sum_1^{k_n} [\tilde \Phi(t_\nu - \mu_l) + \Phi( -t_\nu
   - \mu_l)] + (1 - \eta_n) q_n \nu,
\end{displaymath}
and so, using \eqref{eq:mdot1} and \eqref{eq:dtnudnu},
\begin{displaymath}
   \dot M(\nu;\mu) \leq \frac{2 \phi(0) k_n}{\nu t_\nu} + (1 - \eta_n) q_n.
\end{displaymath}
In particular, if $\nu = a k_n$, then
\begin{displaymath}
   \dot M(a k_n; \mu) \leq \frac{2 \phi(0) }{a t[a k_n]} + q_n.
\end{displaymath}
Finally, if $a = \alpha_n = 1/(b_4 \tau_\eta)$, then \eqref{eq:takn}
shows that for $\eta_n$ sufficiently small,
\begin{equation}
   \label{eq:b4bd}
   \frac{2 \phi(0) }{a t[a k_n]} = b_4 \sqrt{2/\pi}
   \frac{\tau_\eta}{t[\alpha_n k_n]} \leq b_4.
\end{equation}
As a result, uniformly in $\ell_0[\eta_n]$,
\begin{equation}
   \label{eq:l0mdotbd}
   \dot M(\alpha_n k_n; \mu) \leq b_4 + q_n < q' ,
\end{equation}
by the definition of $b_4$; recall assumptions (A) and (Q) of Section \ref{sec:GenAssump}. \qed

\subsection{Weak $\ell_p$: Bounds for the detection function}
\label{sec:lipschitz-result}

For weak $\ell_p$, we do not have such a simple bound on $\dot M$ as
\eqref{eq:l0mdotbd}.
 From the preceding calculations, we know that $\nu \rightarrow \dot
M(\nu;\mu)$ is positive and decreasing. We will need now some sharper
estimates, uniform over $m_p[\eta_n]$ (and $\ell_0[\eta_n]$)
in the scaling $\nu = a k_n,$ with $a$ regarded as variable.
This will lead to bounds on the solution of $M(a k_n ; \mu) = a k_n$
and hence to bounds on $k(\mu)$ (cf. Corollary \ref{cor:kmubound}).


The two key phenomena:

(a) If $a_1$ is fixed, then for $\nu$ in intervals $[ a_1
k_n, a_1^{-1} k_n]$, the slope of $M$ is, for large $n$, essentially
constant and equal to $q_n$. This reflects exclusively the effect of
false detections.

(b) For small $a$ ($\sim 1/ \tau_\eta$, say), the order of
magnitude of $\dot M(a k_n; \mu)$ can be as large as $1/( a
t[ak_n])$. This reflects essentially the effect of true detections.

Since $\mu \rightarrow \dot M(k;\mu)$ is an even function of each
$\mu_l$, we may assume without loss of generality that $\mu_l \geq 0$
for each $l$.

To bound $\dot M$, divide the range of
summation into three regions, defined by the indices
\begin{alignat*}{2}
k_n & = n \eta_n^p \tau_\eta^{-p},  & \qquad
     \bar \mu_{k_n} & = \tau_\eta \\
k_n^{\prime} & = n \eta_n^p \tau_\eta^{p},  & \qquad
     \bar \mu_{k_n^\prime} & = \tau_\eta^{-1}.
\end{alignat*}
Thus, we write
\begin{displaymath}
   \dot M = \dot M_{pos} + \dot M_{trn} + \dot M_{neg},
\end{displaymath}
where the sum in $\dot M_{pos}$ extends over the range $[1,k_n]$ of true
``positives''. The sum in $\dot M_{trn}$ ranges over $(k_n, k_n^\prime]$ and
is ``transitional'', while the sum for $\dot M_{neg}$ ranges over
$(k_n^\prime,n]$ and corresponds to means that are essentially true
``negatives''.

A rough statement of the results to follow is that for $a$ in the
range $[\gamma \tau_\eta^{-1}, 1]$,
\begin{align*}
   \sup_{m_{p}[\eta_n]} \dot M_{pos} (a k_n; \mu) & \asymp \frac{1}{a
   t[ak_n]} \\
   \sup_{m_p[\eta_n]} \dot M_{trn} (a k_n; \mu) & =    O \Bigl( \frac{1}{a
   t^2[ak_n]} \Bigr) \\
                      \dot M_{neg} (a k_n; \mu) &\sim q_n.
\end{align*}

Combining these will establish:

\begin{proposition}
   \label{prop:lipschitz}
Assume (Q), (H) and (A).
For $n > n(b)$, $t_\nu \geq c_0, a \geq 1$ and all $\mu \in
m_p[\eta_n]$,
\begin{equation}
   \label{eq:lpderivbd}
   q_n [1 - \epsilon_n^p] \leq
   \dot M(a k_n; \mu) \leq q_n[1 + c(b) \delta_p(\epsilon_n)] +
      \frac{2 \phi(0)}{a t[ak_n]} \bigl[ 1 + \frac{c_0}{t[a k_n]} \bigr].
\end{equation}
If, in addition, $a \geq \gamma \tau_\eta^{-1},$ then for $\eta^p \leq
\eta(\gamma,p,b_1,b_3)$ sufficiently small
\begin{equation}
   \label{eq:uniflowermp}
   \sup_{m_p[\eta_n]} \dot M (a k_n; \mu)
    \geq q_n [1 - \epsilon_n^p] + \frac{c_0}{a t[ak_n]}.
\end{equation}
\end{proposition}

The proof consists in building estimates for $\dot M$ in the positive,
transition and negative zones. It will also be convenient, for
Corollary \ref{cor:kmubound} below, to obtain estimates at the same time for
the corresponding components of the detection function $M$ itself.

\subsubsection{Positive zone} 
$M_{pos}(\nu;\mu) = \sum_{l=1}^{k_n} 
p_\nu(\mu_l)$
is, for $\mu = \bar \mu_l$, approximately constant
on the interval $[a_1 k_n,
a_1^{-1} k_n]$: for $a_1 \leq a \leq a_1^{-1}$,
\begin{equation}
   \label{eq:sandwich-pos}
   M_{pos}(a k_n; \bar \mu) \in [ 1 - \epsilon_{1n}, 1] k_n.
\end{equation}
\begin{proof}
   The upper bound follows from $p_\nu(\mu_l)\leq 1.$
For the lower bound, since $\tilde \Phi(t_\nu - \bar \mu_l)$ is
decreasing in $l$ and increasing in $\nu$, we have for any $l_1 \leq
k_n$,
\begin{displaymath}
   M_{pos}(\nu;\bar \mu) \geq \sum_1^{k_n} \tilde \Phi (t_\nu - \bar \mu_l)
     \geq l_1 \tilde \Phi(t_\nu - \bar \mu_l)
     \geq l_1 \tilde \Phi( t[\delta k_n] - \bar \mu_{l_1} ).
\end{displaymath}
Choose $\gamma_n = \tau_\eta^{-1}$ and define $l_1$ by the equation
\begin{displaymath}
   \bar \mu_{l_1} = t[\delta k_n] + z(\gamma_n).
\end{displaymath}
 From \eqref{eq:takn} it is clear that $\bar \mu_{l_1} \geq \tau_\eta$
so that $l_1 \leq k_n$. Hence
\begin{equation}
   \label{eq:l1gamma}
   M_{pos}(\nu;\bar \mu) \geq l_1 \tilde \Phi( - z(\gamma_n))
                 = l_1 (1 - \gamma_n).
\end{equation}
We have from \eqref{eq:takn1} that $\bar \mu_l =
\tau_\eta + \sqrt{ 2 \log \tau_\eta} + c(b_1,b_3)$ for $\eta$ small,
and hence
\begin{displaymath}
   \frac{l_1}{k_n} = \Bigl[ \frac{\tau_\eta}{t[\delta k_n] +
   z(\gamma_n)} \Bigr]^p
   \geq \Bigl[ 1 + \frac{\sqrt{2 \log \tau_\eta} + c}{\tau_\eta}
   \Bigr]^{-p}
   \geq 1 - c \tau_\eta^{-1} \sqrt{ 2 \log \tau_\eta}.
\end{displaymath}
Since $\gamma_n = \tau_\eta^{-1}$, the last two displays imply
\eqref{eq:sandwich-pos} with $\epsilon_{1n} = c
\tau_\eta^{-1} \sqrt{ 2 \log \tau_\eta}$.
\end{proof}

{\it Bound for $\dot M_{pos}$.}
Since $M_{pos}(\nu;\mu) = \sum_1^{k_n} \tilde \Phi(t_\nu - \mu_l)
+ \Phi(-t_\nu -\mu_l)$, we
have
\begin{equation}
   \label{eq:mpdotnu}
   \dot M_{pos}(\nu;\mu) = (- \partial t_\nu/ \partial \nu) \sum_1^{k_n}
   \phi(t_\nu - \mu_l) + \phi(t_\nu + \mu_l)
\end{equation}
 From \eqref{eq:dtnudnu}, we obtain $\dot M_{pos}(\nu;\mu) \leq 2 k_n \phi(0)
/ (\nu t_\nu)$. Hence, for $\nu = a k_n$ with $t_\nu \geq 1$ and all
$\mu$,
\begin{equation}
   \label{eq:mpdotakn}
   \dot M_{pos} (ak_n ; \mu) \leq \frac{2 \phi(0)}{a t[ak_n]}.
\end{equation}
Turning to the lower bound, we note that $\bar \mu_l \geq t_\nu$ if
and only if $l \leq n \eta_n^p t_\nu^{-p} =: l_\nu$.
By setting all $\mu_l = t_\nu$ for $l \leq l_\nu$, we find from
\eqref{eq:mpdotnu} and \eqref{eq:dtnudnu} that
\begin{displaymath}
   \sup_{m_p[\eta_n]} \dot M_{pos}(\nu;\mu) \geq \frac{\phi(0)}{2}
   \frac{k_n \wedge l_\nu}{\nu t_\nu}.
\end{displaymath}
(This also holds for $\ell_0[\eta_n]$.)
If $\nu = a k_n,$ then
\begin{displaymath}
   \frac{k_n \wedge l_\nu}{\nu t_\nu} =
   \frac{1}{a t[ak_n]} \Bigl[ \Bigl( \frac{\tau_\eta}{t[ak_n]} \Bigr)^p
   \wedge 1 \Bigr].
\end{displaymath}
If $a \geq \gamma \tau_\eta^{-1}$, then for $\eta$ sufficiently small,
\eqref{eq:takn2} says $\tau_\eta/t[ak_n] \geq \hf$.
Combining the last remarks, we conclude that for $\nu = a k_n$ and $a
\geq \gamma \tau_\eta^{-1}$, then for $\eta$ sufficiently small,
\begin{equation}
   \label{eq:lowermpdot}
   \sup_{m_p[\eta_n]} \dot M_{pos}( ak_n;\mu) \geq \frac{c_0}{a t[ak_n]}.
\end{equation}
Here $c_0$ denotes an absolute constant.

\subsubsection{Transition zone} For $a \leq a_1^{-1}$ and $\eta \leq
\eta_0$ sufficiently small, we have, uniformly in $m_p[\eta_n]$:
\begin{align}
   \label{eq:trans-m}
   0 \leq & M_{trn}( a k_n; \mu) \leq c_0 \tau_\eta^{-1} k_n
   \qquad \qquad \mbox{and} \\
          & \dot M_{trn} (a k_n; \mu) \leq
      \frac{c_0}{a t^2[a k_n]}.  \label{eq:trans-mdot}
\end{align}
To bound $\dot M_{trn} (k;\mu)$, introduce
\begin{displaymath}
   h(\mu,t) = e^{t\mu -\mu^2/2}
\end{displaymath}
which increases from $1$ to a global maximum of $e^{t^2/2} =
\phi(0)/ \phi(t)$ as $\mu$ grows from $0$ to $t$.
We have from \eqref{eq:mdot} that $\dot M(k;\mu) \leq 2 (q_n/n) \sum
h(\mu_l, t_k)$.
The arguments for the two bounds run in parallel. We have
\begin{displaymath}
   M_{trn}(\nu;\mu) \leq \sum_{k_n}^{k_n^\prime} H_1(\mu_l) \qquad
   \mbox{and} \qquad
   \dot M_{trn}(\nu;\mu) \leq (q_n/n) \sum_{k_n}^{k_n^\prime} H_2(\mu_l),
\end{displaymath}
where $H_1(\mu) = 2 \tilde \Phi( t_\nu - \mu)$ and
$H_2(\mu) = h( \mu \wedge t_\nu, t_\nu)$ are both increasing functions
of $\mu \geq 0$. By integral approximation,
\begin{displaymath}
   \sum_{l = k_n}^{k_n^\prime} H(\mu_l) \leq
   \sum_{k_n}^{k_n^\prime} H_1(\bar \mu_l) \leq
   n \eta_n^p \int_{\tau_\eta^{-1}}^{\tau_\eta} H(u) u^{-p-1} du.
\end{displaymath}
For $a \leq a_1^{-1}$, we have $t = t_\nu \geq \tau_\eta - 3/2$ for
$\eta$ sufficiently small by \eqref{eq:takn}, while \eqref{eq:takn2}
shows that $\tau_\eta^{-1} \geq 1/(2 t_\nu)$. Let $\bar H = \sup H$;
we have
\begin{displaymath}
   \int_{\tau_\eta^{-1}}^{\tau_\eta} H(u) u^{-p-1} du
   \leq \int_{1/(2t)}^t H(u) u^{-p-1} du  + (3/2) \bar H t^{-p-1}
   \leq c_0 \bar H t^{-p-1}
\end{displaymath}
after using Lemma \ref{lem:integralbound} below.
For $M_{trn}$, $\bar H = 2$ and we have from the previous displays
\begin{displaymath}
   M_{trn}( a k_n; \mu) \leq c_0 n \eta_n t[ak_n]^{-p-1} \leq c_0
   \tau_\eta^{-1} k_n,
\end{displaymath}
since $t[a k_n]^{-1} \leq 2 \tau_\eta^{-1}.$

For $\dot M_{trn}$, $\bar H = h(t_\nu, t_\nu) = \phi(0) / \phi(t_\nu)$,
and since $\phi(t_\nu) \geq \hf t_\nu \tilde \Phi(t_\nu)$,
\begin{displaymath}
     \dot M_{trn}(a k_n;\mu) \leq q_n \eta_n^p .
    \frac{c_0 h(t_\nu,t_\nu)}{t_\nu^{p+1}}
   \leq \frac{c_0 q_n  \eta_n^p}{ t_\nu^{p+2} \tilde \Phi(t_\nu)}
   \leq c_0 \Bigl( \frac{\tau_\eta}{t[a k_n]}
   \Bigr)^p \frac{1}{a t^2[a k_n]},
\end{displaymath}
and so \eqref{eq:trans-mdot} follows from \eqref{eq:takn}.

\begin{lemma}
   \label{lem:integralbound} For $p \leq 2$ and $t \geq 2$,
and for $h(u,t)$ given by either $\tilde \Phi(t-u)$ or
$\phi(t-u)/\phi(t)$, there is an absolute constant $c_0$ such that
   \begin{displaymath}
      \int_{1/(2t)}^t \frac{h(u,t) \, du}{u^{p+1}} \leq
     \frac{c_0 h(t,t)}{t^{p+1}}.
   \end{displaymath}
\end{lemma}
\begin{proof}
Writing $v$ for $t-u$, we find that in the two cases $h(u,t)/h(t,t)$
equals $2 \tilde \Phi(v)$ or $e^{-v^2/2}$ respectively. By
\eqref{eq:global}, $2 \tilde \Phi(v) \leq e^{-v^2/2}$ for $v \geq 0$,
so  in either case the integral
in question is bounded by
   \begin{displaymath}
     \frac{h(t,t)}{t^{p+1}} \int_0^{t-1/(2t)} \exp \{ -v^2/2 + (p+1) g(v) \} dv,
   \end{displaymath}
   where the convex function $g(v) = \log t - \log (t-v)$ is bounded
   for $0 < v < t-1$ by $4 v (\log t)/t$. Completing the square in
   the exponent, gives an integrand smaller than $\sqrt{2 \pi}$
   times a unit-variance Gaussian density centered at $\mu(p,t) = 4 (p+1)
   (\log t)/t$.  Since $\mu(p,t) \leq c_0$ for $p\leq 2$ and $t
   \geq 2$, the previous integral is bounded by $\sqrt{2 \pi} \exp \{
   \hf \mu^2(p,t) \} \leq c_0$.
\end{proof}

\subsubsection{Negative zone} 
Under conditions  described immediately below,
\begin{align}
   [1 - \epsilon_n^p] q_n \nu & \leq M_{neg}(\nu; \mu) \leq
   [1 + c \delta_p(\epsilon_n)] q_n \nu , \label{eq:mnbd} \\
   [1 - \epsilon_n^p] q_n  & \leq \dot M_{neg}(\nu; \mu) \leq
   [1 + c \delta_p(\epsilon_n)] q_n .  \label{eq:mndotbd}
\end{align}
The lower bound in \eqref{eq:mnbd} holds for all $n,\nu,\mu$. All
other bounds require $\mu \in m_p[\eta_n]$ and $n \geq n(b)$.
The upper bounds further require $\nu$ such that $t_\nu \geq c(b)$.

\bigskip

\textit{Lower Bounds.} \ For $M_{neg}(\nu;\mu) =\sum_{l \geq k'_n} 
p_\nu(\mu_l)$
this is simple because $p_\nu(\mu_l) \geq p_\nu(0) = q_n \nu/n$, so
that $M_{neg}(\nu;\mu) \geq (n - k_n^\prime) q_n \nu/n = [1 - \epsilon_n^p] q_n
\nu$.

For $\dot M_{neg}(\nu; \mu)$, we set $f(\mu,t) = e^{-\mu^2/2} \cosh
(t\mu)$ and check that for given $c_1$, $\mu \rightarrow f(\mu,t)$ is
increasing for $t \mu \leq c_1$ and $t \geq \sqrt{c_1}$.
For $l \geq k_n^\prime$ we have
\begin{equation}
   \label{eq:munegbd}
   \mu_l \leq \bar \mu_l \leq \bar \mu_{k_n^\prime} = \tau_\eta^{-1}
   \leq c_1 t_1^{-1} \leq c_1 t_\nu^{-1} ,
\end{equation}
by \eqref{eq:t1taueta}, for $c_1 = c_1(b)$. Consequently $f(\mu_l,
t_\nu) \geq f(0, t_\nu) = 1$ and so
\begin{displaymath}
   \dot M_{neg}(\nu;\mu) = (q_n/n)\sum_{l \geq k'_n} f(\mu_l,t_\nu)
   \geq (q_n/n) (n- k_n^\prime) = [1 - \epsilon_n^p] q_n.
\end{displaymath}

\textit{Upper Bounds.} \ The arguments run in parallel: we have
\begin{equation}
   \label{eq:bothbds}
     M_{neg}(\nu;\mu) \leq\sum_{l \geq k'_n} H_1(\mu_l) \qquad
   \mbox{and} \qquad
   \dot M_{neg}(\nu;\mu) \leq (q_n/n)\sum_{l \geq k'_n} H_2(\mu_l),
\end{equation}
where $H_1(\mu) = 2 \tilde \Phi(t_\nu - \mu)$ and $H_2(\mu) =
h(\mu,t_\nu)$ are both increasing and convex functions of $\mu \in
[0,1]$, at least when $t>2.$
Using \eqref{eq:munegbd} along with $t = t_\nu$, this convexity
implies
\begin{displaymath}
   H(\mu_l) \leq ( 1 - t \mu_l/c_1) H(0) + (t \mu_l/c_1) H(c_1/t)
\end{displaymath}
and hence, since $t \leq c_b \tau_\eta$ from \eqref{eq:t1taueta},
\begin{equation}
   \label{eq:convexbd}
   n^{-1}\sum_{l \geq k'_n} H(\mu_l)
   \leq H(0) + c_b c_1^{-1} H(c_1/t) \cdot \tau_\eta n^{-1}\sum_{l \geq k'_n}
   \bar \mu_l.
\end{equation}
By an integral approximation, since $\epsilon_n = \eta_n \tau_\eta$
and $k_n^\prime/ n = \epsilon_n^p$, we find
\begin{equation}
   \label{eq:integappr}
     \tau_\eta n^{-1} \sum_{l=k_n^\prime}^n \bar \mu_l \leq \tau_{\eta}
            \eta_n \int_{k_n^\prime/n}^1  x^{-1/p} dx
   =  p \epsilon_n \int_{\epsilon_n}^1 s^{p-2} ds
   =  \delta_p(\epsilon_n).
\end{equation}
To apply these bounds to $M_{neg}(\nu,\mu)$, we note that $H_1(0) = q_n
\nu/n$ while from \eqref{eq:increment}, for $t \geq \sqrt{2 c_1}$,
\begin{displaymath}
   H(c_1 t^{-1}) \leq 8 e^{c_1} \tilde \Phi(t) = 4 e^{c_1} (q_n \nu/n).
\end{displaymath}
 From \eqref{eq:bothbds}, \eqref{eq:convexbd} and \eqref{eq:integappr},
we obtain \eqref{eq:mnbd}.

Turning to $\dot M_{neg}(\nu;\mu)$, we note that $H_2(0) = 1$ and
$H_2(c_1/t) = \exp \{ c_1 - c_1^2/(2 t^2) \} \leq e^{c_1}$ and so the
same bounds combine to yield \eqref{eq:mndotbd}.


\subsubsection{Conclusion} 
\label{sec:kmubound}
The upper bound in 
\eqref{eq:lpderivbd} follows
by combining those in \eqref{eq:mpdotakn}, \eqref{eq:trans-mdot} and
\eqref{eq:mndotbd}. The lower bound \eqref{eq:uniflowermp} follows by combining
\eqref{eq:lowermpdot} and \eqref{eq:mndotbd}.

\begin{corollary}
\label{cor:kmubound}
   Let $d_n = 2 c_0 \tau_\eta^{-1}$ (where $c_0$ is the constant
   in \eqref{eq:trans-m}.) Uniformly in $m_p[\eta_n]$,
   \begin{equation} \label{eq:kmax}
     k(\mu) \leq (1 - q_n - d_n)^{-1} k_n.
   \end{equation}
\end{corollary}
\begin{proof}
   Let $s = (1 - q_n - d_n)^{-1}$. Combining the bounds on $M_{pos},
   M_{trn}$ and $M_{neg}$ in \eqref{eq:sandwich-pos}, \eqref{eq:trans-m} and
   \eqref{eq:mnbd}, we find for $n > n(b)$ and $\eta$ sufficiently
   small that
   \begin{displaymath}
     M(sk_n; \bar \mu)
    \leq [ 1 + c_0 \tau_\eta^{-1} + q_n s (1+ c_b
    \delta_p(\epsilon_n))] k_n
    \leq [ 1 + q_n s + r_n]k_n,
   \end{displaymath}
where, since $s \leq (1-q')^{-1}$,
\begin{displaymath}
   r_n = c_0 \tau_\eta^{-1} + c \delta_p(\epsilon_n) \qquad \qquad
   c = c(b,q').
\end{displaymath}
We have
\begin{displaymath}
   M(sk_n; \bar \mu)/(sk_n)
   \leq 1 - q_n - d_n + q_n + r_n = 1 - d_n +r_n.
\end{displaymath}
Since $\delta_p(\epsilon_n) = o(\tau_\eta^{-1})$ (from the assumptions
on $\eta_n$), clearly $r_n - d_n < 0$ for $n > n(b)$;
for such $n$, $M(sk_n; \bar \mu) < s k_n$, and so
$k(\bar \mu) < s k_n$, as required.
\end{proof}

We draw a consequence for later use. Define
\begin{equation} \label{def:kappa}
   \kappa_n =
   \begin{cases}
     [\alpha_n + (1-q_n)^{-1}] k_n  & \text{for } \ell_0[\eta_n] \\
     [\alpha_n + (1-q_n - d_n)^{-1}] k_n  & \text{for } m_p[\eta_n].
   \end{cases}
\end{equation}
Recall now the notational assumption {\bf (A)}. Clearly, for large $n$,
\begin{equation}
   \label{eq:kappaconsequ}
   \kappa_n \sim (1-q_n)^{-1} k_n, \qquad \text{and} \qquad
   \kappa_n \leq  k_n/ q^{\prime\prime}.
\end{equation}
 From the remark after \eqref{eq:ktilde} (case
$\Theta_n = \ell_0[\eta_n]$), and from Corollary (\ref{cor:kmubound}) (case
$\Theta_n = m_p[\eta_n]$)
\begin{displaymath}
   \sup_{\mu \in \Theta_n} k_+(\mu) \leq \kappa_n.
\end{displaymath}

\section{Large Deviation bounds for $[\hat{k}_G, \hat{k}_F]$}
\label{sec:ldproof}
\setcounter{equation}{0}

We now develop exponential bounds on the FDR
interval $[ \hat k_G, \hat k_F]$ that lead to a proof of Proposition
\ref{prop:sandwichbd}. `Switching' inequalities allow the 
boundary-crossing definitions of $\hat k_G, \hat k_F$ to be 
expressed in terms
of sums of independent Bernoulli variables for which large deviation
inequalities in a `small numbers' regime can be applied.


\subsection{Switching Inequalities}
We will write $Y_l$ for the absolute ordered values $|y|_{(l)}$.  Let
$1 \leq \hat k_G \leq \hat k_F \leq n$ be respectively the smallest and largest
local minima of $k \rightarrow S_k = \sum_{l=k+1}^n Y_l^2 +
\sum_{l=1}^k t_l^2$ for $0 \leq k \leq n.$
The possibility of ties in the sequence $\{ S_k \}$ complicates the
exact description of local minima. Since ties occur with probability
zero, we will for convenience ignore this possibility in the arguments
to follow, lazily omitting explicit mention of ``with probability one''.

Define the exceedance numbers
\[
   N(t_k)  = \# \{ i: |y_i| > t_k \}, \qquad  N_+(t_k)  = \# \{ i: 
|y_i| \geq t_k \}.
\]
Clearly $N(t_k)$ and $N_+(t_k)$ have the same distributions. We now
have
\begin{align}
   \hat k_F & = \max \{ l: Y_l > t_l \}
              = \max \{ l: N(t_l) \geq l \}    \label{eq:kFdef} \\
   \hat k_G + 1 & = \min \{ l: Y_l < t_l \}
              = \min \{ l: N_+(t_l) < l \}   \label{eq:kGdef}.
\end{align}
[We set $\hat k_F =0$ or $\hat k_G = n$ if no such indices $l$ exist.]
To verify the left hand inequalities, note that
$S_k - S_{k-1} = t_k^2 - Y_k^2$, so that
\begin{displaymath}
   S_k \geq S_{k-1} \qquad \Leftrightarrow \qquad Y_k \leq t_k.
\end{displaymath}
The largest local minimum of $S_k$ occurs at $k = \hat k_F$ exactly
when $S_k < S_{k-1}$ but $S_l \geq S_{l-1}$ for all $l > k.$ In other
words $Y_k > t_k$ but $Y_l \leq t_l$ for all larger
$l$, which is precisely (\ref{eq:kFdef}).
Similarly, the smallest local minimum of $S_k$ occurs at $k = \hat
k_G$ exactly when $S_{k+1} > S_k$ but $S_l \leq S_{l-1}$ for all $l
\leq k$, and this leads immediately to (\ref{eq:kGdef}).

For the right-hand inequalities, we simply note that
\[
   N(t_k) \geq k \quad \mbox{iff} \quad Y_k > t_k, \qquad \mbox{and} \qquad
   N_+(t_k) < k \quad \mbox{iff} \quad Y_k < t_k.
\]


\subsection{Exponential Bounds}
First, recall Bennett's exponential inequality
(e.g. \citet[p. 192]{poll84}) in the form which states that for
independent, zero mean random variables $X_1, \dots, X_n$ with $|X_i|
\leq K$ and $V = \sum \mbox{Var} (X_i)$,
\begin{displaymath}
   P \{ X_1 + \ldots + X_n \geq \eta \} \leq \exp \Bigl\{ - \frac{\eta^2}{2 V}
   B(\frac{K \eta}{V}) \Bigr\},
\end{displaymath}
where $B(\lambda) = (2/\lambda^2) [(1+\lambda) \log(1 +\lambda) -
\lambda]$ for $\lambda > 0$ is decreasing in $\lambda$.
We adapt this for settings of
Poisson approximation.

\begin{lemma}   \label{lem:bennetcor}
   Suppose that $Y_l, l=1, \ldots , n$ are independent $0/1$ variables
   with $P(Y_l = 1) = p_l.$ Let $N = \sum_1^n Y_l$ and $M = EN =
   \sum_1^n p_l.$ Then
\begin{align}
   P \{ N \leq k \} & \leq \exp \{ - \tfrac{1}{4} M \, h(k/M) \}  \qquad
    \mbox{if} \ \ k < M,     \label{eq:kbound}  \\
   P \{ N \geq k \} & \leq \exp \{ - \tfrac{1}{4} M \, h(k/M) \}  \qquad
    \mbox{if} \ \ k > M,     \label{eq:kbound1}
\end{align}
where $h(x) = \min \{ |x-1|, |x-1|^2 \}$.
\end{lemma}
\begin{proof}
   In Bennett's inequality, if $M > k$, set $X_l = p_l - Y_l$ and $\eta =
M - k.$ If $M < k$, put $X_l = Y_l - p_l$ and
   $\eta = k - M.$  In each
   case $K=1$, and $V = Var N = \sum p_l(1-p_l)$, and we write $E =
   (\eta^2/V) B(\eta/V)$ for the term in the exponent.

Suppose first $\eta/V \leq 1$. Since $B$ is decreasing, $B(\eta/V)
\geq B(1)$, and since $V \leq M$, we have $\eta/V \geq \eta/M.$ Thus
$E \geq B(1) \eta^2 /M.$

Now if $\eta/V \geq 1$, we note from some calculus that the function
$C(\lambda) = \lambda B(\lambda)$ is increasing for $\lambda \geq 1.$
Hence $(\eta/V) B(\eta/V) \geq B(1) \geq B(1) \min( 1, \eta/M)$.
Thus $E \geq B(1) \min( \eta^2/M,\eta).$

Now $B(1) = 2[ 2 \log 2 - 1] \doteq .77 > 1/2.$
So, in either case $E \geq \hf \min( \eta^2/M,\eta) = (M/2) h( \eta
/M)$ and $ \eta / M  = | (k/M) - 1|$.
\end{proof}


\subsection{Bounds on $k/M_k$}
The Lipschitz properties of $k \rightarrow \dot M(k;\mu)$ established
in Section \ref{sec:mean-detect-funct} are now applied to derive
bounds on the ratios $k/M_k$ appearing in the exponential bounds
 (\ref{eq:kbound})-(\ref{eq:kbound1}).
In the following, we use $b$ to denote the vector of
constants $b = (b_1, ... , b_4,q')$, and the phrase $n > n_0(b)$ to 
indicate that a statement
holds for $n$ sufficiently large, depending on the constants $b$.

\begin{proposition}
   \label{prop:lipbd}
Assume hypotheses (Q), (H) and (A).
If  $\alpha_n k_n \leq k_1,$
then for $n > n(b)$, uniformly in $\mu \in \ell_0[\eta_n]$ and $m_p[\eta_n]$,
   \begin{equation}
     \label{eq:qpbound}
M(k_1 + \alpha_n k_n;\mu) - M(k_1; \mu) \leq  q' \alpha_n k_n.
   \end{equation}
(a) If $k(\mu) \leq k_1 \leq (1-q')^{-1} k_n $ and $k_2 = k_1 + \alpha_n k_n$,
then
\begin{equation}
   \label{eq:k2ratio}
   \frac{k_2}{M_{k_2}} - 1 \geq (1- q')^3 \alpha_n =:   \alpha_n.
\end{equation}
(b) There is $\zeta > 0$ so that,
if $2 \alpha_n k_n \leq k(\mu) \leq 
(1 - q')^{-1} k_n$ and
$k_1 = k(\mu) - \alpha_n k_n$, then
\begin{equation}
   \label{eq:k1ratio}
   1 - \frac{k_1}{M_{k_1}} \geq (1-q')^2 \alpha_n \geq \zeta \alpha_n.
\end{equation}
\end{proposition}

\begin{proof}
Formulas \eqref{eq:mdot}-\eqref{pkdd} show that $\nu \rightarrow \dot M(\nu;\mu)$ is
positive and decreasing and so the left side of \eqref{eq:qpbound} is
positive and bounded above by $\dot M(\alpha_n k_n; \mu).$
For $\mu \in \ell_0[\eta_n]$, \eqref{eq:l0mdotbd} shows that for $n >
n(b)$, $\dot M(\alpha_n k_n;\mu) \leq q'$ for all $\mu$ in
$\ell_0[\eta_n]$.
That the same bound holds uniformly over $m_p[\eta_n]$ also is a
consequence of \eqref{eq:lpderivbd} and \eqref{eq:b4bd}.

(a) To prove (\ref{eq:k2ratio}), note that the assumption $k(\mu) \leq
k_1$ entails $M_{k_1} \leq k_1$, so from (\ref{eq:qpbound}),
\begin{displaymath}
   M_{k_2} = M_{k_1} + M_{k_2} - M_{k_1} \leq k_1 + q' \alpha_n k_n.
\end{displaymath}
Since $k_1 \leq k_n/(1-q')$ and $q' < 1$, we have
\begin{displaymath}
   \frac{M_{k_2}}{k_2} \leq
   \frac{ k_1 + q' \alpha_n k_n}{ k_1 + \alpha_n k_n}
   \leq \frac{1 + q'(1-q') \alpha_n}{1 + (1-q') \alpha_n}.
\end{displaymath}
Thus, since $\alpha_n = O(1/\sqrt{ \log n})$, for $n > n(b)$, we find
\begin{displaymath}
     \frac{k_2}{M_{k_2}} - 1 \geq
     \frac{(1- q')^2 \alpha_n}{1 + q'(1-q') \alpha_n} \geq
     (1-q')^3 \alpha_n.
\end{displaymath}

(b) The assumption that $k(\mu) \geq 2 \alpha_n k_n$ yields $k_1 \geq
\alpha_n k_n$, and so Lipschitz bound (\ref{eq:qpbound}) implies
$M_{k(\mu)} - M_{k_1} \leq q' \alpha_n k_n$.
Hence, since $k(\mu) \leq k_n/(1-q')$,
\begin{displaymath}
   \frac{M_{k_1}}{k_1} \geq \frac{k(\mu) - q' \alpha_n k_n}{k(\mu) -
   \alpha_n k_n} \geq
   \frac{1 - q'(1-q')\alpha_n}{1 - (1-q')\alpha_n}
\end{displaymath}
which leads to (\ref{eq:k1ratio}) by simple rewriting.
\end{proof}



\subsection{Proof of Proposition \ref{prop:sandwichbd}}

   $1^\circ.$ Let $k_1 = k(\mu) \vee \alpha_n k_n$ and $k_2 = k_1 +
   \alpha_n k_n$. From \eqref{eq:kFdef},
\begin{equation}
     \label{eq:boole}
       \{ \hat k_F \geq k_2 \}
       \subset \bigcup_{l \geq k_2} \{ N(t_l) \geq l \}
\end{equation}
For $l > k_2 > k(\mu)$, we necessarily have $E_\mu N(t_l) = M(l;\mu) <
l$, and so from Lemma \ref{lem:bennetcor}
\begin{equation}
   \label{eq:lemma711}
     P_\mu \{ N(t_l) \geq l \} \leq \exp \{ - \tfrac{1}{4} M_l \, h(l/M_l) \},
   \qquad M_l = M(l;\mu).
\end{equation}
For $x \geq 1$, the function $h(x)$ is increasing, and for $l > k_2,$
$l \rightarrow l/M_l$ is increasing and so $h(l/M_l) \geq
h(k_2/M_{k_2})$. Now $k_1$ and $k_2$ satisfy the assumptions of
Proposition \ref{prop:lipbd}(a) and so from \eqref{eq:k2ratio},
$h(k_2/M_{k_2}) \geq \zeta^2 \alpha_n^2$.
Since $l \rightarrow M_l$ is increasing, we have from Proposition
\ref{prop:exceed-bound} that
\begin{equation}
   \label{eq:mlbd}
   M_l \geq M(1;\mu) \geq c (\log n)^{\gamma - 3/2}.
\end{equation}
Combining \eqref{eq:boole}, \eqref{eq:lemma711} and \eqref{eq:mlbd},
we find
\begin{align}
P_\mu \{ \hat k_F > k_2 \}
& \leq \sum_{l>k_2} \exp \{ - \tfrac{1}{4} M_1 \zeta^2 \alpha_n^2
\} \notag \\
& \leq n \exp \{ - c \alpha_n^2 ( \log n)^{\gamma - 3/2} \} \notag \\
& \leq n \exp \{ - c' (\log n)^{\gamma - 5/2} \},  \label{eq:tailbound}
\end{align}
for $c'$ depending on $q'$ and $b_4$.

$2^\circ$. \
Now assume that $k(\mu) \geq 2 \alpha_n k_n$; we establish a high
probability lower bound for $\hat k_G$.
Let $k_1 = k(\mu) - \alpha_n k_n$; from \eqref{eq:kGdef}
\begin{displaymath}
   \{ \hat k_G < k_1 \} = \{ \hat k_G + 1 \leq  k_1 \}
   \subset \bigcup_{l \leq k_1} \{ N_+(t_l) < l \}.
\end{displaymath}
For $l \leq k_1 < k(\mu)$, we necessarily have $M_l > l$ and so
\begin{displaymath}
   P \{ N_+(t_l) < l \} = P \{ N(t_l) < l \}
   \leq \exp \{ - \tfrac{1}{4} M_l h(l/M_l) \}.
\end{displaymath}
Since $l \rightarrow l/M_l \leq 1$ is increasing, and since $k_1$ and
$k(\mu)$ satisfy the assumptions of Proposition \ref{prop:lipbd}(b),
we obtain from \eqref{eq:k1ratio} that
\begin{displaymath}
   h(l/M_l) \geq \Bigl( 1 - \frac{k_1}{M_{k_1}} \Bigr)^2 \geq \zeta^2
   \alpha_n^2.
\end{displaymath}
In addition $l \rightarrow M_l$ is increasing, and so $M_l \geq M_1.$
Since $k(\mu) \geq 2 \alpha_n k_n \geq \alpha_n k_n$, we have from
Proposition \ref{prop:exceed-bound} that \eqref{eq:mlbd} holds here
also.
Hence
\begin{displaymath}
   P_\mu \{ \hat k_G < k_1 \}
   \leq k_1 \exp \{ - \tfrac{1}{4} M_1 \zeta^2 \alpha_n^2 \}
   \leq n   \exp \{ - c' ( \log n)^{\gamma - 5/2} \},
\end{displaymath}
in the same way as for \eqref{eq:tailbound}.

\section{Lemmas on Thresholding} \label{sec:threshlems}
\setcounter{equation}{0}

This section collects some preparatory results on hard (and in some
cases soft) thresholding with both fixed and data dependent
thresholds. These are useful for the analysis and comparison of the
various FDR and penalized rules, and are perhaps of some independent
utility. 

\subsection{Fixed Thresholds}
\label{sec:fixthresh}

First, an elementary decomposition of the $\ell_r$ risk
of hard thresholding.

\begin{lemma} \label{lem:threshdecomp}
Suppose that $x \sim N(\mu,1)$ and that $\eta_H(x,t) = xI\{ |x| \geq
t \}.$ Then
\begin{align}
\rho_H(t,\mu) = E|\eta_H(x,t) -\mu|^r & =
         \int_{-t}^t |\mu|^r \phi(x - \mu) dx + \int_{|x|>t} |x-\mu|^r
         \phi(x-\mu) dx \\
         & = D(\mu,t) + E(\mu,t),
         \label{eq:dedecomp}
\end{align}
where
\begin{align}
         D(\mu,t) & = |\mu|^r[\Phi(t-\mu) - \Phi(-t-\mu) ] \\
         E(\mu,t) & = |t-\mu|^{r-1} \phi(t - \mu) +
                       |t+\mu|^{r-1} \phi(t+\mu) +
                       \epsilon(t-\mu) + \epsilon(t+\mu) \label{eq:emut}\\
         |\epsilon(v)| & = |r-1| \int_v^\infty z^{r-2} \phi(z) dz
        \leq   v^{r-3} \phi(v), \qquad v > 0, 0 < r \leq 2.
         \label{eq:deexpr}
\end{align}
We note that for  $0\leq r \leq 2$
\begin{equation}
   \label{eq:hardriskat0}
   \rho_H(t,0) = 2 \int_t^\infty z^r \phi(z) dz
               = 2 t^r \tilde \Phi(t) (1 + \theta t^{-2}),
               \qquad 0 \leq \theta \leq r,
\end{equation}
and that $E(\mu,t)$ is
(i) globally bounded: $0 \leq E(\mu,t) \leq c_r = \int |z|^r \phi(z)
dz,$ \\
(ii) increasing in $\mu$, at least for $0 \leq \mu \leq t - \sqrt 2$,
\\
(iii) satisfies $E(1,t) \leq c_0 t^r \tilde \Phi(t-1)$ for $t > 1.$

A consequence of \eqref{eq:hardriskat0} is that for some $|\theta_2|
\leq 1$,
\begin{equation}
   \label{eq:penasy}
   \sum_{l=1}^k t_l^k = k t_k^r (1 + \theta t_k^{-2}) + \theta_2 t_1^r
       = k t_k^r (1 + o(1)),
\end{equation}
so long as $k \rightarrow \infty$ and $k/n \rightarrow 0$.
\end{lemma}

\begin{proof}
The risk decomposition is immediate. For \eqref{eq:hardriskat0},
by partial integration
   \begin{displaymath}
     \int_t^\infty z^r \phi(z) dz
      = t^r \tilde \Phi(t) + r \int_t^\infty z^{r-1} \tilde \Phi(z) dz,
   \end{displaymath}
from which the lower bound is clear. For the upper bound, use
\eqref{eq:mills} and $r \leq 2$ to find that the second integral is
bounded by
\begin{displaymath}
    \int_t^\infty z^{r-2} \phi(z) dz
    \leq t^{r-2} \tilde \Phi(t).
\end{displaymath}

Property (i) of $E$ is evident, while for $0 \leq \mu \leq t$, one
finds that
\begin{displaymath}
   \partial E / \partial \mu =
     (t - \mu)^r \phi(t-\mu) - (t + \mu)^r \phi(t + \mu),
\end{displaymath}
and the latter difference is positive if $(t-\mu)^2 > 2$.
For property (iii) we note that
$  E(1,t) \leq 2 \int_{t-1}^\infty z^r \phi(z) dz$
and then appeal to \eqref{eq:hardriskat0}.

Turning to \eqref{eq:penasy}, by integral approximation, we have for
$0 \leq \theta_2 \leq 1$,
\begin{displaymath}
   \sum_1^k t_l^r =
    \sum_1^k z^r( ql/2n) = (2n/q) \int_{t_k}^{t_1} v^r \phi(v) dv +
    \theta_2 t_1^r.
\end{displaymath}
 From \eqref{eq:hardriskat0}, the right hand integral term equals $k
t_k^r - t_1^r + \theta k t_k^{r-2}$ from which the result follows.
\end{proof}


The next lemma, on covariance between the data and hard thresholding,
is simpler in the $\ell_2$ case. The $\ell_r$ analogs are postponed to
Section \ref{subsec.lr-error-term}.

\begin{lemma}
\label{lem:cprop} Let $x \sim N(\mu,1)$.
$\xi(t,\mu) =  E_\mu (x - \mu)[ \eta_H(x,t) - \mu]$ has the properties
\begin{align}
           (i)\qquad  & \xi(t,\mu) = t [ \phi(t-\mu) + \phi(t+\mu)] +
\tilde{\Phi}(t-\mu) +
           \Phi(-t-\mu).  \label{cformula} \\
           (ii) \qquad & \xi(t,\mu) \leq
           \begin{cases}
             2 & \text{for } |\mu| \leq t - \sqrt{2 \log t}, \\
             t+1 & \text{for all } \mu.
           \end{cases}
\label{cbound} \\
           (iii) \qquad  & \mu \rightarrow \xi(t,\mu) \text{ is symmetric
            about } 0, \text{increasing for }  0 \leq \mu \leq t,
           \\
            & \text{and convex for } 0 \leq \mu \leq t/3, \text{ if } t \geq 3.
              \label{cmonotone} \\
           (iv) \qquad & \sup_{|\mu| \leq t/3 } | \xi_{\mu \mu}(t,\mu) |
              \leq c_0.
\label{eq:2ndderiv}
\end{align}
\end{lemma}

\begin{proof}
           Formula \eqref{cformula} follows by direct evaluation, and
inspection shows
           that $\xi(t,\mu)$ is symmetric about $\mu =0.$
The global bound in \eqref{cbound} follows since $2 \phi(0) \leq 1$.
If $0 \leq \mu \leq t - \sqrt{ 2 \log t}$, then
\begin{displaymath}
   \xi(t,\mu) \leq 2 t \phi(t-\mu) + 1
            \leq 2 t \phi( \sqrt{2 \log t}) + 1
            = 2 \phi(0) + 1 \leq 2.
\end{displaymath}
           For the monotonicity, writing $\xi_\mu$ for $\partial
           \xi/ \partial \mu$,
           \[
           \xi_\mu(t,\mu) = \phi(t-\mu) [t(t-\mu)+1] - 
\phi(t+\mu)[t(t+\mu) + 1].
           \]
           Since $\phi(t-\mu)/\phi(t+\mu) = e^{2t \mu} \geq 1 + 2 t
\mu$, it follows for $0 \leq t \leq \mu$ that
           \[
           \xi_\mu(t,\mu) / \phi(t+ \mu) \geq (1 + 2 t \mu)[t(t-\mu)+1]
- t(t+\mu) - 1
           = 2 t^2 \mu(t - \mu) \geq 0.
           \]
Finally, differentiating $\xi_\mu$ again with respect to $\mu$ yields
\begin{displaymath}
   \xi_{\mu \mu}(t,\mu) =  g(\mu;t) - g(-\mu;t),
\end{displaymath}
where $g(\mu;t) = [ t(t-\mu)^2 - \mu] \phi(t-\mu)$ and
\begin{displaymath}
    g_{\mu}(\mu; t) =
      [ (t-\mu) ( (t-\mu)^2 t - 2t -\mu ) - 1] \phi(t-\mu) \geq 0
\end{displaymath}
if $t>3$ and $|\mu| \leq t/3$. For such $\mu$, $\xi_{\mu \mu}(t,\mu)
\geq 0$, which establishes the convexity claimed in \eqref{cmonotone}.
Finally, for $0 \leq \mu \leq t/3$,
\begin{displaymath}
   \xi_{\mu \mu}(t,\mu) \leq g(t/3;t) \leq t^3 \phi(2t/3) \leq c_0.
\qedhere
\end{displaymath}
\end{proof}

\subsection{Data-dependent thresholds}
\label{sec:datdepthresh}

\begin{lemma}
           \label{lem:threshbd}
           Let $x = \mu + z \sim N(\mu,1)$ and $\eta(x,\hat{t} \, )$
denote soft or hard
           thresholding at $\hat{t}.$ For $r > 0$,
           \begin{equation}
                   | \eta(x,\hat{t} \,) - \mu|^r \leq 2^{(r-1)_+} ( \,
|z|^r + |\hat{t}|^r ).
                   \label{threshbd}
           \end{equation}
\end{lemma}
\begin{proof}
           Check cases and use $|a + b|^r \leq 2^{(r-1)_+} ( |a|^r + |b|^r )$.
\end{proof}

\begin{lemma}
           \label{lem:randthresh}
           Suppose that $y \sim N_n(\mu,I)$ and that $\hat{\mu}(y)$ corresponds
           to soft or hard thresholding at random level $\hat{t}:$
           $\hat{\mu}(y)_i = \eta (y_i,\hat{t} \, ).$ Suppose that $\hat{t}
           \leq t$ almost surely  on the event $A$ (with $t \geq [E
           |z|^{2r}]^{1/2r}$).
           Then for $r > 0$,
           \begin{equation}
                   E_{\mu} [ \| \hat{\mu} - \mu \|^r, A] \leq 2^{r \vee
                   1/2} t^{r}
                   n P_{\mu}(A)^{1/2}.
                   \label{randthresh}
           \end{equation}
\end{lemma}
\noindent \textit{Remark:} The notation $E[X,A]$ denotes $E X I_A$
where $I_A$ is the indicator function of the event $A$.
\begin{proof}
           Rewrite the left side and use Cauchy-Schwartz:
           \begin{displaymath}
                   \sum_{i=1}^n E [ \, | \eta(y_i,\hat{t} \, ) - \mu_i
|^r,A ] \leq
                   P(A)^{1/2} \sum_1^n \{ E [ | \eta_i - \mu_i |^{2r}, A
]\}^{1/2}.
           \end{displaymath}
           Now \eqref{threshbd} and the bound on $\hat{t}$ imply
           \begin{displaymath}
                   E[ \, | \eta_i - \mu_i |^{2r},A] \leq 2^{(2r-1)_+}  [
E |z_i|^{2r} +
                   \hat{t}^{2r}] \leq 2^{2r \vee 1} t^{2r}.
                   \qedhere
           \end{displaymath}
\end{proof}

\medskip
Continuing with $y \sim N_n(\mu,I)$,
the next lemma matches the $\ell_r$ risks of two hard threshold estimators
$\hat \mu(y)_i = \eta_H (y_i; \hat t)$ and $\hat{\mu}'$ with data-dependent
thresholds $\hat t$ and $\hat{t}'$ if those thresholds are close.  Assume also
that there is a non-random bound $t$ such that $\hat t, \hat{t}'
\leq t$ with probability one. Then
\begin{displaymath}
     | \eta_H(y_i,\hat t) - \eta_H(y_i, \hat t') |\
     \leq
     \begin{cases}
       t & \text{if} \ |y_i| \ \ \text{lies between} \ \hat t, \hat t', \\
       0 & \text{otherwise}.
     \end{cases}
\end{displaymath}
Let $N' = \# \{ i: |y_i| \in [ \hat t , \hat t' ] \}$-- clearly
$  \| \hat \mu - \hat \mu' \|_r^r \leq t^r N'$.
In various cases, $N'$ can be bounded
on a high probability event, yielding

\begin{lemma}
     \label{lem:riskdel}
     Let $\beta_n$ be a specified sequence, and
     with the previous definitions, set $B_n = \{ N' \leq \beta_n
     \}$. For $0 < r \leq 2$,
     \begin{equation}
       \label{eq:riskdel}
       \begin{split}
       |\rho(\hat \mu, \mu) - \rho(\hat \mu', \mu)| \leq
          2 \beta_n t^r & +
          r I \{ r > 1 \} \rho (\hat \mu', \mu)^{1 - 1/r} \beta_n^{1/r} t \\
          & + 8 t^r n P_\mu (B_n^c)^{1/2}.
       \end{split}
     \end{equation}
\end{lemma}

\begin{proof}
     To develop an $\ell_r$ analog of \eqref{eq:lossident}, we note a
simple bound valid for all $a, z \in \mathbb{R}$:
\begin{equation}
           \Bigl\lvert \ |a+z|^r - |a|^r \Bigr\rvert \leq
                   \begin{cases}
                   |z|^r & 0 < r \leq 1 \\
                   r( \, |a| + |z| )^{r-1} |z| & 1 < r.
           \end{cases}
           \label{eq:auxilbd}
\end{equation}
[For $r > 1$, use derivative bounds for $y \rightarrow |y|^r$].
We consider here only $1 < r \leq 2$: the case $r \leq 1$ is similar
and easier. Thus, setting $\epsilon = \hat{\mu} - \hat{\mu}'$,
$\Delta = \hat \mu - \mu$ and similarly for $\Delta'$
\begin{displaymath}
           | E \{ \sum_i |\Delta_{i}|^r - |\Delta_{i}'|^r, B_n \}  |
           \leq r E \{ \sum_i | \Delta_{i}'|^{r-1} |\epsilon_i| +
           |\epsilon_i|^r, B_n \}.
\end{displaymath}
Using H\"older's inequality and defining $\varepsilon_n = E \{ \|
\hat{\mu} - \hat{\mu}' \|_r^r, B_n \},$ we obtain
\begin{equation}
           \bigl\lvert  E \{ \| \Delta \|_r^r - \| \Delta' \|_r^r , B_n \}
           \bigr\rvert \leq  r \rho (\hat{\mu}',\mu)^{(r-1)/r}
           \varepsilon_n^{1/r} + r \varepsilon_n.
           \label{eq:bdiff}
\end{equation}
 From the definition of event $B_n$ and the remarks preceding the
lemma, $\varepsilon_n \leq \beta_n t^r$.  On $B_n^c$, apply Lemma
\ref{lem:randthresh} to obtain \eqref{eq:riskdel}.
\end{proof}


\section{Upper Bound Result: $\ell_2$ error}
\label{sec:proof}
\setcounter{equation}{0}

We now turn to the upper bound,
Theorem \ref{th:mainres}.
We begin with the simplest case: squared-error loss.
Only the outline of the argument is presented in this
section, with details provided in the next section.
The extensions to $\ell_r$ error
measures, of considerable importance to the conclusions of the
paper, are not straightforward. The proofs
are postponed until Section \ref{sec.ellr}.

The approach taken in this section was sketched in the
introduction, see \eqref{eq:penalized} and \eqref{eq:keybd}.
We define certain
empirical and theoretical complexity functions -- the empirical
complexity being minimized by $\hat{\mu}_2$.  A basic inequality
bounds the theoretical complexity of $\hat{\mu}_2$ by the minimal
theoretical complexity plus
an ``error term'' of covariance type.  When
maximized over a sparse parameter set $\Theta_n$, the minimum theoretical
complexity has the same leading asymptotics as
the minimax estimation risk for $\Theta_n$.
To complete the proof, the error term is bounded.
This analysis is sketched in Section \ref{ssec:error-term}; the
detailed proof relies on an average case and large deviations analysis
of the penalized FDR index $\hat{k}_2$.
The upshot is that these terms are negligible if $q \leq 1/2$, but
add substantially to the maximum risk when $q  > 1/2$; this
was 
foreshadowed Proposition \ref{prop:riskmualpha}
and its discussion.
Finally, a risk comparison is
used to extend the conclusion from the penalized estimate
$\hat{\mu}_2$ to the original FDR estimate $\hat{\mu}_F$.

\subsection{Empirical and Theoretical Complexities}
\label{ssec:complexities}

In Section \ref{subsec.variational}, we have defined
$\hat{\mu}_2$ as the minimizer of
the \textit{empirical} complexity $K(\tilde{\mu},y) = \| y -
\tilde{\mu} \|^2 + Pen(\tilde{\mu})$ (note that now we set $\sigma^2
=1).$
Substituting $y = \mu + z$ into $K(\hat{\mu}_2,y)$ yields the
decomposition
           \begin{equation}
                   K(\hat{\mu}_2,y) = K(\hat{\mu}_2,\mu) + 2 \langle z, \mu -
                           \hat{\mu}_2 \rangle + \| z \|^2.
                   \label{eq:kdecomp}
           \end{equation}
Now let $\mu_0 = \mu_0(\mu)$ denote the minimizer over $\tilde{\mu}$ of the
\textit{theoretical} complexity $K(\tilde{\mu},\mu)$ corresponding to the
unknown mean vector $\mu:$
\begin{equation}
           K(\mu_{0},\mu) = \inf_{\tilde{\mu}} \| \mu - \tilde{\mu} \|^2 +
          Pen(\tilde{\mu})
           \label{eq:theorcompl}
\end{equation}
There is a decomposition for $K(\mu_0,y)$ that is exactly analogous to
\eqref{eq:kdecomp}:
\begin{displaymath}
           K(\mu_0,y) = K(\mu_0,\mu) + 2 \langle z, \mu - \mu_0 \rangle
+ \| z \|^2.
\end{displaymath}
Since by definition of $\hat{\mu}_2$, $ K(\mu_0,y) \geq K(\hat{\mu}_2,y)$,
we obtain, after noting the cancellation of the quadratic error terms
and rearranging,
           \begin{equation}
                   K(\hat{\mu}_2,\mu) \leq K(\mu_0,\mu) + 2 \langle z,
\hat{\mu}_2 - \mu_0
                   \rangle.
                   \label{complexbound}
           \end{equation}
Thus the complexity of $\hat{\mu}_2$ is bounded by the minimum
theoretical complexity plus an error term.
Up to this point, the development is close to that of \citet{dojo94c},
as well as work of other authors (e.g.
\citet{geer90}).
Since
         \begin{equation}
                 K(\hat{\mu}_2,\mu)  =  \| \hat{\mu}_2 - \mu \|^2 +
                 Pen(\hat{\mu}_2),
         \label{eq:kmuhat}
\end{equation}
we obtain a bound for $\rho(\hat{\mu}_2,\mu) = E_{ \mu} \| \hat{\mu}_2 -
\mu \|^2$ by taking expectations in \eqref{complexbound}.  Since the
error term has zero mean, we may replace $\mu_0$ by $\mu$ and obtain
the basic bound
         \begin{equation}
                 \rho(\hat{\mu}_2,\mu) \leq K(\mu_0, \mu) +
                 2 E_{\mu} \langle z, \hat{\mu}_2 - \mu \rangle
                 - E_\mu Pen(\hat{\mu}_2).
                 \label{basicbound}
         \end{equation}
   
We view the right-hand side as containing
a `leading term' 
$K(\mu_0,\mu)$ -- the theoretical complexity --
and an `error term' 

\begin{equation} \label{sim-z}
Err_2(\mu,\hat{\mu}_2) \equiv 
2E_{\mu} \langle z, \hat{\mu}_2 - \mu \rangle
                 - E_\mu Pen(\hat{\mu}_2).
\end{equation}

We now 
claim that the maximum theoretical complexity
over sparsity classes $\Theta_n$ is asymptotic
to the minimax risk.

\begin{proposition} \label{prop:theocomplex}
Assume (Q), (H).
          \begin{equation}
                   \sup_{\mu \in \Theta_n} K(\mu_0(\mu),\mu) \sim
                   R_n(\Theta_n), \qquad \qquad n
                           \rightarrow \infty.
                   \label{eq:maxtheor}
           \end{equation}
\end{proposition}

The same minimax risk bounds the error 
term:
\begin{proposition} \label{prop:errmax}
Assume (Q), (H). With 
$u_{2p} = 1$ for $\ell_0$ and strong $\ell_p$,
and $u_{2p} = 1 - p/2$ 
for weak $\ell_p$.
          \begin{equation}
                   \sup_{\mu \in \Theta_n} Err_2(\mu,\hat{\mu}_2)   = 
     \left\{ \begin{array}{ll}
                  R_n(\Theta_n) \cdot 
u_{2p} \cdot  \frac{(2q -1)}{1-q}  & q > 1/2, \\
                   o(R_n(\Theta_n)) & q \leq 1/2 . 
 
\end{array} \right.                  \label{eq:maxerror}
           \end{equation}
\end{proposition}

Together, these propositions give the
Upper Bound 
result in the squared-error case. 

\subsection{Maximal Theoretical Complexity}
\label{ssec:max-theoretical-complexity}

We prove Proposition \ref{prop:theocomplex},
beginning with the
nearly-black classes $\Theta_n = \ell_0[\eta_n]$.

As in Section \ref{subsec.variational},
decompose the optimization problem \eqref{eq:theorcompl} defining
$K(\mu_0,\mu)$ over the number of non-zero components in
$\tilde{\mu}.$ Assign these to the largest components of $\mu$: hence
\begin{equation}
           K(\mu_0,\mu) = \inf_{k} \ \sum_{l=k+1}^{n}
           \mu_{(l)}^2 + \sum_{l=1}^k t_l^2.
           \label{eq:kmu0mu}
\end{equation}
On $\ell_0[\eta_n],$ at most $k_n = [n \eta_n]$ components of $\mu$
can be non-zero.  Hence the infimum over $k$ may be restricted to $0
\leq k \leq k_n$.  This implies
\begin{equation}
               \sup_{\mu \in \ell_0[\eta_n]} K(\mu_0,\mu) = \sum_1^{k_n} t_l^2.
              \label{eq:equality}
          \end{equation}
Indeed, choosing $k = k_n$ in \eqref{eq:kmu0mu} shows the left side
to be smaller than the right side in \eqref{eq:equality},
while equality occurs for any $\mu$ with non-zero entries
$\mu_1 = \ldots = \mu_{k_n}>t_1.$
Noting 
       \begin{equation}
                   \sum_1^k t_l^2 \sim k t_k^2, \qquad \qquad
                  t_k^2 \sim 2 \log ( \tfrac{n}{k} \tfrac{2}{q} )
                  \qquad \qquad \text{if}  \ \ k = o(n).
                   \label{slowlyvar}
        \end{equation}
(cf. \eqref{eq:penasy} and Lemma \ref{lem:quantile}),
along with $\eta_n = O(n^{-\delta})$, we get
\begin{displaymath}
           \sum_1^{k_n} t_l^2 \sim k_n t_{k_n}^2 \sim n \eta_n \cdot 2 \log
           \eta_n^{-1} \sim R_n (\ell_0[\eta_n])
\end{displaymath}
which establishes \eqref{eq:maxtheor} in the $\ell_0$ case.

\medskip

\textit{Remark.\ } Using a fixed penalty $Pen_{fix}(\mu) = t^2 \| \mu \|_0$
in the above argument would yield $\sup K = k_n t^2 \approx n \eta_n
t^{2},$ but the $t^{2}$ term is unable to adapt to varying
signal complexity.

\medskip

\small
\textbf{Weak $\ell_p$.} \
The maximum of \eqref{eq:kmu0mu} over $\mu \in m_p[\eta_n]$ occurs at the
extremal vector $\bar{\mu}_l = C_{n} l^{-1/p}$, where $C_{n} = n^{1/p}
\eta_n.$
Define $\underline{k}_n$ to be the solution of $C_{n}^2
\underline{k}_n ^{-2/p} =
t_{\underline{k}_n }^2.$
Using \eqref{slowlyvar}, we obtain
\begin{align}
           \sup_{\mu \in m_p[\eta_n]} K(\mu_0,\mu)
           & = \inf_{k} C_{n}^2 \sum_{k+1}^{n}  l^{-2/p} + \sum_1^k
t_l^2  \notag \\
           & \sim C_{n}^2 \tau_p \underline{k}_n ^{1 - 2/p} +
\underline{k}_n  t_{\underline{k}_n }^2
          \qquad \qquad \tau_p = \tfrac{p}{2-p}   \notag \\
           & =     (1 + \tau_p ) \underline{k}_n  t_{\underline{k}_n
}^2. \label{supstep}
\end{align}
Thus $\underline{k}_n $ is
the optimal number of non-zero components and may be rewritten as
           \begin{equation}
                   \underline{k}_n  = C_{n}^p t_{\underline{k}_n }^{-p}
= n \eta_n^p t_{\underline{k}_n }^{-p}.
                   \label{eq:kstareqn}
           \end{equation}
A little analysis using \eqref{slowlyvar} and the equation for
$\underline{k}_n $ shows that $t_{\underline{k}_n }^{2} \sim 2 \log
\eta_n^{-p}$.
[For this reason, we define $k_n = n \eta_n^p \tau_\eta^{-p}.$]
From \eqref{munasymp} we then conclude $t_{\underline{k}_n }^{2} \sim
\mu_n^2,$
which via \eqref{eq:lpmmxrisk} and \eqref{eq:rmpn} shows that the right side
of \eqref{supstep} is asymptotically equivalent to $R(m_p[\eta_n])$,
as claimed.

\medskip

\textit{Remark.}
The least favorable configuration for $\mu$ is thus given by $\mu_l =
\min(C_{n} l^{-1/p}, t_l) $ $= \min(\eta_n (l/n)^{-1/p}, t_l) $, which, after
noting \eqref{slowlyvar}, is
essentially identical with the least favorable distribution
\eqref{weak-least-fav}. In addition, the maximisation has exactly the
same structure as the approximate evaluation of the Bayes risk of
soft thresholding over this least favorable distribution; compare
\eqref{risk1} - \eqref{risk2}.
Replacing the slowly varying boundary $l \rightarrow t_l$ by $m_k =
t[k_n] + \alpha$ leads to the configurations \eqref{eq:mualpham}.

\small
\textbf{Strong $\ell_p$.} \
The maximal theoretical complexity
is the value of the optimization 
problem
\[
   (Q(n,\eta_n^p)) \qquad \max \sum \min( \mu_{(l)}^2, 
t_\ell^2 ) \mbox{  subject to } \sum_\ell \mu_{(\ell)}^p \leq n 
\eta_n^p .
\]
The change of variables  $x_\ell = \mu_{(\ell)}^p$ 
allows to write this as
\[
   \max \sum_\ell \min( x_\ell^{2/p} , 
t_\ell^2 ) \mbox{  subject to } \sum_\ell x_{\ell} \leq n \eta_n^p, 
\quad x_1 \geq x_2 \geq \dots  .
\]
Since $p < 2$, the objective 
function is strictly convex on $ \Pi_\ell [0,t_\ell^2]$,
and the constraint set 
is convex.  The maximum will be obtained
at an extreme point of the 
constraint set, i.e. roughly at a sequence vanishing
for $\ell > k$ (for some $k$) 
and equal to $t_\ell^{1/p}$ for $\ell \leq k$.
Let $\tilde{k}_n$ be the 
largest $k$ for which
\[  
     \sum_{\ell=1}^k t_\ell^p \leq n 
\eta_n^p .
\]
Then the maximal theoretical complexity obeys
\[
\sum_{\ell=1}^{\tilde{k}_n} t_\ell^2 \leq val(Q(n,\eta_n^p)) \leq 
\sum_{\ell=1}^{\tilde{k}_n+1} t_\ell^2 .
\]
Again using 
(\ref{slowlyvar}), we get
\[
   \sup \{  K(\mu_0,\mu) : \mu \in 
\ell_p[\eta_n] \} \sim \tilde{k}_n t[\tilde{k}_n]^2  \sim n \eta_n^p 
\tau_\eta^{2-p}.
\]
So (\ref{eq:maxtheor}) follows in the 
$\ell_p[\eta_n]$ case.
\normalsize

\subsection{The Error term}
\label{ssec:error-term}

We outline the proof of \eqref{eq:maxerror}.
Recall the definitions 
$k_\pm$ and $t_\pm$ from Section 5,
at equations 
(\ref{eq:kplusdef}),(\ref{eq:kminusdef}),(\ref{eq:tplusminusdef}),
and 
their use in Proposition 5,1.
We will rely on the fact that it is $\Theta_n$-likely under $P_\mu$
that $\hat{k}_2 \leq k_+(\mu)$ and hence that $\hat{t}_2 \geq t_{-}(\mu).$
First write
\begin{equation}
                \langle z, \hat{\mu}_2 - \mu \rangle = \sum_1^n z_i [
                 \eta_H(y_i,\hat{t}_2) - \mu_i].
                       \label{eq:errorterm}
\end{equation}
We exploit monotonicity of the error term for \textit{small}
components $\mu_i$. Indeed, if $|\mu_i| \leq t_-(\mu)  \leq \hat{t}_2$, then
(cf. Lemma \ref{lem:monotonicity})
         \begin{equation}
                 z_i [ \eta_H(y_i,\hat{t}_2) - \mu_i] \leq
                 z_i [ \eta_H(y_i, t_-(\mu) ) - \mu_i],
                 \label{monotone}
         \end{equation}
as may be seen by checking cases.
This permits us to replace $\hat{t}_2$ by the fixed threshold value
$t_-(\mu) $  for the vast majority of
components $\mu_i,$ with Proposition \ref{prop:sandwichbd}
providing assurance that $\hat{t}_2 \geq t_-(\mu) $ with high probability.
We recall the defiinition of the covariance
kernel $\xi$ in Lemma 
\ref{lem:cprop}.
For $z_1 \sim N(0,1)$ and scalar mean $\mu_1$,
\begin{displaymath}
         \xi(t,\mu_1) = E \  z_1 [ \eta_H( z_1 + \mu_1 , t) - \mu_1].
\end{displaymath}

The function $\xi(t,\mu_1)$ is the covariance between $y_1$ and
$\eta_{H}(y_1,t)$ when the data $y_1 \sim N(\mu_1,1)$.  Lemma
\ref{lem:cprop} shows that $\xi$ is even in $\mu_1$, has a minimum of $2 t
\phi(t)$ at $\mu_1 =0,$ rising to a maximum near $\mu_1 = t$ (though
always bounded by $t + 1)$, and dropping quickly to 1 for large $\mu_1.$
It turns out that, uniformly on nearly-black sequences $\mu \in
\Theta_n,$ the main contribution to the sum comes from components
$\mu_1$ near $0$.

Similarly, it is $\Theta_n$-likely that
\begin{displaymath}
     Pen(\hat{\mu}_2) = \sum_1^{\hat{k}_2} t_l^2 \geq k_-(\mu)
     t_+^2(\mu).
\end{displaymath}

We proceed heuristically here, leaving the (necessary!) careful
verification to Section \ref{sec:errtermproof}.
Interpreting "$\approx$" to mean ``up to terms whose positive part
is $o(R_n(\Theta))$'', we have, uniformly on $\Theta_n$,
\begin{align}
    Err_2(\mu,\hat{\mu}_2)
     &  \approx 2 \sum \xi(t_-(\mu),\mu_i) - k_-(\mu) t_+^2(\mu) 
\label{sim-a}\\
      & \approx 4 n t_-(\mu) \phi(t_-(\mu)) - k_-(\mu) t_+^2(\mu) 
\label{sim-b} \\
      & \approx [ 4n \tilde{\Phi}(t_-(\mu)) - k_-(\mu)] t_+^2(\mu) 
\label{sim-c}  \\
      & \approx (2 q_n - 1) k_-(\mu) t_+^2(\mu)  .  \label{sim-d}
\end{align}
At \eqref{sim-b} we have first used the fact that $\xi(t,0) = 2 t
\phi(t).$
Second, for the at most $k_n$ non-zero terms, we use the bound
$\xi(t,\mu) \leq t+ 1 \leq t_1 + 1$ \-- compare (\ref{cbound}) \-- and note
that their contribution is at most $O(k_n t_1) = o(R_n)$.
At \eqref{sim-c} we used $\phi(t) \sim t \tilde \Phi(t)$
as $t \rightarrow \infty$, and at \eqref{sim-d} the definitions of
$t_+(\mu)$, $k_-$ and the asymptotics of each.

Expression \eqref{sim-d} is negative if $q_n< 1/2$,
making $Err_2$ 
for our purposes negligible.
For $q_n \geq 1/2,$ since $k \rightarrow k t_k^2$ is increasing,
cf. \eqref{eq:nutnur}, we use the bound of \eqref{eq:ktilde},
namely $k(\mu) \leq \tilde k_n$ on $\ell_0[\eta_n]$, along with
(\ref{slowlyvar}) to conclude that  (\ref{sim-d}) is not larger than
\begin{equation}
\label{eq:ktk2}
     (2 q_n - 1) \tilde{k} t_{\tilde{k}}^2 \sim \frac{2 q_n - 1}{ 1 -
     q_n} R_n( \ell_0[\eta_n] ).
\end{equation}
This motivates (\ref{eq:maxerror}) in the $\ell_0$ case.

\small
\textbf{Weak $\ell_p$ }
The outline is much as above, although there is detailed technical
work since all means may be non-zero (subject
to the weak-$\ell_p$ sparsity constraint); for example in the
transition \eqref{sim-a} to \eqref{sim-b}. Again, after \eqref{sim-d}
we are led to maximize $k_-(\mu)$ over $m_p[\eta]$ and from
\eqref{eq:kmubar}, find $k_-(\mu) \leq k( \bar \mu) \leq k_n/(1 - q_n) =
\tilde k$, say.
Here $k_n = n \eta_n^p \tau_\eta^{-p}$ is the effective non-zero index
for weak $\ell_p$ defined after \eqref{eq:kstareqn}.

Consequently, since $t_{\tilde k} \sim \tau_\eta$, and using the
expressions \eqref{eq:lpmmxrisk} and \eqref{eq:rmpn} for minimax
risks, we obtain
\begin{align*}
   (2q_n-1 ) k_{-}(\mu) t_{+}^2(\mu)
     & \leq (2q_n-1 ) \tilde k t^2_{\tilde k} (1 + o(1)) \\
     & \sim \frac{2q_n-1}{1-q} k_n\tau_\eta^2  \sim 
\frac{2q_n-1}{1-q}  n \eta_n^p \tau_\eta^{2-p} \\
     & \sim  \frac{2q_n-1}{1-q}  R_n( \ell_p[\eta_n])
       \sim u_{2p} \cdot  \frac{2q_n-1}{1-q}  \cdot R( m_p[\eta_n]),
\end{align*}
with $u_{2p} = (1-p/2)$.

\indent
\textbf{Strong $\ell_p$ } The inclusion $\ell_p[\eta_n] \subset 
m_p[\eta_n]$
and the previous display give
\[
  (2q_n-1 ) k_-(\mu) 
t_{+}^2(\mu) \leq   \frac{2q_n-1}{1-q}  R_n( \ell_p[\eta_n])  (1 + 
o(1)).
\]
\normalsize

If the above arguments were complete  \-- rather than just sketches 
\--
we would now have the right side of
\eqref{eq:maxerror} in Proposition \ref{prop:errmax}.  Details to 
come in Section \ref{sec:errtermproof}.

\subsection{From penalized to original FDR}

We extend the adaptive minimaxity result for the penalized estimator
$\hat{\mu}_2$ which thresholds at $\hat{t}_2$ to any threshold
$\hat{t}$ in the range $[\hat{t}_F, \hat{t}_{G}]$ defined in
Section \ref{subsec.variational}. In
particular the adaptive minimaxity will apply to the original FDR
estimator $\hat{\mu}_F$.

First compare the squared error of a deviation
$\hat{\delta}_2 = \hat{\mu}_2 - \mu$ with that of $\hat{\delta} =
\hat{\mu} - \mu$:
\begin{equation}
                   \| \hat{\delta} \|_2^2 - \| \hat{\delta}_2 \|_2^2 =
                   \| \hat{\mu} - \hat{\mu}_2 \|_2^2 + 2 (\hat{\mu} -
\hat{\mu}_2)
                   \cdot (\hat{\mu}_2 - \mu).
           \label{eq:lossident}
\end{equation}
Now suppose $\hat{\mu}_D$ (``$D$'' for ``data dependent) has the 
form
(\ref{eq:muhatlim})
All such estimators differ from $\hat{\mu}_2$ at most in those
co-ordinates $y_l$ with $\hat{k}_G \leq l \leq \hat{k}_F$, and on
such co-ordinates the difference between the two estimates is at most
$\hat{t}_G \leq t_1 = z(q/2n).$ Hence
\begin{displaymath}
           \| \hat{\mu} - \hat{\mu}_2 \|_2^2 \leq t_1^2 ( \hat{k}_F -
           \hat{k}_G ).
\end{displaymath}
Proposition \ref{prop:sandwichbd} and \eqref{eq:kpminuskm}
show that FDR control, combined
with sparsity, forces the ``crossover interval'' $[\hat{k}_G,
\hat{k}_F]$ to be relatively small, having length bounded by $\gamma
\alpha_n k_n.$

On the event described in Proposition \ref{prop:sandwichbd}, we have
\begin{equation}
       \| \hat{\mu} - \hat{\mu}_2  \|_2^2 \leq \gamma \alpha_n k_{n} t_1^2 =
                   o( R_n (\Theta_n) ).
  \label{eq:mumuA}
\end{equation}
We summarize, with remaining details deferred to Section
\ref{sec:pen-to-orig-again}.

\begin{theorem} \label{th:pen-to-fdr}
    If $\hat{\mu}_D$ satisfies \eqref{eq:muhatlim}, then
for each $r \in (0,2]$
        \begin{displaymath}
                   \sup_{\mu \in \Theta_n} | \rho(\hat{\mu}_D,\mu) -
\rho(\hat{\mu}_r,\mu)| =
                   o( R_n (\Theta_n) ),
           \end{displaymath}
           so that asymptotic minimaxity of $\hat{\mu}_r$ implies the same
           property for any such $\hat{\mu}_D.$
\end{theorem}


\section{Error term: Quadratic Loss}
\label{sec:errtermproof}
\setcounter{equation}{0}

We now formalize the error term analysis of
Section \ref{ssec:error-term}, collecting and
applying the tools built up in earlier sections.

\begin{lemma}
   \label{lem:monotonicity}
   If $|\mu| \leq t^1 \leq t^2$, then
   \begin{displaymath}
     (x-\mu) [ \eta_H(x,t^2) - \mu] \leq (x-\mu) [ \eta_H(x,t^1) - \mu].
   \end{displaymath}
\end{lemma}
\begin{proof}
   The difference RHS - LHS equals
   \begin{displaymath}
     (x-\mu) [ \eta_H(x,t^1) -  \eta_H(x,t^2)] =
     (x-\mu) x I \{ t^1 \leq |x| \leq t^2 \} \geq 0,
   \end{displaymath}
since $\mbox{sgn} \, x = \mbox{sgn} (x-\mu)$ if $|x| \geq t^1 \geq |\mu|.$
\end{proof}

We proceed with the formal analysis of the error term \eqref{sim-z}.
Set
$$\hat{e}_i = e_i(\hat{t}_2) = 2 (y_i - \mu_i)
[\eta_H(y_i,\hat{t}_2) - \mu_i],$$
and
$$A_n = A_n(\mu) = \{ t_- \leq \hat{t}_2 \leq t_+  \}, \qquad
S_n(\mu) = \{ i: | \mu_i | \leq t_-  \}.
$$
We have
\begin{align}
         2 E \langle z, \hat{\mu}_2 - \mu \rangle   =  E \sum \hat{e}_i
         & = E[ \sum_{S_n(\mu)} \hat{e}_i, A_n ] + E [ \sum_{S_n^c(\mu)}
         \hat{e}_i, A_n ] + E [ \sum \hat{e}_i, A_n^c].
         \label{eq:err-decomp} \\
         & = D_{an} + T_{2n} + T_{3n},
\end{align}
where we use $D_{an}, D_{bn}$ etc to denote `dominant' terms, and
$T_{jn}$ to denote terms that will be shown to be negligible.

Let $e_i = e_i(t_-):$ the monotonicity of errors for small
components (shown in Lemma \ref{lem:monotonicity})
says that the first term on the
right side is bounded above by
\begin{displaymath}
         E[ \sum_{S_n(\mu)} e_i, A_n ]
          = E[ \sum_{S_n(\mu)} e_i ] - E[ \sum_{S_n(\mu)} e_i, A_n^c ]
          = D_{bn} + T_{4n}.
\end{displaymath}
Recalling the definition of $\xi(t,\mu)$ from Section \ref{lem:cprop},
we have $E e_i = 2 \xi(t_-, \mu_i)$ and
\begin{align*}
   D_{bn} & = 2 |S_n(\mu)| \xi(t_-,0) + 2 \sum_{S_n(\mu)} [ \xi(t_-,\mu_i) -
   \xi(t_-,0)] \\
          & \leq 2 n \xi(t_-,0) + T_{1n}(\mu),
\end{align*}
say. To summarize, we obtain the following decomposition for the error term
\eqref{sim-z}:
\begin{displaymath}
   Err_2(\mu,\hat{\mu}_2)
     \leq D_{cn}(\mu) + \sum_{j=1}^4 T_{jn}(\mu),
\end{displaymath}
where
\begin{displaymath}
   D_{cn}(\mu) = 2 n \xi(t_-,0) - E_\mu Pen(\hat \mu_2).
\end{displaymath}

Recall that $R_n(\Theta_n) \asymp k_n \tau_\eta^2$ for both
$\ell_0[\eta_n]$ and $m_p[\eta_n]$. In the following, we will show
negligibility of error terms by establishing that they are $O(k_n
\tau_\eta)$, (or, in one case, $o(k_n \tau_\eta^2)$) uniformly over
$\ell_0[\eta_n]$ or $m_p[\eta_n]$ respectively.

\bigskip
\textit{Dominant term.}
Using \eqref{eq:phiupper},
\begin{displaymath}
   \xi(t,0) = 2 t \phi(t) + 2 \tilde \Phi(t) \leq 2 t^2 \tilde \Phi(t) +
   6 \tilde \Phi(t).
\end{displaymath}
Since $2 \tilde \Phi(t_-) = q_n k_+ n^{-1}$, we obtain
\begin{equation}
   \label{eq:cat0bd}
   2n \xi(t_-,0) \leq 2 q_n k_+ t_-^2 + 6 q_n k_+
               \leq 2 q_n k_- t_-^2 + c k_n \tau_\eta,
\end{equation}
after observing that $k_+ - k_- \leq 3 \alpha_n k_n \leq c k_n
\tau_\eta^{-1}$ by (A), and that $t_- \leq c \tau_\eta$ from
\eqref{eq:t1taueta}.

For the penalty term in $D_{cn}$,
we note that on $A_n$, $\hat k_2 \geq k_-$, and
so $Pen(\hat \mu_2) = \sum_1^{\hat k_2} t_l^2 \geq k_- t_+^2 \geq k_-
t_-^2.$
On the other hand, since $A_n$ is $\Theta_n-$likely,
\begin{displaymath}
   E_\mu \bigl[ \sum_1^{\hat k_2} t_l^2 , A_n^c \bigr]
   \leq n t_1^2 P_\mu(A_n^c) \leq c k_n \tau_\eta.
\end{displaymath}
As a result
\begin{equation}
   \label{eq:penbd}
   E_\mu Pen(\hat \mu_2) \geq k_- t_-^2 + O(k_n \tau_\eta).
\end{equation}

\bigskip
Combining \eqref{eq:cat0bd} and \eqref{eq:penbd}, we obtain
\begin{displaymath}
   D_{cn}(\mu) \leq (2 q_n - 1) k_- t_-^2 + O(k_n \tau_\eta).
\end{displaymath}
If $q_n \leq 1/2$, then of course the leading term is non-positive,
while in the case $1/2 \leq q_n < 1,$ we note from
the monotonicity of $k \rightarrow k t_k^2$ (cf. \eqref{eq:nutnur})
and the definition \eqref{eq:kappaconsequ} of $\kappa_n$ that
\begin{displaymath}
   k_- t_-^2 \leq k_+ t^2[k_+] \leq \kappa_n t^2[\kappa_n]
             \sim (1 - q_n)^{-1} k_n \tau_\eta^2,
\end{displaymath}
which leads to the second term in the upper bound of
\eqref{eq:mainbound_mp}.

\bigskip

\textit{Negligibility of $T_{1n} - T_{4n}$.}  Consider first the term
$T_{1n}(\mu)$. For the non-zero $\mu_l$ (of which there are at most $n
\eta_n$), use the bound \eqref{cbound} to get
\begin{displaymath}
   T_{1n}(\mu) \leq k_n (t_1 + 1) \leq k_n \tau_\eta.
\end{displaymath}

For the large signal component term $T_{2n}$, we have, using Lemma
\ref{lem:randthresh} and the bound \( \hat{t}_2  \leq t_{1} =
O(\log^{1/2} n) \),
\begin{equation}
   \label{eq:t2nbd}
         T_{2n} \leq  \sum_{S_n^c(\mu)} E [\hat{e}_i, A_{n}] \leq
         2 \sum_{S_n^c(\mu)} \{ E [ \eta(y_i,\hat{t}_2) - \mu_i ]^2 \}^{1/2}
         \leq c_0 t_1 | S_n^c(\mu) |
          \leq | S_n^c(\mu) | (c_0 \log n )^{1/2}.
\end{equation}
On $\ell_0[\eta_n],$ clearly $| S_n^c(\mu) | \leq n \eta_n$ and so
$ T_{2n}(\mu) \leq c_1 k_n \tau_\eta.$

For the small threshold term \( T_{3n} \), note first that
\( \sum \hat{e}_{i} \leq 2 \| z \| \| \hat{\mu}_2  - \mu \|, \)
so that
\begin{displaymath}
         T_{3n}(\mu) \leq 2 ( E \| z \|^{2} )^{1/2} \{ E [ \| \hat{\mu}_2  -
         \mu \|^{2}, A_{n}^{c}] \}^{1/2}.
\end{displaymath}
Now \( E \| z \|^{2} = n \), and since \( A_{n}^{c} \) is a rare
event, apply Lemma \ref{lem:randthresh}, noting
the bound \( \hat{t}_2  \leq t_1 . \) Thus
\begin{displaymath}
          T_{3n}(\mu)  \leq 8 n t_1  P_{\mu}(A_{n}^{c})^{1/4} \leq
                         c_{1} n t_1  \exp \{ - c_{2} \log^2 n \}
                         = o( R_n(\ell_0[\eta_n]) ).
                 \end{displaymath}
uniformly on \( \ell_0[\eta_{n}] \) after applying Proposition
\ref{prop:sandwichbd}.

The remaining term \( T_{4n} \) is handled exactly as was \( T_{3n}
\): if we let \( \hat{\mu}_{F} \) denote hard thresholding at the
(fixed) threshold \( t_-  \), then \( \sum_{S_{n}(\mu)} e_{i} \leq
\| z \| \| \hat{\mu}_{F} - \mu \| \); and now Lemma \ref{lem:randthresh}
and Proposition \ref{prop:sandwichbd} can be used as before.

\bigskip

\textbf{Weak $\ell_p$.}  A little extra work is required to analyze
$T_{1n}(\mu)$, so we first dispose of $T_{2n}-T_{4n}.$ The analysis
of $T_{3n}$ and $T_{4n}$ is essentially as above. For $T_{2n},$ we
bound  $| S_n^c(\mu) |$ using the extremal element of $m_p[\eta_n]$,
namely $\bar{\mu}_l = \eta_n (n/l)^{1/p}.$
Thus, for all $\mu \in m_p[\eta_n],$
\begin{displaymath}
   S_n^c(\mu) \subset \{ l : \bar \mu_l > t_-(\mu) \},
\end{displaymath}
and since, for $\eta$ sufficiently small,
\begin{equation}
   \label{eq:tminusbd}
   t_-(\mu) = t[k_+(\mu)] \geq t[k_+(\bar \mu)] \geq \tau_\eta - 3/2
    = \tau_\eta( 1 + o(1) ),
\end{equation}
we have
\begin{equation}
     \label{eq:sncmu-lp}
           | S_n^c(\mu) |  \leq  n \eta_n^p t_-^{-p}(\mu)
                           \leq  n \eta_n^p \tau_\eta ^{-p}( 1 + o(1)).
\end{equation}

 From \eqref{eq:lpmmxrisk} we have $R_n = R(\ell_p[\eta_n]) \sim n \eta_n^p
\tau_\eta^{2-p}$, and so from \eqref{eq:t2nbd} and
\eqref{eq:sncmu-lp}, we get
\begin{displaymath}
   T_{2n}(\mu) \leq c_0 t_1 n \eta_n^p \tau_\eta^{-p} \leq c k_n \tau_\eta.
\end{displaymath}

For the $T_{1n}$ term, we obtain from Lemma \ref{lem:cprop} that
\begin{displaymath}
   \xi(t,\mu) - \xi(t,0) \leq c (\mu^2 \wedge t),
\end{displaymath}
and so $T_{1n}(\mu) \leq 2 c \sum \bar \mu_l^2 \wedge t_1.$
The negligibility of $T_{1n}$ is a consequence of the following.

\begin{lemma}
   \label{lem-t1bd}
   For $0 < p < r \leq 2$, we have
   \begin{displaymath}
     \sum_l \bar \mu_l^r \wedge t_1^{(r-1)_+} = o( k_n \tau_\eta^r).
   \end{displaymath}
\end{lemma}
\begin{proof}
   Define $\tilde k$ by $\bar \mu_{\tilde k} = \tau_\eta/ \log
   \tau_\eta$ so that $\tilde k = n \eta_n^p \tau_\eta^{-p} \log^p
   \tau_\eta$. We have
   \begin{displaymath}
     \sum_{l=1}^n \bar \mu_l^r \wedge t_1^{(r-1)_+}
       \leq \tilde k t_1^{(r-1)_+} + \sum_{l > \tilde k} \bar \mu_l^r.
   \end{displaymath}
Since $t_1 \leq c \tau_\eta$,
\begin{displaymath}
   \tilde k t_1^{(r-1)_+} \leq c k_n \tau_\eta^{(r-1)_+} \log^p
   \tau_\eta = o( k_n \tau_\eta^r).
\end{displaymath}
By integral approximation, since $0 < p < r$,
\begin{align*}
   \sum_{l > \tilde k}  \bar \mu_l^r
& \leq n \eta_n^r \int_{\tilde k/n}^1 x^{-r/p}
   dx \leq c_{rp} n \eta_n^r (\tilde k/n)^{1 - r/p} \\
      & = c_{rp} n \eta_n^p \tau_\eta^{r-p} \log^{p-r} \tau_\eta
        = o( k_n \tau_\eta^r).    \qedhere
\end{align*}
\end{proof}

\section{$\ell_{r}$ losses}\label{sec.ellr}
\setcounter{equation}{0}

This section retraces for $\ell_r$ loss the steps used for squared
error in Section \ref{sec:proof}, making adjustments for the fact that
the quadratic decomposition \eqref{eq:kdecomp} is no longer available.
It turns out that this decomposition is merely a convenience - the
asymptotic result of Theorem \ref{th:mainres} is as sharp for all
$0 < r \leq 2.$ However the analysis of the error term is more complex
than in Section \ref{sec:errtermproof}, requiring bounds developed in
Lemmas \ref{lem:b-of-a} and
\ref{lem:crbd}.

\subsection{Empirical Complexity for $\ell_r$
loss}\label{subsec.emp-comp-r}

For an $\ell_{r}$ loss function, we use a modified empirical complexity
\begin{displaymath}
           K(\tilde{\mu},y;r) = \| y - \tilde{\mu} \|_r^r +
           \sum_{l=1}^{N(\tilde{\mu})} t_l^r.
\end{displaymath}
The minimizers of empirical and theoretical complexity are defined,
respectively, by
\begin{align*}
           \hat{\mu}_r & = \argmin_{\tilde{\mu}} K(\tilde{\mu},y;r),  \\
           \mu_0       & = \argmin_{\tilde{\mu}} K(\tilde{\mu},\mu;r) .
\end{align*}
For $\ell_r$ loss, the quadratic decomposition of \eqref{eq:kdecomp}
is replaced by
\begin{equation}
           K(\tilde{\mu},\mu+z) = K(\tilde{\mu},\mu) + \| \mu - \tilde{\mu} + z
           \|_r^r - \| \mu - \tilde{\mu} \|_r^r .
           \label{eq:kdecomp-r}
\end{equation}
The key inequality
\begin{displaymath}
           K(\hat{\mu}_r,y) \leq K(\mu_0,y),
\end{displaymath}
when combined with \eqref{eq:kdecomp-r},  applied to both
$\tilde{\mu} = \hat{\mu}_r$ and $\mu_0$,  yields the analog of
\eqref{complexbound}:
\begin{displaymath}
           K(\hat{\mu}_r,\mu) \leq K(\mu_0,\mu) + D(\hat{\mu}_r,\mu_0,\mu,y).
\end{displaymath}
Setting
\begin{equation}
           \hat{\delta} = \mu - \hat{\mu}_r, \quad
           \delta_0 = \mu - \mu_0, \quad y = \mu + z,
           \label{eq:deviations}
\end{equation}
we have for the error term
\begin{equation}
           \hat{D} = D(\hat{\mu}_r,\mu_0,\mu,y) =
           \| \delta_0 + z \|_r^r - \| \delta_0 \|_r^r - \| \hat{\delta}  + z
\|_r^r
           + \| \hat{\delta} \|_r^r = \sum_1^n \hat{d}_l
           \label{eq:dhat}
\end{equation}
Thus, with $Err_r \equiv E_\mu
           \hat{D} - E_\mu Pen_r(\hat \mu_r)$,
\begin{equation}
           E \| \hat{\mu}_r - \mu \|_r^r \leq K(\mu_0,\mu; r) + 
Err_r(\mu,\hat{\mu_r}).
           \label{eq:lr-lossbd}
\end{equation}


\subsection{Maximum theoretical complexity}\label{subsec.maxtheor-r}

The theoretical complexity corresponding to $\mu$ is given by
\begin{equation}
                   K(\mu_0,\mu;r) = \inf_k \sum_{k+1}^{n} |\mu|_{(l)}^r
+ \sum_1^k t_l^r.
           \label{eq:theor-r}
\end{equation}
We may argue as in Section \ref{ssec:max-theoretical-complexity}
that for $\Theta_n = \ell_0[\eta_n]$,
\begin{equation}
   \label{eq:rpowasy}
           \sup_{\mu \in \Theta_n} K(\mu_0,\mu) \leq
           \sum_1^{k_n} t_l^r
           \sim k_n t_{k_n}^r
           \sim n \eta_n ( 2 \log \eta_n^{-1} )^{r/2} \sim
           R_n (\Theta_n;r).
\end{equation}
and that for $\Theta_n = m_p[\eta_n]$,
\begin{displaymath}
     \sup_{\mu \in \Theta_n} K(\mu_0,\mu) = \inf_k \ C_n^r \sum_{k+1}^n
l^{-r/p} +
     \sum_1^k t_l^r \sim R_n(\Theta_n;r).
\end{displaymath}
Finally, for $\Theta_n = \ell_p[\eta_n]$, we may 
argue that
\begin{displaymath}
     \sup_{\mu \in \Theta_n} K(\mu_0,\mu) \sim  \max_k  \{
     \sum_1^k t_l^r :   \sum_1^k t_l^p \leq n \eta_n^p \}   \sim 
R_n(\Theta_n;r).
\end{displaymath}

We remark that if $\underline{k}(\mu)$ is an index minimizing 
\eqref{eq:theor-r},
then $\mu_{0i}$ is obtained from hard thresholding of $\mu_i$ at
$t_0 = t[\underline{k}(\mu)]$ (interpreted as $t_1$ if $k(\mu)=0$).
In any event, this implies
\begin{equation}
           | \delta_{0i} | = | \mu_i - \mu_{0i} | \leq | \mu_i | \wedge t_1
           \label{eq:del0ibd}
\end{equation}

\subsection{The $\ell_r$ error term}\label{subsec.lr-error-term}

There is an $\ell_r$ analog of bound
\eqref{monotone}; this allows us to replace the random threshold
$\hat{t}_r$ by the fixed threshold value $t_\kappa $ for the most important
cases.

Indeed, suppose that $|\mu_i| \leq t_\kappa  \leq \hat{t}_r$.
Let $\bar{\mu}_i(y) = \eta_H(y_i,t_\kappa )$ denote hard thresholding at
$t_\kappa ,$ and let $\bar{\delta}_i = \mu_i - \bar{\mu}_i$ denote the
corresponding deviation. We claim that
\begin{equation}
           |\hat{\delta}_i|^r - |\hat{\delta}_i + z_i|^r \leq
           |\bar{\delta}_i|^r - |\bar{\delta}_i + z_i|^r
           \label{eq:errorupper}
\end{equation}
Indeed, $\hat{\delta}_i = \bar{\delta}_i$ unless $t_\kappa  \leq |y_i| \leq
\hat{t}_r.$ In this case, we have $\hat{\mu}_i = 0$ so that
$\hat{\delta}_i = \mu_i$ while $\bar{\mu}_i = y_i$ so that
$\bar{\delta}_i = -z_i.$ In this case, \eqref{eq:errorupper} reduces
to
\begin{displaymath}
           |\mu_i|^r - |y_i|^r \leq |z_i|^r
\end{displaymath}
which is trivially true since $|\mu_i| \leq t_\kappa    \leq |y_i| .$

We now derive the $\ell_r$ analog of the error decomposition
\eqref{eq:err-decomp}.
Recalling the notation \eqref{eq:deviations} - \eqref{eq:lr-lossbd},
we have
\begin{equation}
     \label{eq:dihatdef}
\hat{d}_i = d_i(\hat{t}_r) = |\delta_{0i} + z_i|^r - |\delta_{0i}|^r
- |\hat{\delta}_i + z_i|^r + |\hat{\delta}_i |^r.
\end{equation}
Defining as in Section \ref{ssec:error-term}
    the sets $A_n = \{ t_- \leq \hat{t}_r \leq t_+  \}$ and
$S_n(\mu) = \{ i: | \mu_i | \leq t_-  \},$ we obtain
\begin{align*}
           E_\mu \hat{D}   =  E \sum \hat{d}_i
            & = E[ \sum_{S_n(\mu)} \hat{d}_i, A_n ] + E [ \sum_{S_n^c(\mu)}
           \hat{d}_i, A_n ] + E [ \sum \hat{d}_i, A_n^c] \\
           & = D_{an} + T_{2n} + T_{3n}.
\end{align*}
Let $d_i = d_i(t_-)$: the monotonicity of errors for small
components (compare \eqref{eq:errorupper}) says that the leading term
           \begin{displaymath}
                   D_{an} \leq E[ \sum_{S_n(\mu)} d_i, A_n ]
                    = E[ \sum_{S_n(\mu)} d_i ] - E[ \sum_{S_n(\mu)} d_i, A_n^c ]
            = D_{bn} + T_{4n}.
           \end{displaymath}
Consider first the dominant term $D_{bn}.$ First, write
\begin{align}
           E d_i & = E_\mu |\delta_{0i} + z_i|^r - |\delta_{0i}|^r
- |\bar{\delta}_i + z_i|^r + |\bar{\delta}_i |^r \notag  \\
             & = \psi_r(\delta_{0i}) + \xi_r(t_- ,\mu_i),  \label{eq:edi}
\end{align}
where, for $y=\mu + z,$ \ $z \sim N(0,1)$, and $0 < r \leq 2,$ we define
\begin{align}
           \psi_r(a) & = E[ \ |a+z|^r - |a|^r - |z|^r ] ,  \label{eq:bra} \\
           \xi_r(t,\mu) & = E[ \ | \eta_H(y,t) - \mu|^r - | 
\eta_H(y,t) - y |^r +
           |y-\mu|^r].   \label{eq:crt-mu}
\end{align}
[Note that a term $E |z|^r$ has been introduced in both $\psi_r$ and
$\xi_r$ -- as a result $\psi_2(a) \equiv 0$ and $\xi_2(t,\mu) = 2 \xi(t,\mu)$
as defined at \eqref{cformula}.]  The next lemmas, proved in
Section \ref{subsec.proofs-r-nblack} below, 
play the same role as Lemma
\ref{lem:cprop} for the $\ell_2$ case.

\begin{lemma}  \label{lem:b-of-a}
           The function $\psi_r(a)$ defined at \eqref{eq:bra} is even in $a$ and
           \begin{equation}
                                   |\psi_r(a)| \leq    \begin{cases}
                                           C_1 |a|^r &  \text{for all } a   \\
                                   C_2 |a|^{(r-1)_+}  & \text{for } |a|
\ \text{large}.
                                   \end{cases}
           \label{eq:b-of-a}
\end{equation}
\end{lemma}

\begin{lemma}  \label{lem:crbd}
           The function $\xi_{r}(t,\mu)$ defined at \eqref{eq:crt-mu} is even in
           $\mu$ and satisfies
           \begin{align}
                           \xi_r(t,0) & = 2 \int_{|z| > t} |z|^r \phi(z)
dz,         \label{eq:crbd} \\
                           |\xi_r(t,\mu)| & \leq C [ t^{(r-1)_+} + 1]
\qquad \qquad \mu \in
                           \mathbb{R}, t>0  . \label{eq:crbound}
       \end{align}
\end{lemma}

With the preceding notation, we may therefore write
\begin{displaymath}
    D_{bn}(\mu)  = \sum_{S_n(\mu)} \psi_r(\delta_{0i}) + \xi_r(t_-,\mu_i)
      = |S_n(\mu)| \xi_r(t_- ,0) + T_{1n}(\mu),
\end{displaymath}
say.  To summarize, we obtain the following decomposition for the error term
in \eqref{eq:lr-lossbd}:
\begin{align*}
   E_\mu \hat D - E_\mu Pen_r(\hat \mu)
      & \leq D_{cn} (\mu) + \sum_{j=1}^4 T_{jn}(\mu)  \\
   D_{cn}(\mu)  & =  n \xi_r(t_-,0) - E_\mu Pen_r(\hat \mu_r).
\end{align*}

\textit{Dominant term.}
Using \eqref{eq:hardriskat0},
\begin{displaymath}
   \xi_r(t,0)
   \leq 4 t^r \tilde \Phi(t) [1 + 2 t^{-2}].
\end{displaymath}
Since $2 \tilde \Phi(t_-) = q_n k_+ n^{-1}$, we obtain
\begin{displaymath}
   n \xi_r(t_-,0) \leq 2 q_n k_+ t_-^r + 8 q_n k_+ t_-^{r-2}
        \leq 2 q_n k_- t_-^r + c k_n \tau_\eta^{r-1},
\end{displaymath}
since $k_+ - k_- \leq c k_n \tau_\eta^{-1}$ and $t_-(\mu) \asymp
\tau_\eta$.

For the penalty term, arguing as before yields
\begin{displaymath}
   E_\mu Pen_r(\hat \mu_r) \geq k_- t_-^r + O( k_n \tau_\eta^{r-1}).
\end{displaymath}
Combining the two previous displays, we obtain
\begin{displaymath}
   D_{cn}(\mu) \leq (2 q_n - 1) k_- t_-^r + O( k_n \tau_\eta^{r-1})
\end{displaymath}

If $q_n \leq 1/2$, of course the leading term is non-positive, while in
the case  $1/2 \leq q_n < 1$,  we note from \eqref{eq:nutnur}
and the definition \eqref{eq:kappaconsequ} of $\kappa_n$ that
\begin{displaymath}
   k_-(\mu) t_-^r(\mu) \leq k_+ t^r[k_+]
      \leq \kappa_n t^r[\kappa_n] \sim  (1 - q_n)^{-1} k_n \tau_\eta^r,
  \end{displaymath}
which shows that $D_{cn}(\mu)$ is bounded by
the second term  in the upper bound of
(\ref{eq:mainbound_mp}).

\bigskip

\textit{Negligibility of $T_{1n} - T_{4n}$}.
Since there are at most $k_n$ non-zero terms in $T_{1n}$, from Lemmas
\ref{lem:b-of-a} and \ref{lem:crbd}, we obtain
\begin{displaymath}
   T_{1n}(\mu) \leq C k_n t_1^{(r-1)_+} = o(k_n \tau_\eta^r).
\end{displaymath}

To bound $T_{2n}$, we first note from the properties of hard
thresholding (compare \eqref{threshbd}) that
\begin{displaymath}
           |\hat{\delta}_i| = |\hat{\mu}_{r,i} - \mu_i | \leq \hat{t}_r + |z_i|
           \leq t_1 + |z_i|.
\end{displaymath}
Inequality \eqref{eq:del0ibd} shows that $| \delta_{0i}| \leq t_1$.

Combined with \eqref{eq:dihatdef} and \eqref{eq:auxilbd}, this
shows, for $1 < r \leq 2$
\begin{equation}
     | \hat{d}_i | \leq 3 t_1^{r-1} |z_i| + 6 |z_i|^r,
           \label{eq:dihatbd}
\end{equation}
while for $0 < r \leq 1$ only the $|z_i|^r$ term is needed.
Consequently, there exist constants $C_i$ such that for $0 < r \leq 2$
\begin{equation}
   \label{eq:di2bd}
           E |\hat{d}_l| \leq C_1 t_1^{(r-1)_+},  \qquad \text{and}
   \qquad  E \hat d_i^2 \leq C_2 t_2^{2(r-1)_+}.
\end{equation}
Thus
\begin{equation}
   \label{eq:t2nlr}
           |T_{2n}| \leq \sum_{S_n^c(\mu)} E | \hat{d}_i | \leq C | S_n^c(\mu) |
           t_1^{(r-1)_+}.
\end{equation}
And so on $\ell_0[\eta_n]$,
\begin{displaymath}
           |T_{2n}| \leq  C n \eta_n t_1^{(r-1)_+} = o( k_n \tau_\eta^r).
\end{displaymath}

To bound $T_{3n}(\mu)$, use \eqref{eq:di2bd} and Cauchy-Schwartz:
\begin{displaymath}
   T_{3n}(\mu) \leq P(A_n^c)^{1/2} \sum_i ( E \hat d_i^2 )^{1/2}
               \leq c n t_1^{(r-1)_+} \exp \{ - c_0 ( \log n)^2/2 \}
               = o(R_n(\Theta_n)),
\end{displaymath}
since $A_n$ is $\Theta_n-$likely.
Argument similarly for $T_{4n}$,
with threshold at $t_-(\mu) $ instead of $\hat{t}_r$.


\textit{Weak $\ell_p$.}
For the $T_{1n}$ term, we use a consequence of Lemmas
\ref{lem:b-of-a} and \ref{lem:crbd}, proved in 
Section \ref{subsec.proofs-r-nblack}, 
to bound the summands in $T_{1n}(\mu)$:
\begin{equation}
     \psi_r(\delta_{0i}) + |\xi_r(t_- ,\mu_i) - \xi_r(t_- ,0) | \leq C [ |
           \mu_i|^r \wedge t_1^{(r-1)_+}]
           \label{eq:T0-err}
\end{equation}
Combined with Lemma \ref{lem-t1bd}, this shows that
$\sup_{m_p[\eta_n]} T_{1n}(\mu) = o(k_n \tau_\eta^r)$.

For $T_{2n}$, we use  \eqref{eq:sncmu-lp} in \eqref{eq:t2nlr} to
obtain
\begin{displaymath}
   |T_{2n}| \leq C n \eta_n^p t_\eta ^{r-p} t_1^{(r-1)_+ -r} =
      o(R_n(\Theta_n)).
\end{displaymath}
The analysis of $T_{3n}$ and $T_{4n}$ is as for the $\ell_0$ case.

\subsection{From penalized to original FDR}
\label{sec:pen-to-orig-again}

\textit{Proof of Theorem \ref{th:pen-to-fdr}.}
Apply Lemma \ref{lem:riskdel} with $\hat \mu = \hat \mu_D$ and
$\hat \mu^\prime$ and $t = t_1$.
We abbreviate the minimax risk $R_n(\Theta_n)$ by $R_n$.
Theorem \ref{th:mainres} shows that for sufficiently large $n$,
$\sup_{\Theta_n} \rho (\hat \mu_r, \mu) \leq c_0 R_n$ so that the bound
established by Lemma \ref{lem:riskdel} yields
\begin{equation}
   \label{eq:riskdelta}
   \sup_{\mu \in \Theta_n} | \rho( \hat \mu_D, \mu) - \rho( \hat \mu_r,
   \mu )|
    \leq 2 \beta_n t_1^r + 2 c_0 I\{r>1 \} R_n^{1 - 1/r} (\beta_n
   t_1^r)^{1/r} + 8 n t_1^r \sup_{\Theta_n} P_\mu (B_n^c)^{1/2}.
\end{equation}

The thresholds $\hat t_D$ and $\hat t_r$ corresponding to $\hat \mu_D$
and $\hat \mu_r$ both lie in $[\hat t_G, \hat t_F]$, and so
with probability one,
\begin{displaymath}
   N'  = \# \{ i : |y_i| \in [ \hat t_D,  \hat t_r] \}
      \ \leq \ \# \{ i: \hat t_F \leq |y_i| < \hat t_G \}
       \leq \hat k_F - \hat k_G.
\end{displaymath}
[The first inequality is valid except possibly on a zero probability
event in which some $|y_i| = \hat t_G$.
To see the second inequality, note from \eqref{eq:kFdef} that
$l > \hat k_F$ implies $Y_l \leq t_l < \hat t_F$, while
\eqref{eq:kGdef} entails that $l \leq \hat k_G$ implies
$Y_l \geq t_l \geq \hat t_G$.
Consequently $\hat t_F \leq Y_l < \hat t_G$ implies $\hat k_G < l \leq
\hat k_F$ which yields the required inequality.]

If we take $\beta_n = 3 \alpha_n k_n$, then Proposition
\ref{prop:sandwichbd} implies that
\begin{displaymath}
   P_\mu (B_n^c) = P_\mu(N' > \beta_n)
      \leq P_\mu( \hat k_F - \hat k_G  > 3 \alpha_n k_n)
      \leq c_0 \exp \{ - c_1 \log^2 n \},
\end{displaymath}
so that the third term in \eqref{eq:riskdelta} is $o(R_n)$.
Finally $\beta_n t_1^r \leq c \alpha_n k_n \tau_\eta^r \leq c \alpha_n
R_n$ so that the first two terms are also $o(R_n)$, completing the proof.

\subsection{Remaining proofs for the $\ell_r$ case}
\label{subsec.proofs-r-nblack}

\begin{proof}[Proof of Lemma \ref{lem:b-of-a}.]
           The function $a \rightarrow |z+a|^r$ is H\"older$(r)$ with constant
           $C=2$ uniformly in $z \in \mathbb{R}$ and $r \in (0,2]$.
           Consequently, so is $f(a) = E| Z + a|^r.$ Since $f'(0)=0$,
           \begin{displaymath}
                   \psi_r(a) = \begin{cases}
                           f(a) - f(0) - |a|^r & r \leq 1 \\
                           \int_0^{a} [ f'(s) - f'(0)] ds - |a|^r  & 1 < r
\leq 2,
                   \end{cases}
           \end{displaymath}
           and the global bound in \eqref{eq:b-of-a} follows from the
           H\"older properties of $f$.

           For $a$ large, $\int_{- \infty}^{-a} [|z+a|^r + |z|^r]
\phi(z) dz = o(1)$
           and so it suffices to consider
           \begin{displaymath}
                   \int_{-a}^{\infty} |(a+z)^r - a^r| \phi(z) dz \leq
                   a^r \int_{-a}^a  \frac{C|z|}{a} \phi(z) dz +
\int_a^{\infty} (2z)^r
                   \phi(z) dz \leq C a^{r-1},
           \end{displaymath}
using \eqref{eq:hardriskat0}.
         Thus, for $a$ large,
                   $|\psi_r(a)| \leq C a^{r-1} + E|Z|^{r} \leq C_2 a^{(r-1)_+}.$
\end{proof}

\bigskip

\begin{proof}[Proof of Lemma \ref{lem:crbd}]
           Making the threshold zone explicit leads to
           \begin{equation}
                   \xi_r(t,\mu) = \int_{|y|<t} [|y-\mu|^r - |y|^r +  |\mu|^r ]
\ \phi(y-\mu)
                   dy + 2 \int_{|y|>t} |y-\mu|^r \phi(y-\mu) dy.
           \label{eq:cr-thresh}
       \end{equation}
 From this, evenness and \eqref{eq:crbd} are straightforward.
           Evenness means that we may take $\mu \geq 0.$
           Let $c_r = \int |y-\mu|^r \phi(y-\mu) dy,$ and note that
for $y > 0$,  $\phi(-y-\mu) \leq \phi(y-\mu).$ Thus
                   \begin{equation}
                   |\xi_r(t,\mu)| \leq 2 c_r + 2 \int_0^t |\mu^r - y^r| \
\phi(y-\mu) dy.
                   \label{eq:inter1}
           \end{equation}
           If $r \leq 1$, then $|\mu^r - y^r| \leq |y-\mu|+1$ as follows by
           checking cases with $y,\mu \in [0,1]$ and $[1,\infty]$ respectively.
           Hence $|\xi_r(t,\mu)| \leq 2 c_r + 3$ and \eqref{eq:crbd} follows.

           Now suppose $1 < r \leq 2,$ and that $0 \leq \mu \leq t$. Break the
           integral in \eqref{eq:inter1} into two pieces. Consider first $0 \leq
           y \leq \mu$:
           \begin{align}
                   \int_0^\mu (\mu^r - y^r) \ \phi(y-\mu) dy & =
                   r \int_0^{\mu} dv \, v^{r-1} \int_0^v \phi(y-\mu) dy
                  \label{8-18a} \\
                   & \leq r \int_0^{\mu} (\mu - w)^{r-1} \tilde{\Phi}(w) dw
\qquad
                   \qquad (w = \mu - v) \notag \\
                   & \leq r \mu^{r-1} \int_0^{\infty} \tilde{\Phi}(w) dw
= c \mu^{r-1}
                   \leq c t^{r-1}. \notag
           \end{align}
           Arguing similarly on the interval $\mu \leq y \leq t$,
                   \begin{align*}
                   \int_\mu^t (y^r - \mu^r) \phi(y-\mu) dy & =
                   r \int_{\mu}^t dv \, v^{r-1} \int_v^t \phi(y-\mu) dy \\
                   & \leq r \int_0^{t-\mu} (\mu+w)^{r-1} \tilde{\Phi}(w) dw
                   \leq c t^{r-1}.
           \end{align*}
           Combining the last two displays with \eqref{eq:inter1} yields
\eqref{eq:crbd}.

           Finally suppose that $\mu > t.$ In \eqref{eq:inter1} we bound
                   \begin{displaymath}
                   \int_0^t  (\mu^r - y^r) \ \phi(y-\mu) dy  \leq
                   (\mu^r - t^r) \Phi(t-\mu) + \int_0^t (t^r-y^r) \phi(y-t) dy.
           \end{displaymath}
           The last integral is bounded as for the integral $0 \leq y \leq \mu$
           at \eqref{8-18a} above (setting $\mu = t$). Finally, write
$\mu = t+x$: since
                   \begin{displaymath}
                   (t+x)^r - t^r \leq
                                   \begin{cases}
                           2^r x^r & \text{if} \ t \leq x \\
                           r \int_t^{t+x} v^{r-1} dv \leq C
(2t)^{(r-1)_+}x  & \text{if} \ t \geq x
                   \end{cases}
           \end{displaymath}
           and $x \rightarrow x^a \tilde{\Phi}(x)$ is bounded
           it follows that for some $C,$  $[ (t+x)^r - t^r]
\tilde{\Phi}(x) \leq C t^{(r-1)_+}.$
\end{proof}

\bigskip

\begin{proof}[Proof of \eqref{eq:T0-err}]
           First, recall from \eqref{eq:del0ibd} that $| \delta_{0i} |
\leq \mu_i \wedge t_1.$
Lemma \ref{lem:b-of-a} implies the existence of a constant $C = C(r)$
such that $\psi_r(\Delta_{0i}) \leq C |\mu_i|^r \wedge t_1^{(r-1)_+}.$
In view of \eqref{eq:crbound}, it remains to show that
\begin{equation}
           | \xi_r(t,\mu) - \xi_r(t,0) | \leq C | \mu |^r.
           \label{eq:dev-bound}
\end{equation}
Rearranging \eqref{eq:cr-thresh} leads to
\begin{displaymath}
           \xi_r(t,\mu) = c_r - \int_{|y| < t} |y|^r \phi(y-\mu) dy + |\mu|^r [
           \Phi(t-\mu) - \Phi(-t - \mu) ] +
           \int_{|z+\mu| > t} |z|^r \phi(z) dz.
\end{displaymath}
Write $\xi_{r}(t,\mu) - \xi_r(t,0) = U_1(\mu) + U_2 (\mu) + U_3 (\mu)$
corresponding to the three terms above.
Consider first $U_1$.  If $|\mu| \geq 1,$ then
$|U_1(\mu)| \leq \xi_r + \mu^r \int | 1 + z/\mu |^r \phi(z) dz \leq c
|\mu|^r.$
For $|\mu| \leq 1,$  $U_1(\mu)$ is $C^2$ with $U_1(0) =
U_1^{\prime}(0) = 0$ and
$U''$ uniformly bounded in $t$. Hence $|U_1(\mu)| \leq C \mu^2 \leq
C |\mu|^r$ for $|\mu| \leq 1$ also.
Clearly $|U_2(\mu)| \leq C_2 |\mu|^r,$ while
$U_3 (\mu)$ is bounded and $C^1$ with
$U_3^{\prime}(\mu) = |t-\mu|^r \phi(t-\mu) - |t+\mu|^r \phi(t+ \mu).$
Thus $U_3^{\prime}(0) = 0$ and for $|\mu| \leq 1$ and $t \geq 2$,
$U_3^{\prime \prime} \leq C_3$.
Hence $U_3(\mu) \leq C_4 (\mu^2 \wedge 1) \leq C_5 |\mu|^r,$ and
\eqref{eq:dev-bound} is proved.
\end{proof}

\section{Gaussian tails and quantiles}
\label{sec:gauss-quantiles}
\setcounter{equation}{0}

We collect in this Appendix some results about the normal density 
$\phi$,
the normal CDF $\Phi$ and the normal quantile function $z()$; 
these have
been used extensively above.  

\begin{lemma}  \label{lem:gausstails}
We have
\begin{alignat}{2}
   \phi(v) & \leq (v + 2 v^{-1}) \tilde \Phi(v),
   & \qquad \qquad & v \geq \sqrt{2}.   \label{eq:phiupper} \\
\intertext{More generally, Mills' ratio $M(y) = y \tilde{\Phi}(y)/ \phi(y)$
increases from $0$ to $1$ as $y$ increases from $0$ to $\infty$. In
particular,}
\phi(v)/(2v) \leq \tilde \Phi(v) & \leq \phi(v)/v,
& \qquad \qquad & v \geq 1,    \label{eq:mills} \\
2 \tilde \Phi (v) & \leq e^{-v^2/2},
& \qquad \qquad & v \geq 0,       \label{eq:global} \\
\tilde \Phi( v - c/v) & \leq 4 e^c \tilde \Phi(v),
& \qquad \qquad & v \geq \sqrt{2c}.   \label{eq:increment}
\end{alignat}
\end{lemma}
\begin{proof}
Bound \eqref{eq:phiupper} follows from the standard lower bound
$\tilde \Phi(v) \geq (v^{-1} - v^{-3}) \phi(v)$, compare
\citet[p. 175]{fell68}, and $(1 - v^{-2})^{-1} \leq 1 + 2 v^{-2}$ for
$v^2 \geq 2.$
Monotonicity of Mills' ratio follows from differentiation and partial
integration, and then evaluating
$M(1) = \tilde{\Phi}(1)/\phi(1) \doteq 0.655 > 1/2$
yields \eqref{eq:mills}.
The difference $g(v) = 2 \tilde \Phi(v) - e^{-v^2/2}$
has $g'(v) = [2 \phi(0) - v] e^{-v^2/2} \leq 0$ for $v \leq 2
\phi(0)$
and vanishes at $0$.
But when $v \geq 2 \phi(0)$, \eqref{eq:global} follows from
\eqref{eq:mills}. Finally, \eqref{eq:increment} follows by applying the right
and then left sides of \eqref{eq:mills}, noting that $t - c/t \geq
t/2$ and that $\phi(t - c/t) \leq \phi(t) e^c.$
\end{proof}

\begin{lemma}
   \label{lem:logbound}
Suppose that $k$ and $\alpha$ are such that
$\max \{ \alpha, 1/\alpha \} \leq  C \log k$. Then
\begin{equation}
   \label{eq:logbound}
   \sqrt{2 \log k \alpha} = \sqrt{2 \log k} + \theta \sqrt{C},
\end{equation}
and if $\alpha \geq 1$, then $0 \leq \theta \leq \sqrt{2}/e \leq 1$,
while if $\alpha \leq 1$, then $ -1.1 \leq - \sqrt{8}/e \leq \theta
\leq 0$.
\end{lemma}
\begin{proof}
   If $\alpha \geq 1,$ then from the inequality
$(1+x)^{1/2} \leq 1 + x/2$, valid for nonnegative $x$,
\begin{displaymath}
   0 \leq   \sqrt{2 \log k \alpha} - \sqrt{2 \log k}
   \leq \frac{\log \alpha}{\sqrt{2 \log k}}
   \leq \sqrt{\frac{C}{2}} \frac{\log (C \log k)}{\sqrt{C \log k}}
   \leq \frac{\sqrt{2}}{e} \sqrt{C},
\end{displaymath}
since $u \rightarrow u^{-1/2} \log u \leq 2/e $ for $u \geq 1.$
When $\alpha \leq 1,$ essentially the same argument applies, with
inequalities reversed, and using $(1+x)^{1/2} \geq 1 + x$, valid for
$-1 \leq x \leq 0.$
\end{proof}
As an immediate corollary, we note that if $k_n = b_n n^{1-\beta}$ and
$\max \{ b_n, 1/b_n \} \leq c \log^\rho n$ for $0 \leq \rho \leq 1$,
then
\begin{equation}
   \label{eq:logbd}
   | \sqrt{2 \log (n/k_n)} - \sqrt{2 \beta \log n}| \leq 1.1
   ( c^{-1} \beta \log^{1-\rho} n )^{-1/2}.
\end{equation}

\begin{lemma} \label{lem:quantile}
(1) Let $z(\eta) = \tilde{\Phi}^{-1}(\eta)$ denote the upper
$(1-\eta)^{\text{th}}$ percentile of the Gaussian distribution. If $\eta \leq
.01,$ then
\begin{align}
  z^{2}(\eta) & =   2 \log \eta^{-1} - \log \log \eta^{-1} - r_2(\eta), &
                 r_2(\eta) \in [1.8, 3], \label{eq:quantilebracket}
                 \\
  z(\eta)  & =  \sqrt{2 \log \eta^{-1}} - r_1(\eta)  &
                 r_1(\eta) \in [0, 1.5].   \label{eq:zetabd}
\end{align}

(2) We have $z'(\eta) = -1/\phi(z(\eta))$, 
and hence if $\eta_2 > \eta_1 > 0,$ then
\begin{equation}
       z(\eta_1) - z(\eta_2) \leq \frac{\eta_2 - \eta_1}{\eta_1 z(\eta_1)}.
       \label{eq:lipschitz}
\end{equation}
In addition, if $t_\nu = z(\nu q/2n) \geq 1$, then
\begin{equation}
   \label{eq:dtnudnu}
   - \partial t_\nu / \partial \nu = \theta /(\nu t_\nu), \qquad
     \theta \in [ \hf, 1],
\end{equation}
and for $0 \leq r \leq 2$ and $t_\nu^2 > 2$,
\begin{equation}
   \label{eq:nutnur}
   \partial ( \nu t_\nu^r) / \partial \nu =
     t_\nu^{r-2} [ 1 - r \theta t_\nu^{-2}] > 0.
\end{equation}

(3) If $n^{-1} \log ^{5} n \leq \eta_{n}^{p} \leq b_2 n^{- b_3},$ then
\begin{equation}
   \label{eq:taubds}
   2 b_3 \log n - 2 \log b_2 \leq \tau_\eta^2 \leq 2 \log n - 10 \log
   \log n,
\end{equation}
and so for $n > n(b)$, we have
\begin{equation}
   \label{eq:taucrude}
   \tau_\eta^2 = 2 \gamma_n \log n \qquad \qquad
   0 < c(b) \leq \gamma_n \leq 1.
\end{equation}

(4) If $q_n \geq b_1 / \log n$ and $\nu \leq b_2 n^{1-b_3}$,
     then for $n \geq n(b_2,b_3)$,
     \begin{equation}
       \label{eq:tnubd}
       -3/2 \leq t_\nu - \sqrt{2 \log (n/\nu)} \leq 2 (b_1 b_3)^{-1/2}.
     \end{equation}

(5) (i) For $n > n(b)$,
\begin{equation}
   \label{eq:t1taueta}
   t_1 / \tau_\eta \leq c(b).
\end{equation}
\indent (ii) If $a \leq 1$ and $\eta_n^p \leq e^{-1/2}$, then
\begin{equation}
   \label{eq:takn}
   t[a k_n] \geq \tau_\eta - 3/2.
\end{equation}
If $a \leq \delta^{-1},$ the same inequality holds for $\eta <
\eta(p,\delta)$ sufficiently small.

(iii) If $a \geq \gamma \tau_\eta^{-1}$, then for $\eta_n^p \leq
\eta(\gamma,p,b_1,b_3)$ sufficiently small (and not the same at each
appearance),
\begin{align}
   t[a k_n] & \leq \tau_\eta + c(b_1,b_3),  \qquad \mbox{and}
     \label{eq:takn1} \\
   t[a k_n] & \leq 2 \tau_\eta.    \label{eq:takn2}
\end{align}

(iv) In particular, for $a \in [\gamma \tau_\eta^{-1}, \delta^{-1}]$,
then as $\eta_n \rightarrow 0$,
\begin{equation}
   \label{eq:takntilde}
   t[a k_n] \sim \tau_\eta.
\end{equation}
\end{lemma}

\begin{proof}
     (1)   Taking logarithms in the inequality
     \begin{displaymath}
       \frac{\phi(z(\eta))}{z(\eta)} \bigl( 1 - z(\eta)^{-2}) \bigr)
       \leq \tilde \Phi (z(\eta))  \leq \frac{\phi(z(\eta))}{z(\eta)},
     \end{displaymath}
and using $\tilde \Phi (z(\eta)) = \eta$ leads to
\begin{equation}
     \label{eq:sandwichz}
     2 \log ( 1 - z(\eta)^{-2}) \leq z^2(\eta) - 2 \log \eta^{-1} + 2
     \log z(\eta) + \log 2 \pi \leq 0.
\end{equation}
Now $ z^2 (\eta) \leq 2 \log \eta^{-1}$ for $\eta \leq 1/2$ and
$\eta \rightarrow \log( 1 - z(\eta)^{-2})$ is decreasing, so for
$\eta \leq \eta_0$,
\begin{displaymath}
     z^2(\eta) \geq 2 \log \eta^{-1} - \log \log \eta^{-1} +
     [2 \log (1 - z(\eta_0)^{-2} ) - \log 4 \pi ].
\end{displaymath}
For $\eta = 0.01,$ the quantity in square brackets equals -2.94.
Substituting this into the upper bound half of \eqref{eq:sandwichz}
yields
\begin{displaymath}
     z^2(\eta) \leq 2 \log \eta^{-1} - \log \log \eta^{-1} - g( \log
     \eta^{-1} ),
\end{displaymath}
where the increasing function
$g(v) = \log [ 2 \pi v^{-1} ( 2v - \log v - 3) ] \geq 1.8$
for $v \geq \log \eta_0^{-1} \doteq 4.60.$

Turning to the proof of (\ref{eq:zetabd}), note first that $\eta \leq
0.01$ entails $\log \log \eta^{-1} \geq 1.527$ which implies the right
inequality in view of the bound on $r(\eta)$ in
(\ref{eq:quantilebracket}).  For the left inequality, we rewrite
(\ref{eq:quantilebracket}) and appeal to (\ref{eq:logbound}), by setting
$k = \eta^{-1}$ and $\alpha^{-1} = \sqrt{\log \eta^{-1}}
e^{-r(\eta)/2}$.  The error bound in (\ref{eq:logbound}) is
$(\sqrt{8}/e) \sqrt{C}$, where
\begin{displaymath}
   C = \frac{1}{\alpha \log k} \leq \frac{e^{-r(\eta)/2}}{\sqrt{\log \eta^{-1}}}
\end{displaymath}
whenever $\eta \leq \eta_0$. With $\eta_0 = 0.01$, calculation shows
that $(\sqrt{8}/e) \sqrt{C} \leq 1.4850 \leq 3/2.$

(2) Differentiating the equation $\eta = \tilde{\Phi}(z(\eta))$
yields $z'(\eta) = -1 / \phi(z(\eta))$ which is decreasing in $\eta.$
Hence
\begin{displaymath}
       z(\eta_1) - z(\eta_2) \leq (\eta_2 - \eta_1) / \phi(z(\eta_1)).
\end{displaymath}
Since $\eta = \tilde{\Phi}(z(\eta)) \leq \phi(z(\eta))/ z(\eta),$
\eqref{eq:lipschitz} follows.

 From the Mills' ratio remark in the proof of Lemma
\ref{lem:gausstails}, we have $z \tilde \Phi (z) / \phi(z) \geq 1/2$
for $z \geq 1$, and so for such $z(\eta)$,
\begin{displaymath}
   -z'(\eta) = \theta / (\eta z(\eta) ), \qquad \theta \in [1/2,1].
\end{displaymath}
Since $\partial t_\nu / \partial \nu = (q_n/2n) z'(\nu q/ 2n)$, this
yields \eqref{eq:dtnudnu}, from which \eqref{eq:nutnur} is immediate.

%

(3,4) Displays \eqref{eq:taubds} and \eqref{eq:taucrude} are immedate.
  Bound (\ref{eq:zetabd}) says that $t_\nu - \sqrt{2 \log (2n/q \nu)}
     \in [-3/2,0].$ Apply inequality (\ref{eq:logbound}) with $\alpha =
     2/q \geq 1$ and $k = n / \nu$ to find that
     \begin{displaymath}
       0 \leq \sqrt{2 \log (2n/q \nu)} - \sqrt{2
         \log (n/\nu)} \leq c,
     \end{displaymath}
where for $n \geq n(b_2, b_3)$,
\begin{displaymath}
   c^2 = \frac{2}{q_n \log (n/\nu)} \leq \frac{2 \log n}{b_1 (b_3 \log n -
   b_2 )} \leq \frac{4}{b_1 b_3}.
\end{displaymath}

(5) Part (i) follows from \eqref{eq:tnubd} and
\eqref{eq:taucrude}. For part (ii) we first note, from
\eqref{eq:tnubd} that for $n > n(b_2, b_3)$,
\begin{displaymath}
   t[ak_n] = \sqrt{ 2 \log n /(a k_n)} + \theta_1, \qquad
   -3/2 \leq \theta_1 \leq 2 (b_1 b_3)^{-1/2}.
\end{displaymath}
Since $n/k_n = \eta_n^{-p} \tau_\eta^p$, we can write
\begin{displaymath}
   2 \log [n/(a k_n)] = \tau_\eta^2 + L(\eta,a), \qquad
   L(\eta,a) = p \log \tau_\eta^2 + 2 \log a^{-1}.
\end{displaymath}
If $\eta^p \leq e^{-1/2}$, then $\tau_\eta^2 = 2 \log \eta^{-p} \geq
1,$ and if also $a \leq 1$, then $L(\eta,a) \geq 0$, so that
\eqref{eq:takn} follows.
If however we assume only that $a \leq (1-q')^{-1}$, then $L(\eta,a)
\geq 0$ if $\tau_\eta^2$ is sufficiently large, i.e. if $\eta^p \leq
\eta(p,q')$ is sufficiently small.

For part (iii), from the bound $\sqrt{x+\epsilon} \leq \sqrt{x} +
\epsilon/(2 \sqrt{x})$ (valid for $0 \leq \epsilon \leq 3x$), we find
\begin{displaymath}
   \sqrt{\tau_\eta^2 + L(\eta,a)} \leq \tau_\eta + L(\eta,a)/(2 \tau_\eta),
\end{displaymath}
and clearly $L(\eta,a)/\tau_\eta \leq [ 2(p+1) \log \tau_\eta + 2 \log
\gamma^{-1} ]/\tau_\eta \leq 1$ for $\eta \leq \eta(\gamma,p)$.
This establishes \eqref{eq:takn1} and \eqref{eq:takn2} is a direct
consequence.
\end{proof}

\section{Proofs of Lower Bounds} \label{sec:lowerbds}
\setcounter{equation}{0}

This final appendix  combines ideas from Sections 6-8 
to finish the Lower Bound result.

\subsection{Proof of Proposition \ref{prop:hardthreshl0}}
\label{sec:proof-proposition-4.5}

 From the structure \eqref{eq:mualphal} of the configuration
$\mu_\alpha$, the total risk can be written in terms of the univariate
component risks as
\begin{displaymath}
   \rho(\hat\mu_{H,t}, \mu_\alpha) =
     k_n \rho_H(t, m_\alpha) + (n-k_n) \rho_H(t,0).
\end{displaymath}
 From \eqref{eq:hardriskat0} together with the definition of $t = t[a
k_n]$, we obtain
\begin{displaymath}
   n \rho_H(t[a k_n],0) = a k_n q_n t^r[a k_n](1 + \theta),
\end{displaymath}
with $0 \leq \theta \leq 2 t[a k_n]^{-2}.$ Since $t[a k_n] \sim \tau_\eta$ by
\eqref{eq:takntilde}, we conclude that
\begin{displaymath}
   (n-k_n) \rho_H( t[a k_n], 0) \sim a q_n k_n \tau_\eta^r.
\end{displaymath}

In the notation of \eqref{eq:dedecomp},
\begin{displaymath}
   \rho_H(t,m_\alpha) = m_\alpha^r [ \Phi(t-m_\alpha) -
   \Phi(-t-m_\alpha)] + E( m_\alpha,t),
\end{displaymath}
and as noted there, $0 \leq E(m_\alpha,t) \leq c_r$.
In addition, $m_\alpha^r \Phi(-t-m_\alpha) \leq c_r^\prime$, and as
noted at \eqref{eq:malpha}, $m_\alpha - t = \alpha +
O(\tau_\eta^{-1})$ and so $\Phi(t - m_\alpha) = \tilde \Phi(\alpha) +
O(\tau_\eta^{-1})$. Consequently, since $m_\alpha \sim \tau_\eta$, we
conclude that
\begin{displaymath}
   k_n \rho_H(t,m_\alpha) \sim k_n \tau_\eta^r \tilde \Phi(\alpha) (1 +
   o(1)).
\end{displaymath}

\subsection{Proof of Proposition \ref{prop:risk-hard-thresh}}
\label{sec:proof-proposition-4.6}

We use \eqref{eq:dedecomp} to decompose the total risk
\begin{displaymath}
   \rho( \hat \mu_{H,t}, \mu_\alpha)
    = \sum_l D(\mu_{\alpha l}, t) + E( \mu_{\alpha l}, t)
    = D + E,
\end{displaymath}
say. To bound the ``bias'' or ``false negative'' term $D$, we choose
an index $l_\alpha = n \eta_n^p m_\alpha^{-p}$ so that
$\bar \mu_{l_\alpha} = m_\alpha \sim \tau_\eta$.
Now decompose $D$ into $D_1 + D_2$ according as $l \leq l_\alpha$ or
$l > l_\alpha$. For $l \leq l_\alpha,$ we have identically
$\mu_{\alpha l} \equiv m_\alpha$, and so
\begin{displaymath}
   D_1 = l_\alpha m_\alpha^r [ \Phi(t - m_\alpha) - \Phi( -t -
   m_\alpha)] = l_\alpha m_\alpha^r \tilde \Phi(\alpha) (1 + o(1))
\end{displaymath}
by the same arguments as for Proposition \ref{prop:hardthreshl0}.
And since $m_\alpha \sim \tau_\eta,$ we have
\begin{displaymath}
   l_\alpha m_\alpha^r = n \eta_n^p m_\alpha^{r-p}
         \sim n \eta_n^p \tau_\eta^{r-p}
        = k_n \tau_\eta^r.
\end{displaymath}
The novelty with the weak-$\ell_p$ risk comes in the analysis of
\begin{displaymath}
   D_2 = \sum_{l > l_\alpha} \bar \mu_l^r [ \Phi(t - \bar \mu_l) - \Phi(
   -t - \bar \mu_l)].
\end{displaymath}
The second term is negligible, being bounded in absolute value by
$\tilde \Phi(t) \sum_{l > l_\alpha} \bar \mu_l^r = o( k_n
\tau_\eta^r)$.
For the first term we have an integral approximation (with $\bar \mu
(x) = \eta_n (n/x)^{1/p}$)
\begin{displaymath}
   D_2 \sim \int_{l_\alpha}^n \bar \mu^r(x) \Phi( t - \bar \mu (x)) dx
     = p l_\alpha m_\alpha^r \int_{\eta_n/m_\alpha}^1 v^{r-p-1} \Phi(
     t- m_\alpha v) dv,
\end{displaymath}
after setting $x = l_\alpha v^{-p}$.
[Remark: to bound the error in the integral approximation,
observe that if $f'(x)$ is smooth with at most one zero in $[a,b]$, then
the difference between $\sum_a^b f(l)$ and  $\int_a^b f$
is bounded by $\sup_{[a,b]} |f|.$]

For $0 < v < 1$, we have $t - m_\alpha v \sim (1-v) \tau_\eta
\rightarrow \infty$, and so from the dominated convergence theorem,
the integral converges to $\int_0^1 v^{r-p-1} dv$, so
\begin{displaymath}
   D_2 \sim [p/(r-p)] l_\alpha m_\alpha^r \sim [p/(r-p)] k_n
   \tau_\eta^r.
\end{displaymath}
Putting together the analyses of $D_1$ and $D_2$, we find that the
false-negative term
\begin{displaymath}
   D = [ \tilde \Phi(\alpha) + [p/(r-p)] ] k_n \tau_\eta^r (1 + o(1)).
\end{displaymath}

To bound the ``variance'' or ``false positive'' term $E$, decompose
the sum into three terms, corresponding to indices in the ranges
$[1, l_1], (l_1, l_2]$ and $(l_2,n]$, where
\begin{displaymath}
   l_1 = n \eta_n^p \ \  \leftrightarrow \ \ \bar \mu_{l_1} = 1, \qquad
   \text{and} \qquad
l_2 = n \eta_n^p t_\eta^{2p} \ \  \leftrightarrow \ \ \bar \mu_{l_2} =
   \tau_\eta^{-2}.
\end{displaymath}
We use (i) - (iii) of Lemma \ref{lem:threshdecomp} to show that terms
$E_1$ and $E_2$ are negligible. For $E_1$, the global bound (i)
gives $E_1 \leq c_r l_1 = o( k_n \tau_\eta^r)$.
For $E_2$, properties (ii) (monotonicity) and (iii) show that
\begin{displaymath}
   E(\bar \mu_l,t) \leq E(\bar \mu_{l_1},t) = E(1,t)
   \leq c_0 t^r \tilde \Phi(t-1)
   \leq c_1 \tau_\eta^r \tilde \Phi( \tau_\eta - 5/2),
\end{displaymath}
where the last inequality uses
$t = t[a k_n] \sim \tau_\eta$ and
that $t \geq \tau_\eta - 3/2$ from \eqref{eq:takn}.
Hence
\begin{displaymath}
   E_2 \leq l_2 E(1,t) \leq c_0 k_n \tau_\eta^r \cdot \tau_\eta^{3p}
   \tilde \Phi( \tau_\eta - 5/2) = o( k_n \tau_\eta^r).
\end{displaymath}

Finally, we focus on the dominant term $E_3 =\sum_{l > l_2} E( \bar
\mu_l, t)$.
For $l > l_2$, we have $\bar \mu_l \leq \tau_\eta^{-2}$ and $t \leq c
\tau_\eta$ so that
\begin{displaymath}
   \phi(t \pm \bar \mu_l)
    = \phi(t) \exp \{ \mp t \bar \mu_l - \bar \mu_l^2 /2 \}
    = \phi(t) ( 1 + O(\tau_\eta^{-1}) ).
\end{displaymath}
Now if $|\mu| \leq \tau_\eta^{-1}$ and $t = t[a k_n] \leq 2 \tau_\eta$
(from \eqref{eq:takn2}), then
\begin{displaymath}
   \phi(t-\mu) = \phi(t) \exp (t \mu - \mu^2/2 )
              = \phi(t) ( 1 + O(\tau_\eta^{-1}) ),
\end{displaymath}
and so
\begin{align*}
   \gamma(t-\mu) & = |t-\mu|^{r-1} \phi(t-\mu) =
        t^{r-1} \phi(t) (1 + O(\tau_\eta^{-1}) ), \ \  \text{and} \\
        \epsilon(t-\mu) & \leq |t-\mu|^{r-3} \phi(t-\mu) \leq
        t^{r-3} \phi(t) (1 + O(\tau_\eta^{-1}) ).
\end{align*}
Consequently, using \eqref{eq:emut} and $t^{-1} \phi(t) \sim \tilde
\Phi(t)$, we get
\begin{displaymath}
   E_3 = 2(n-l_2) t^{r-1} \phi(t) ( 1 + O(\tau_\eta^{-1}) )
       = a q_n k_n \tau_\eta^r (1+ o(1)),
\end{displaymath}
using the manipulations of \eqref{eq:zerocptrisk} -
\eqref{eq:useasymp}.

\subsection{Proof of Proposition \ref{prop:kmualpha}.}
\label{sec:proof-proposition-4.7}

Fix $a_1 > 0.$
The idea, both for $\ell_0[\eta]$ and for $m_p[\eta_n]$, is to obtain
bounds
\begin{displaymath}
   M_-(k) \leq M(k; \mu_\alpha) \leq M_+(k),
   \qquad \qquad k \in [ a_1 k_n,  a_1^{-1} k_n],
\end{displaymath}
for which solutions $k_\pm$ to $M_\pm (k) = k$ can be easily found.
 From monotonicity of $k \rightarrow M(k;\mu_\alpha)$, it then follows
that $k_- \leq k(\mu_\alpha) \leq k_+$.

In both cases we establish a representation of the form
\begin{equation}
   \label{eq:maknrep}
   M(ak_n;\mu_\alpha) = k_n [\Phi(\alpha) + \theta_1 \tau_\eta^{-1}]
                        + [1 + \theta_2 \delta(\eta_n)] q_n k,
\end{equation}
where $|\theta_i| \leq c(\alpha,a,b)$ and $\delta(\eta_n) \rightarrow
0$ as $\eta_n \rightarrow 0$.
 From this, expressions for $k_\pm$ are easily found and
$k_\pm \sim k_n \Phi(\alpha)/(1-q_n)$ is easily checked.

For $\ell_0[\eta_n]$, \eqref{eq:maknrep} follows from
\eqref{eq:mkmualpha} and \eqref{eq:malpha}.
For $m_p[\eta_n]$, we first formulate a lemma.

\begin{lemma}
   \label{lem:ave-behavior}
Let $a$ and $\alpha$ be fixed. If $|f|$ and $|f'|$ are bounded by 1, then
\begin{displaymath}
   k_n^{-1} \sum_{l=1}^{k_n} f( \mu_{\alpha l} - t[a k_n]) =
      f(\alpha) + \theta \tau_\eta^{-1},
      \qquad |\theta| \leq c(\alpha,a,b).
\end{displaymath}
\end{lemma}
\begin{proof}
   For $l \leq l_\alpha = n \eta_n^p m_\alpha^{-p}$, we have
$\mu_{\alpha l} = m_\alpha$ and so with $t = t[a k_n]$,
\begin{displaymath}
   LHS = k_n^{-1} \sum_1^{k_n \wedge l_\alpha} f( m_\alpha - t)
          + \theta_1 \| f \|_\infty [ 1 - 1 \wedge (l_\alpha/k_n)],
         \qquad |\theta_1| \leq 1.
\end{displaymath}
 From \eqref{eq:takkbd}, we have
\begin{displaymath}
   f(m_\alpha -t)
     = f(\alpha) + \theta_2 \| f' \|_\infty (t[k_n] - t[a k_n])
     = f(\alpha) + \theta_3 \tau_\eta^{-1},
\end{displaymath}
with $|\theta_2| \leq 1$ and $|\theta_3| \leq c(a)$.
 From \eqref{eq:takn} and \eqref{eq:takn1} we have
$|m_\alpha - \tau_\eta| \leq c(b)$ and so
\begin{displaymath}
   k_n/l_\alpha = (m_\alpha/\tau_\eta)^p = 1 + \theta_5 \tau_\eta^{-1},
    \qquad \quad |\theta_5| \leq c(b).
\end{displaymath}
Now combine the last three displays.
\end{proof}

We exploit the division into the `positive',
`transition' and `negative' zones defined in Section
\ref{sec:lipschitz-result}.
Applying the preceding Lemma with $f = \Phi$ yields
\begin{displaymath}
   M_{pos}(k; \mu_\alpha) = k_n[ \Phi(\alpha) + \theta_3 \tau_\eta^{-1}],
\end{displaymath}
while from \eqref{eq:trans-m} and \eqref{eq:mnbd} we obtain
\begin{displaymath}
   M_{trn}(k;\mu_\alpha) = \theta_4 k_n \tau_\eta^{-1}, \quad \text{and}
   \quad M_{neg}(k;\mu_\alpha) = [ 1 + \theta_5 \delta_p (\epsilon_n)] q_n k.
\end{displaymath}
Putting together the last two displays, we recover \eqref{eq:maknrep}
and the result.

\subsection{Proof of Proposition \ref{prop:riskmualpha}.}
\label{sec:proof-proposition-4.8}


Let $t'$ denote the fixed threshold $t' = t[a_0 k_n]$
where $a_0 = (1-q)^{-1} \Phi(\alpha)$.
We will use Lemma \ref{lem:riskdel} to show that
\begin{displaymath}
   | \rho( \hat \mu_F, \mu_\alpha) - \rho( \hat \mu_{H,t'},
     \mu_\alpha)| = o(k_n \tau_\eta^r )
\end{displaymath}
so that the conclusion will follow from Proposition
\ref{prop:risk-hard-thresh}.

To apply Lemma \ref{lem:riskdel}, set $\hat \mu = \hat \mu_F$ and
$\hat \mu' = \hat \mu_{H,t'}$, so that the thresholds
$\hat t = \hat t_F$ and $\hat t' = t'$ respectively, which are both
bounded by $t_1$.

Let $k_\pm = k_\pm(\mu_\alpha) = k(\mu_\alpha) \pm \alpha_n k_n$ and
recall that the event $A_n = \{ k_- \leq \hat k_F \leq k_+ \}$ is
$m_p[\eta_n]-$likely.
If we set $t_\pm = t[k_\mp]$, then on event $A_n$,
\begin{displaymath}
   N^\prime = \# \{ i: |y_i| \in [ \hat t_F, t^\prime] \}
      \leq N^{\prime \prime} = \# \{ i: |y_i| \in [ t_-, t_+] \}.
\end{displaymath}
Hence $P(N^\prime > \beta_n) \leq P(A_n^c) + P(N^{\prime \prime} >
\beta_n )$.

We use the exponential bound of Lemma \ref{lem:bennetcor} to choose
$\beta_n$ so that $P(N^{\prime \prime} > \beta_n )$ is small.
 From the definition of the threshold function, $M(k) = M(k;
\mu_\alpha)$, we have
\begin{displaymath}
   M^{\prime \prime} := E N^{\prime \prime}
       = M(k_+) - M(k_-).
\end{displaymath}
Since $k \rightarrow \dot M_k$ is decreasing and using the derivative
bounds of Proposition \ref{prop:lipschitz}, we find that for $n$
sufficiently large,
\begin{equation}
   \label{eq:mppbds}
   (k_+ - k_-) \dot M(k_+) \leq M^{\prime \prime}
   \leq (k_+ - k_-) \dot M(k_-).
\end{equation}

\noindent
\textit{Claim.} $\dot M (k_\pm; \mu_\alpha) = \theta [q_n + c(\alpha)
\tau_\eta^{-1}] \quad $ for some $\theta \in [1/2,2]$.
\begin{proof}
   We again use the positive-transition-negative decomposition, this
   time of $\dot M(a k_n; \mu_\alpha)$.
Write, with $\nu = a k_n$,
\begin{displaymath}
   \dot M_{pos}( a k_n; \mu_\alpha)
    = (- \partial t_\nu/ \partial \nu) \sum_1^{k_n}
            \phi(t - \mu_{\alpha l}) + \phi( t + \mu_{\alpha l}).
\end{displaymath}
 From \eqref{eq:lipschitz} and \eqref{eq:dtnudnu}, we have
$- k_n ( \partial t_\nu/ \partial \nu) \sim 1/(a \tau_\eta).$
Applying Lemma \ref{lem:ave-behavior} to $f(x) = \phi(-x)$, we
conclude that for $n$ sufficiently large,
\begin{displaymath}
   \dot M_{pos} (a k_n; \mu_\alpha) = (\theta_1/(a \tau_\eta))
   [\phi(\alpha) + \theta_2 \tau_\eta^{-1}].
\end{displaymath}
Appealing to \eqref{eq:trans-mdot} and \eqref{eq:mndotbd},
\begin{displaymath}
   \dot M_{trn}( a k_n; \mu_\alpha) = \theta_3 \tau_\eta^{-2}, \qquad
      \text{and} \qquad
   \dot M_{neg}( a k_n; \mu_\alpha) = [1 + \theta_4 \delta_p(\epsilon_n)] q_n.
\end{displaymath}
Combining the last two displays and noting that $k_\pm = k(\mu_\alpha)
\pm \alpha_n k_n$ correspond to $a = \phi(\alpha)(1 -q_n)^{-1}$, we
obtain the claim.
\end{proof}

\medskip
Now set $q_{\alpha n} = q_n + c(\alpha) \tau_\eta^{-1}$ and select
$\beta_n = 8 \alpha_n k_n q_{\alpha n}$.
 From \eqref{eq:mppbds} and the claim, we have
\begin{displaymath}
   \alpha_n k_n q_{\alpha n} \leq M^{\prime \prime}
      \leq 4 \alpha_n k_n q_{\alpha n},
\end{displaymath}
and so $\beta_n/  M^{\prime \prime} \geq 2.$ Consequently,
\begin{align*}
   P( N^{\prime \prime} > \beta_n)
     & \leq \exp \{ - (1/4) M^{\prime \prime} h(\beta_n/ M^{\prime
     \prime}) \}  \\
     & \leq \exp \{ - (1/4) \alpha_n k_n q_{\alpha n} h(2) \}
       \leq c_0 \exp \{ - c_1 \log^2 n \}
\end{align*}

Now $\beta_n t_1^r \asymp \alpha_n q_{\alpha n} k_n \tau_\eta^r
  = o( k_n \tau_\eta^r)$, while Proposition \ref{prop:risk-hard-thresh}
  shows that $\rho( \hat \mu_{H,t^\prime}, \mu_\alpha) = O( k_n
  \tau_\eta^r).$ So Lemma \ref{lem:riskdel} applies and we are done.



\end{document}